\documentclass[11pt,a4paper]{article}

\usepackage[english]{babel}
\usepackage{amsmath}
\usepackage{amssymb}
\usepackage{graphicx}        

\usepackage{cases}
\usepackage{amsfonts}
\usepackage{color}
\usepackage{verbatim}
\usepackage{caption}
\usepackage{subcaption}
\usepackage{fullpage}

\usepackage{array}

\usepackage{epsfig}
\usepackage{float}
\usepackage{url}

\usepackage{theorem}

\usepackage{fullpage}

\usepackage{hyperref}

\usepackage{lmodern} 
\usepackage{amsmath}
\usepackage{amssymb}
\usepackage{graphicx,scalerel}

\newcommand\nwhat[1]{\hstretch{2}{\hat{\hstretch{.5}{#1}}}}

\def\overbigdot#1{\overset{\hbox{\tiny$\bullet$}}{#1}}
\def\twobigdot#1{\overset{\hbox{\tiny$\bullet \bullet$}}{#1}}

\definecolor{viola}{rgb}{0.3,0,0.7}
\definecolor{ciclamino}{rgb}{0.5,0,0.5}
\definecolor{rosso}{rgb}{0.85,0,0}

\def\pier #1{{\color{black} #1}}
\def\abramo #1{{\color{black} #1}}
\def\abramonew #1{{\color{black} #1}}
\def\pco #1{{\color{black} #1}}
\def\abramohh #1{{\color{black} #1}}
\def\abramohhb #1{{\color{black} #1}}
\def\pierhhb #1{{\color{black} #1}}
\def\michhhb #1{{\color{black} #1}}
\def\ap #1{{\color{black} #1}}
\def\reviewa #1{{\color{black} #1}}
\def\reviewm #1{{\color{black} #1}}

\newtheorem{rem}{Remark}[section]
\newtheorem{lem}{Lemma}[section]
\newtheorem{thm}{Theorem}[section]

\newenvironment{pf}{\noindent{\bf Proof. \/}\noindent%
}{\hfill\EndProofMarker}
\newcommand{\EndProofMarker}{$\Box$}
\newcommand{\diver}{\ensuremath{\operatorname{div}}}
\newcommand{\curl}{\ensuremath{\operatorname{curl}}}

\newcommand{\grad}{\ensuremath{\operatorname{grad}}}

\begin{document}
\title{ Large deformations  in terms of stretch and rotation\\ and global solution to the quasi-stationary problem}
\iftrue
\date{}
\maketitle 
\begin{center}
\vskip-1.2cm
{\large\sc Abramo Agosti$^{(1)}$}\\
{\normalsize e-mail: {\tt abramo.agosti@unipv.it}}\\[0.25cm]
{\large\sc Pierluigi Colli $^{(1),(2)}$}\\
{\normalsize e-mail: {\tt pierluigi.colli@unipv.it}}\\[0.25cm]
{\large\sc Michel Fr\'emond$^{(3)}$}\\
{\normalsize e-mail: {\tt michel.fremond@uniroma2.it}}\\[.5cm]
$^{(1)}$
{\small Department of Mathematics, University of Pavia}\\
{\small via Ferrata 5, I-27100 Pavia, Italy}\\[.3cm] 
$^{(2)}$
{\small Research Associate at the IMATI--C.N.R, Pavia}\\
{\small via Ferrata 5, I-27100 Pavia, Italy}\\[.3cm] 
$^{(3)}$
{\small Lagrange Laboratory,}\\
{\small \michhhb{Department of Civil Engineering and Informatics Engineering}, \\University of Roma ``Tor Vergata"}\\
{\small Via del Politecnico, 1, 00163 Roma, Italy}
\end{center}
\fi
\maketitle

\date{}

\begin{abstract}
In this paper we derive a new model for \michhhb{visco-elasticity} with large deformations where the independent variables are the stretch and the rotation tensors \michhhb{which intervene with second gradients terms accounting for physical properties  in} the principle of virtual power.
\michhhb{Another} basic feature of our model is that there is conditional compatibility, entering the model as \pierhhb{kinematic} constraints and depending on the magnitude of an internal \ap{force} associated to dislocations. Moreover, due to the \pierhhb{kinematic} constraints, the virtual velocities depend on the solutions of the problem. As a consequence, the variational formulation of the problem and the related mathematical analysis are neither standard nor straightforward. We adopt the strategy to invert the \pierhhb{kinematic} constraints through Green propagators, obtaining a system of integro-differential coupled equations.
As a first mathematical step, we develop the analysis of the model in a simplified setting, i.e. considering the quasi-stationary version of the full system where we neglect inertia. In this context, we prove the existence of a global in time strong solution in three space dimensions for the system, employing techniques from PDEs and convex analysis, thus obtaining a novel \reviewa{breakthrough} in the field of three-dimensional finite \michhhb{visco-elasticity} described in terms of the stretch and rotation variables.
We also study a limit problem, letting the magnitude of the internal \ap{force} associated to dislocations tend to zero, in which case the deformation becomes incompatible and the equations takes the form of a coupled system of PDEs. For the limit problem we obtain global existence, uniqueness and continuous dependence from data in three space dimensions. 
\end{abstract}


\section{Introduction}
Large deformation theory introduced by John Ball relies on the gradient matrix%
\begin{equation}
 \label{eqn:intro1}
\mathbf{F}=\grad \boldsymbol{\Phi} =\mathbf{RW}, 
\end{equation}
with position function $\boldsymbol{\Phi}$, stretch matrix $\mathbf{W}$\ and rotation
matrix $\mathbf{R}$.
\reviewa{For what concerns the decomposition of the gradient of position in terms of the stretch and the rotation tensors, we point out that the analysis of models in nonlinear three-dimensional {visco-elasticity} where the rotation field is considered as one of the primary unknowns was posed in \cite{ciarlet} as an open problem. {In the latter reference, the authors highlighted} the interest in this approach to describe elasticity since it posseses ``a more geometrical flavour than the classical approach'' with the deformation gradient.}
\reviewa{In the present work we choose to describe the motion with matrices $\mathbf{R}$ and $\mathbf{W}$
which can be experimented and we aim at obtaining existence and regularity results for the derived model}.

To satisfy relationship \eqref{eqn:intro1}, stretch matrix $\mathbf{W}$ has to
satisfy the compatibility conditions. These conditions can be infringed in
case there are dislocations, \abramohh{in which case the matrix $\mathbf{RW}$ is no longer a gradient \cite[\pierhhb{p.~389}]{germain}}, \ap{see also} \michhhb{\cite{acharya}}. In this situation we have \pierhhb{\begin{gather*}
\mathbf{RW=}\grad \boldsymbol{\Phi} \mathbf{+}\curl \mathbf{Z,} \quad 
\diver \mathbf{Z}=0,
\end{gather*}}%
where $\mathbf{Z}$ accounts for the dislocations. The occurrence of
dislocations is conditional depending on the intensity of an internal force.

We prove an existence theorem in the framework of these ideas. \abramohh{We start by deriving a model, \ap{which \reviewa{is further}} \michhhb{detailed in \cite{fremond}}, with conditional compatibility and inertia\ap{. We start} from a generalized form of the principle of virtual power, together with kinematic conditions for the deformation map and the dislocations, choosing constitutive assumptions for the  \michhhb{internal forces} in the system in order \ap{to satisfy} the Clausius--Duhem dissipative equality. 
\ap{We assume a quadratic expression for the free energy density of the system, depending on the stretch, the rotation and the dislocation tensors, and a quadratic form also for the dissipation potential, containing viscous contributions in terms of the time derivative of the stretch tensor and on the angular velocity tensor.}}

\abramohhb{The generalized virtual velocities are associated to the main variables of the problem, i.e.\pierhhb{,} to $\boldsymbol{\Phi}, \mathbf{W,R,Z}$, which are linked by \pierhhb{kinematic} constraints. As a consequence, the virtual velocities themselves satisfy internal constraints depending on the solutions of the problem. This feature of the model makes the definition of weak solutions involved, since we should deal with a variational formulation with test functions depending on the solutions themselves. Hence, we adopt the strategy to express the virtual velocities associated to the deformation map and the dislocations in terms of the virtual velocities associated to the stretch matrix and to the rotation through integral operators \ap{related} to the \pierhhb{kinematic} constraints\pierhhb{. Thus, we reduce the set of independent virtual velocities and eliminate} their internal constraints, obtaining a system of integro-differential coupled equations.}

\abramohhb{We consider the inertia of the system expressed by a virtual power of acceleration forces containing second-order interaction terms in space, which allows us to obtain sufficient regularity of weak solutions to be able to take into account a 
point with inertia, in agreement with experiments.}
\reviewa{Indeed, consider a nail firmly hammered in a wall and assume it is a point which has
inertia. When experimenting the motion of the nail, its angular velocity as
well as its angular acceleration are equal to the angular velocity and
acceleration of the wall: these quantities are continuous with respect to
space. In mathematical parlance this mechanical property is: the angular
velocity and acceleration of the nail are the trace on the point of the
volume quantities. Thus it is reasonable {that} they are functions with such
mathematical properties: for instance, with  second
order derivatives {that are volume square integrable}. Note that in this point of view the bilateral contact of
the wall with a wooden stick, a beam involving second order derivatives\pierhhb{,} is
straightforward.}

{We also associate to dislocations an internal \ap{force}, which is a new independent variable of the system\reviewa{, a stress, which activates dislocations when its magnitude is greater than a certain threshold $k>0$}. \ap{Also}, we impose the positive definiteness of the stretch matrix as an internal constraint in the free energy of the system, \abramohhb{which implies that the material is neither flattening nor crushing and that a point which is inside its domain at a certain time remains in the interior of the domain at later times \cite[Theorem 5.5.1]{ciarlet2}.}}

\pier{In the present contribution, as a first step and in order to expose all the technicalities for a simplified problem, we develop the analysis for the quasi-stationary approximation of the full system, i.e., neglecting inertia. In this context, we obtain the existence of a global in time strong solution in three space dimensions. Together with the derivation of the model, \abramohhb{to the best of our knowledge this analytical result is a novel contribution in the framework of \pierhhb{3D} finite visco-elastic problems solved in the stretch matrix and rotation variables. Hence, it is a first step in the analysis of models in nonlinear three-dimensional \michhhb{visco-elasticity} where the rotation field is considered as one of the primary unknowns, which, \reviewa{as already observed}, was firstly posed as an open problem in \cite{ciarlet}.}}

{{We point out that the available analytical results in two and three spatial dimensions for visco-elastic problems with large deformations described in the standard approach with the deformation gradients entail the existence of local in time weak solutions \cite{bonetti1,bonetti2,bonetti3}.} We also study the limit problem as $k\to 0$, in which case it becomes a coupled system of PDEs where the incompatibility is always active. In this latter situation we obtain global existence, uniqueness and continuous dependence from data (\pierhhb{i.e., well-posedeness}) in three space dimensions. The study of the full case with \pierhhb{inertial terms} will be the subject of a second contribution.}

\pier{The paper is organized as follows. In Section \ref{sec:notation} we introduce the necessary notation and some preliminary results. In Section \ref{sec:model} we derive the full model with inertia and  \michhhb{the new internal forces terms}. In Section \ref{sec:analysis} we study the existence problem for the quasi-stationary approximation of the full problem. In Section \ref{sec:limit} we complete the analysis by studying the limiting case as $k\to 0$. We conclude with some observations and future perspectives in Section \ref{sec:conclusions}.}

\section{Notations and preliminaries}
\label{sec:notation}
In this section we introduce the notation and the preliminary results about the functional setting which will be necessary for the model derivation.

\subsection{Geometrical and functional setting}
Let $\mathcal{D}_a\subset \mathbb{R}^3$ be an open bounded and simply connected domain with Lipschitz boundary $\Gamma_a:=\partial \mathcal{D}_a$, and let $[0,T]$ be a finite time interval, with $T>0$. We introduce the notation $\mathcal{D}_{aT}:= \mathcal{D}_a\times [0,T]$. \reviewa{In the following, we use the bold notation to indicate quantities which are not scalars, i.e. vectors and tensors}. We indicate as $M(\mathbb{R}^{3\times 3})$ the linear space of square matrices, endowed with the Frobenius inner product
\[
\mathbf{A}: \mathbf{B}=\sum_{i,j=1}^3\mathbf{A}_{ij}\mathbf{B}_{ij},
\]
for any $\mathbf{A},\mathbf{B}\in M(\mathbb{R}^{3\times 3})$.
We also indicate with the notation $: :$ the Frobenius inner product in $M(\mathbb{R}^{3\times 3\times 3})$, and with the notation $: : :$ the Frobenius inner product in $M(\mathbb{R}^{3\times 3\times 3\times 3})$. \pier{The orthogonal subspaces of symmetric and antisymmetric matrices are denoted by} \pier{$Sym(\mathbb{R}^{3\times 3})\subset M(\mathbb{R}^{3\times 3})$ and $Skew(\mathbb{R}^{3\times 3})\subset M(\mathbb{R}^{3\times 3})$}, 
respectively. We indicate the set of special orthogonal matrices as $SO(\mathbb{R}^{3\times 3})$ and the set of positive definite symmetric matrices as  $Sym^+(\mathbb{R}^{3\times 3})$. \abramo{We recall that for any $\mathbf{R}\in SO(\mathbb{R}^{3\times 3})$ there exists a unique $\mathbf{A}\in Skew(\mathbb{R}^{3\times 3})$ such that $\mathbf{R}=\mathrm{e}^{\mathbf{A}}$, where the exponential of a matrix must be intended as $\mathrm{e}^{\mathbf{A}}=\sum_{n=0}^{\infty}\frac{\mathbf{A}^n}{n!}$.}
\pier{For a generic subset $K \subset M(\mathbb{R}^{3\times 3})$, let $I_K : M(\mathbb{R}^{3\times 3}) \to \{0, +\infty\}$  denote the indicator function of $K$, which is defined, for any $\mathbf{A} \in M(\mathbb{R}^{3\times 3})$,by  $ I_k(\mathbf{A})=0 $ if $\mathbf{A} \in K$, $ I_k(\mathbf{A})= +\infty $ if $\mathbf{A} \not\in K$.}

We introduce the space of vector fields $\mathcal{V}:=(\mathbb{R}^3)^{\mathcal{D}_{aT}}$, whose elements are functions from $\mathcal{D}_{aT}$ to $\mathbb{R}^3$. We further introduce the spaces of tensor fields $\mathcal{M}:=(M(\mathbb{R}^{3\times 3}))^{\mathcal{D}_{aT}}$, $\mathcal{SO}:=(SO(\mathbb{R}^{3\times 3}))^{\mathcal{D}_{aT}}$, $\mathcal{S}:=(Sym(\mathbb{R}^{3\times 3}))^{\mathcal{D}_{aT}}$ and $\mathcal{A}:=(Skew(\mathbb{R}^{3\times 3}))^{\mathcal{D}_{aT}}$, with $\mathcal{M}=\mathcal{S}\oplus \mathcal{A}$. Given a tensor $\mathbf{A}\in \mathcal{M}$, we denote by $\rm{Sym}(\mathbf{A}):=\frac{\mathbf{A}+\mathbf{A}^T}{2}$ its symmetric part and by $\rm{Skew}(\mathbf{A}):=\frac{\mathbf{A}-\mathbf{A}^T}{2}$ its antisymmetric part. We also need to introduce the space of tensor fields $\mathcal{M}_{\diver}:=\{\mathbf{A}\in \mathcal{M}|\diver \mathbf{A}=\mathbf{0}\}$, where the \pco{\textit{divergence}} of a second order tensor is defined row wise. In the following, we will operate also with the \textit{curl} of second order tensors, which is defined row wise. 

We denote by $L^p(\mathcal{D}_a;K)$ and $W^{r,p}(\mathcal{D}_a;K)$ the standard Lebesgue and Sobolev spaces of functions defined on $\mathcal{D}_a$ with values in a set $K$, {where $K$ may be $\mathbb{R}$ or a vector subspace of a multiple power of $\mathbb{R}$,} and by  $L^p(0,t;V)$ the Bochner space of functions defined on $(0,t)$ with values in the functional space $V$, with $1\leq p \leq \infty$. If $K\equiv \mathbb{R}$, we simply write  $L^p(\mathcal{D}_a)$ and $W^{r,p}(\mathcal{D}_a)$. 
For a normed space $X$, the associated norm is denoted by $\pier{\|}\cdot\pier{\|}_X$. In the case $p=2$, we use the notations $H^1:=W^{1,2}$ and $H^2:=W^{2,2}$, and we denote by $(\cdot,\cdot)$ and $\pier{\|}\cdot\pier{\|}$ the $L^2$ scalar product and induced norm between functions with scalar, vectorial or tensorial values. Moreover, we denote by $C^k(\mathcal{D}_a;K),C_c^k(\mathcal{D}_a;K)$ the spaces of continuously differentiable functions (respectively with compact support) up to order $k$ defined on $\mathcal{D}_a$ with values in a set $K$\pco{; by}  
{$C^k([0,t];V)$, $k\geq 0$, the spaces of {continuously differentiable functions up to} order $k$ from $[0,t]$ to the space $V$. The dual space of a Banach space $Y$ is denoted by $Y'$. Finally, we denote by $W_0^{r,p}
(\mathcal{D}_a;K)$ the closure of $C_c^{\infty}(\mathcal{D}_a;K)$ with respect to the norm $\pier{\|}\cdot\pier{\|}_{W^{r,p}(\mathcal{D}_a;K)}$, and \pco{by} $W^{-r,p'}(\mathcal{D}_a;K)$ the dual space of $W_0^{r,p}(\mathcal{D}_a;K)$, with $p\geq 1$ and $p'\geq 1$ conjugate exponents. As before, when $p=2$ we will indicate 
the latter functional spaces as $H_0^r(\mathcal{D}_a;K)$ and $H^{-r}(\mathcal{D}_a;K)$. The duality pairing between $H_0^1(\mathcal{D}_a;K)$ and $H^{-1}(\mathcal{D}_a;K)$ is denoted by $<\cdot,\cdot>$. We endow the space 
$H_0^1(\mathcal{D}_a;K)$ with the inner product $(A,B)_{H_0^1(\mathcal{D}_a;K)}:=(\grad A,\grad B)$, for all $A,B \in H_0^1(\mathcal{D}_a;K)$, and we introduce the Riesz isomorphism $\mathcal{R}:H_0^1(\mathcal{D}_a;K)\rightarrow 
H^{-1}(\mathcal{D}_a;K)$, defined by
\[
<\mathcal{R}A,B>=(A,B)_{H_0^1(\mathcal{D}_a;K)}, \quad \forall A,B\in H_0^1(\mathcal{D}_a;K).
\]
The operator $\mathcal{R}=-\Delta$ is the negative weak Laplace operator with homogeneous Dirichlet boundary conditions, which is positive definite and self adjoint. As a consequence of the Lax--Milgram theorem and the Poincar\'e inequality, the inverse operator $(-\Delta)^{-1}:H^{-1}(\mathcal{D}_a;K)\rightarrow H_0^1(\mathcal{D}_a;K)$ is well defined, and we set $A:=(-\Delta)^{-1}F=\mathcal{G}_L\ast F$, for $F\in H^{-1}(\mathcal{D}_a;K)$, where $\mathcal{G}_L$ is the Green propagator associated to the Laplace operator with 
homogeneous Dirichlet boundary conditions and $\ast$ denotes the convolution operation, if $-\Delta A=F$ in $\mathcal{D}_a$ in the weak sense, and $A=0$ on $\Gamma_a$ in the sense of traces. We note that, if $A\in H_0^1(\mathcal{D}_a;K)$ solves $-\Delta A=F$ for some $F\in W^{m,p}(\mathcal{D}_a;K)$, $1<p<\infty$, $m\in \mathbb{N}$, and $\Gamma_a$ is of class $C^{m+2}$, then from elliptic regularity theory $A\in W^{m+2,p}(\mathcal{D}_a;K)$ and $-\Delta A=F$ a.e. in $\mathcal{D}_a$, with 
\begin{equation}
\label{eqn:1}
\pier{\|}A\pier{\|}_{W^{m+2,p}(\mathcal{D}_a;K)}\leq C\pier{\|}F\pier{\|}_{W^{m,p}(\mathcal{D}_a;K)}.
\end{equation}
 We also need to introduce the spaces
\begin{align*}
& L_{\diver}^2(\mathcal{D}_a,K):=\overline{\{\mathbf{u}\in C_c^{\infty}(\mathcal{D}_a,K): \, \diver\mathbf{u}=0 \; \text{in} \; \mathcal{D}_a\}}^{\pier{\|}\cdot\pier{\|}_{L^2(\mathcal{D}_a;K)}},\\
& H_{0,\diver}^1(\mathcal{D}_a,K):=\overline{\{\mathbf{u}\in C_c^{\infty}(\mathcal{D}_a,K): \, \diver\mathbf{u}=0 \; \text{in} \; \mathcal{D}_a\}}^{\pier{\|}\cdot\pier{\|}_{H^1(\mathcal{D}_a;K)}}.
\end{align*}
\noindent
The duality pairing between $H_{0,\diver}^1(\mathcal{D}_a;K)$ and $(H_{0,\diver}^1\left(\mathcal{D}_a;K)\right)^{\prime}$ is still denoted by $<\cdot,\cdot>$. 
We can introduce, in a similar manner as before, the Riesz isomorphism $\mathcal{R}_{\diver}:H_{0,\diver}^1(\mathcal{D}_a;K)\rightarrow \left(H_{0,\diver}^1(\mathcal{D}_a;K)\right)^{\prime}$, defined by
\[
<\mathcal{R}_{\diver}A,B>=(\grad A,\grad B), \quad \forall A,B\in H_{0,\diver}^1(\mathcal{D}_a;K).
\]
The operator $\mathcal{R}_{\diver}=-P_L\Delta$, where $P_L:L^2(\mathcal{D}_a;K)\to L_{\diver}^2(\mathcal{D}_a;K)$ \pier{denotes the Leray projector, is the} negative projected Laplace operator with homogeneous Dirichlet boundary conditions, which is positive definite and self adjoint. As a consequence of the Lax--Milgram theorem and the Poincar\'e inequality, the inverse operator $(-P_L\Delta)^{-1}:\left(H_{0,\diver}^1(\mathcal{D}_a;K)\right)^{\prime}\rightarrow H_{0,\diver}^1(\mathcal{D}_a;K)$ is well defined, and we set $A:=(-P_L\Delta)^{-1}F=\mathcal{G}_{L,\diver}\ast F$, for $F\in \left(H_{0,\diver}^1(\mathcal{D}_a;K)\right)^{\prime}$, where $\mathcal{G}_{L,\diver}$ is the Green propagator associated to the projected 
Laplace operator with homogeneous Dirichlet boundary conditions, if $-P_L\Delta A=F$ in $\mathcal{D}_a$ in the weak sense, and $A=0$ on $\Gamma_a$ in the sense of traces. We again note that, if $A\in H_{0,\diver}^1(\mathcal{D}_a;K)$ solves $-P_L\Delta A=F$ for some $F\in W^{m,p}(\mathcal{D}_a;K)\cap L_{\diver}^2(\mathcal{D}_a,K)$, $1<p<\infty$, $m\in \mathbb{N}$, and $\Gamma_a$ is of class $C^{m+2}$, then from elliptic regularity theory $A\in W^{m+2,p}(\mathcal{D}_a;K)\cap L_{\diver}^2(\mathcal{D}_a,K)$ and $-P_L\Delta A=F$ a.e. in $\mathcal{D}_a$.

In the following, $C$ denotes a generic positive constant independent of the unknown variables, the discretization and the \michhhb{physical parameters}, the value of which might change from line to line; $C_1, C_2, \dots$ indicate generic positive constants whose particular value must be tracked through the calculations; $C(a,b,\dots)$ denotes a constant depending on the nonnegative parameters $a,b,\dots$.

\reviewa{
\subsection{Helmholtz--Hodge decomposition for vector fields.}
We now recall specific forms of the Helmholtz--Hodge decomposition for vector fields which will be useful in the forthcoming sections. We refer the reader to \cite{chorin,dautraylions,denaro} for their proofs. Here, we  are interested in specific decompositions obtained through a constructive procedure by means of the solution of elliptic problems with Dirichlet boundary conditions. As we will see, the aforementioned elliptic problems will be crucial in our model derivation to define the \textit{kinematic constraints} between the model variables.  
\begin{thm}
\label{thm:hhd}
Let $\mathcal{D}_a\subset \mathbb{R}^3$ be an open bounded and simply connected domain with Lipschitz and connected boundary $\Gamma_a:=\partial \mathcal{D}_a$. 
Let us introduce the spaces
\begin{align*}
& \grad H^1:=\{\mathbf{w}\in L^2(\mathcal{D}_a,\mathbb{R}^3): \exists p \in H^1(\mathcal{D}_a) \; \text{such that}\; \mathbf{w}=\grad p\}, \\
& \grad H^1_{c}:=\{\mathbf{w}\in L^2(\mathcal{D}_a,\mathbb{R}^3): \exists p \in H^1(\mathcal{D}_a) \; \text{such that}\; \mathbf{w}=\grad p, \;\; p|_{\Gamma_a}=c\}, \\
& \curl H^1:=\{\mathbf{w}\in L^2(\mathcal{D}_a,\mathbb{R}^3): \exists \mathbf{v} \in H^1(\mathcal{D}_a,\mathbb{R}^3) \; \text{such that}\; \mathbf{w}=\curl \mathbf{v}\}, \\
& H_{0,\rm{div}}:=\{\mathbf{v}\in L^2(\mathcal{D}_a,\mathbb{R}^3): \rm{div}\mathbf{v}=0, \; \mathbf{v}\cdot \mathbf{n}|_{\Gamma_a}=0\},
\end{align*}
where $c\in \mathbb{R}$ is an arbitrary constant.
For any $\boldsymbol{\xi}\in L^2(\mathcal{D}_a,\mathbb{R}^3)$, there exist a unique $\mathbf{v}\in \grad H^1$, with $\mathbf{v}=\grad p$, and a unique $\mathbf{r}\in H_{0,div}$ such that
\begin{equation}
\label{hh5}
\boldsymbol{\xi}=\mathbf{v}+\mathbf{r}=\grad p + \mathbf{r},
\end{equation}
i.e. the following decomposition is valid 
\begin{equation}
\label{hh5bis}
L^2(\mathcal{D}_a,\mathbb{R}^3)=\grad H^1\oplus H_{0,\rm{div}}.
\end{equation}
Moreover, there exist a unique $\mathbf{w}\in \grad H^1_{c}$, with $\mathbf{w}=\grad u$, and a unique $\mathbf{q}\in \curl H^1$, with $\mathbf{q}=\curl \boldsymbol{\alpha}$, such that 
\begin{equation}
\label{hh6}
\boldsymbol{\xi}=\mathbf{w}+\mathbf{q}=\grad u + \curl \boldsymbol{\alpha},
\end{equation}
i.e. the following decomposition is valid 
\begin{equation}
\label{hh6bis}
L^2(\mathcal{D}_a,\mathbb{R}^3)= \grad H^1_{c}\oplus \curl H^1.
\end{equation}
\end{thm}
\begin{rem}
\label{rem:hhsc}
The hypotheses that  $\mathcal{D}_a$ is simply connected and that $\Gamma$ is connected are made to simplify the presentation of the results. The theorem could be extended in a standard manner to a (not simply) connected domain $\mathcal{D}_a$ with boundary constituted by a finite number of connected components by topological arguments as done in \cite[Chapter IX]{dautraylions}. Anyhow, the situation in which the initial form of the body is topologically simply connected until it develops cuts or holes is mechanically meaningfull.
\end{rem}
The decomposition \eqref{hh5bis} is proved in \cite[Chapter IX]{dautraylions} by proving the closedness of $\grad H^1$ in $L^2(\mathcal{D}_a,\mathbb{R}^3)$ and the orthogonality between $\grad H^1$ and $H_{0,\rm{div}}$ in the $L^2(\mathcal{D}_a,\mathbb{R}^3)$ topology. The existence of the decomposition can be proved also in a constructive way by solving an elliptic problem with Neumann boundary conditions for $p$ in \eqref{hh5} \cite{chorin,denaro}, i.e. by solving
\[
\begin{cases}
\Delta p = \rm{div}\boldsymbol{\xi},\\
\grad p \cdot \mathbf{n}|_{\Gamma_a}=\boldsymbol{\xi}\cdot \mathbf{n}|_{\Gamma_a},
\end{cases}
\]
where $\mathbf{n}$ is the outward unit normal vector to $\Gamma_a$, and then setting $\mathbf{r}=\boldsymbol{\xi}-\grad p\in H_{0,\rm{div}}$. The decomposition \eqref{hh6bis} is proved in \cite[Chapter IX]{dautraylions} by first proving the decomposition $L^2(\mathcal{D}_a,\mathbb{R}^3)= \grad H^1_{0}\oplus H_{\rm{div}}$, where $H_{\rm{div}}:=\{\mathbf{v}\in L^2(\mathcal{D}_a,\mathbb{R}^3): \rm{div}\mathbf{v}=0\}$, and then identifying $\curl H^1$ as a proper subspace of $H_{\rm{div}}$. The existence of the decomposition can be proved also in a constructive way by solving an elliptic problem with Neumann boundary conditions for $\boldsymbol{\alpha}$ in \eqref{hh6} \cite{denaro}, considering moreover the constraint $\rm{div}\boldsymbol{\alpha}=0$, which is not reductive since $\boldsymbol{\alpha}$ is uniquely defined up to the gradient of a scalar function. Indeed, existence can be proved by solving
\begin{equation}
\label{ellneum}
\begin{cases}
-\Delta \boldsymbol{\alpha} = \curl \boldsymbol{\xi},\\
\curl \boldsymbol{\alpha} \wedge \mathbf{n}|_{\Gamma_a}=\boldsymbol{\xi}\wedge \mathbf{n}|_{\Gamma_a},
\end{cases}
\end{equation}
and then setting $\mathbf{w}=\grad u=\boldsymbol{\xi}-\curl \boldsymbol{\alpha}\in \grad H^1_{c}$. Note that $\mathbf{w}\wedge \mathbf{n}|_{\Gamma_a}=\mathbf{0}$ implies that $u|_{\Gamma_a}=c$, for any given $c\in \mathbb{R}$. By expressing the vector component $\boldsymbol{\alpha}$ in \eqref{hh6} through a Helmholtz--Hodge decomposition of type \eqref{hh5bis}, it's possible to obtain a further generalized Helmholtz--Hodge decomposition of the form (see e.g. \cite{sprossig}): for any $\boldsymbol{\xi}\in L^2(\mathcal{D}_a,\mathbb{R}^3)$, there exist a unique $\mathbf{w}\in \grad H^1_{c}$, with $\mathbf{w}=\grad u$, and a unique $\mathbf{q}\in \curl \left(H_{0,\rm{div}}\cap H^1(\mathcal{D}_a,\mathbb{R}^3)\right)$, with $\mathbf{q}=\curl \mathbf{d}$ and $\mathbf{d}\in H_{0,\rm{div}}\cap H^1(\mathcal{D}_a,\mathbb{R}^3)$, such that 
\begin{equation}
\label{hhdir}
\boldsymbol{\xi}=\mathbf{w}+\mathbf{q}=\grad u + \curl \mathbf{d}.
\end{equation}
Taking the divergence and the curl of \eqref{hhdir}, choosing $c=0$ and substituting the slip boundary condition $\mathbf{d}\cdot \mathbf{n}|_{\Gamma_a}=0$ with the no-slip condition $\mathbf{d}|_{\Gamma_a}=\boldsymbol{0}$, we can construct the decomposition \eqref{hhdir} by solving the following elliptic problems with Dirichlet boundary conditions
\begin{equation}
\label{elldir}
\begin{cases}
\Delta u = \rm{div}\boldsymbol{\xi},\\
u|_{\Gamma_a}=0,
\end{cases}
\begin{cases}
-\Delta \mathbf{d} = \rm{curl}\boldsymbol{\xi},\\
\mathbf{d}|_{\Gamma_a}=\mathbf{0}.
\end{cases}
\end{equation}
}
\subsection{Functional inequalities}
We recall the Gagliardo-Nirenberg inequality {(see e.g. \cite{gagliardo,nirenberg,leoni})}.
\begin{lem}
\label{lem:gagliardoniremberg}
Let $\mathcal{D} \subset \mathbb{R}^3$ be a bounded domain with Lipschitz boundary and $f\in W^{m,r}\cap L^q$, $q\geq 1$, $r\leq \infty$, where $f$ can be a function with scalar, vectorial or tensorial values. For any integer $j$ with $0 \leq j < m$, suppose there is $\alpha \in \mathbb{R}$ such that
\[
j-\frac{3}{p}=\left(m-\frac{3}{r}\right)\alpha+(1-\alpha)\left(-\frac{3}{q}\right), \quad \frac{j}{m}\leq \alpha \leq 1.
\]
Then, there exists a positive constant $C$ depending on $\Omega$, m, j, q, r, and $\alpha$ such that
\begin{equation}
\label{eqn:2}
\pier{\|}D^jf\pier{\|}_{L^p}\leq C\pier{\|}f\pier{\|}_{W^{m,r}}^{\alpha}\pier{\|}f\pier{\|}_{L^q}^{1-\alpha}.
\end{equation}
\end{lem}
%
Finally, we will use the following result.
\begin{lem}
 \label{lem:omegalp}
 Let $p\geq 1$ and $\boldsymbol{\Omega}_1,\boldsymbol{\Omega}_2\in L^p\left(\mathcal{D}_a,Skew\left(\mathbb{R}^{3\times 3}\right)\right)$. There exists a positive constant $C$ such that
 \begin{equation}
  \label{eqn:omegalp}
  \left\|\mathrm{e}^{\boldsymbol{\Omega}_1}-\mathrm{e}^{\boldsymbol{\Omega}_2}\right\|_{L^p\left(\mathcal{D}_a,\mathbb{R}^{3\times 3}\right)}\leq C\left\|{\boldsymbol{\Omega}_1}-{\boldsymbol{\Omega}_2}\right\|_{L^p\left(\mathcal{D}_a,Skew\left(\mathbb{R}^{3\times 3}\right)\right)}.
 \end{equation}
\end{lem}
\begin{pf}
 We introduce the three Euler angles $\theta,\phi,\chi$, associated to a skew symmetric tensor $\boldsymbol{\Omega}\in \mathcal{A}$, and the three matrices $\mathbf{A},\mathbf{B},\mathbf{C}\in Skew(\mathbf{R}^{3\times 3})$ which are elements of the canonical basis for $Skew(\mathbf{R}^{3\times 3})$, \pier{i.e.,}
\[\mathbf{A}=
\begin{pmatrix}
0 & -1 & 0 \\
1 & 0 & 0\\
0 & 0 & 0
\end{pmatrix}
, \quad 
\mathbf{B}=
\begin{pmatrix}
0 & 0 & 0 \\
0 & 0 & -1\\
0 & 1 & 0
\end{pmatrix}
,\quad 
\mathbf{C}=
\begin{pmatrix}
0 & 0 & -1 \\
0 & 0 & 0\\
1 & 0 & 0
\end{pmatrix}.
\]
Then we may write, for any $\mathbf{x}\in \mathcal{D}_a$, 
\[
\boldsymbol{\Omega}(\mathbf{x})=\theta(\mathbf{x}) \mathbf{A}+\phi(\mathbf{x}) \mathbf{B}+\chi(\mathbf{x}) \mathbf{C}.
\]

\noindent
\abramo{
Observing the fact that, for any $n\in \mathbb{N}$,  
\[
 \mathbf{A}^{2n+1}=(-1)^n\mathbf{A},\quad \mathbf{A}^{2n+2}=(-1)^{n+1}\begin{pmatrix}
1 & 0 & 0 \\
0 & 1 & 0\\
0 & 0 & 0
\end{pmatrix},
\]
with similar relations for $\mathbf{B}$ and $\mathbf{C}$, we have that
}
\begin{align*}
&\mathrm{e}^{\boldsymbol{\Omega}_1(\mathbf{x})}-\mathrm{e}^{\boldsymbol{\Omega}_2(\mathbf{x})}=\mathrm{e}^{\theta_1(\mathbf{x}) \mathbf{A}}\mathrm{e}^{\phi_1(\mathbf{x}) \mathbf{B}}\mathrm{e}^{\chi_1(\mathbf{x}) \mathbf{C}}-\mathrm{e}^{\theta_2(\mathbf{x}) \mathbf{A}}\mathrm{e}^{\phi_2(\mathbf{x}) \mathbf{B}}\mathrm{e}^{\chi_2(\mathbf{x}) \mathbf{C}}\\
&= \begin{pmatrix}
\cos(\theta_1(\mathbf{x})) & -\sin(\theta_1(\mathbf{x})) & 0 \\
\sin(\theta_1(\mathbf{x})) & \cos(\theta_1(\mathbf{x})) & 0\\
0 & 0 & 0
\end{pmatrix}
\begin{pmatrix}
0 & 0 & 0 \\
0 & \cos(\phi_1(\mathbf{x})) & -\sin(\phi_1(\mathbf{x}))\\
0 & \sin(\phi_1(\mathbf{x})) & \cos(\phi_1(\mathbf{x}))
\end{pmatrix}
\\
&\times
\begin{pmatrix}
\cos(\chi_1(\mathbf{x})) & 0 & -\sin(\chi_1(\mathbf{x})) \\
0 & 0 & 0\\
\sin(\chi_1(\mathbf{x})) & 0 & \cos(\chi_1(\mathbf{x}))
\end{pmatrix}
-
\begin{pmatrix}
\cos(\theta_2(\mathbf{x})) & -\sin(\theta_2(\mathbf{x})) & 0 \\
\sin(\theta_2(\mathbf{x})) & \cos(\theta_2(\mathbf{x})) & 0\\
0 & 0 & 0
\end{pmatrix}
\\
&\times 
\begin{pmatrix}
0 & 0 & 0 \\
0 & \cos(\phi_2(\mathbf{x})) & -\sin(\phi_2(\mathbf{x}))\\
0 & \sin(\phi_2(\mathbf{x})) & \cos(\phi_2(\mathbf{x}))
\end{pmatrix}
\begin{pmatrix}
\cos(\chi_2(\mathbf{x})) & 0 & -\sin(\chi_2(\mathbf{x})) \\
0 & 0 & 0\\
\sin(\chi_2(\mathbf{x})) & 0 & \cos(\chi_2(\mathbf{x}))
\end{pmatrix}.
\end{align*}
The bound \eqref{eqn:omegalp} is thus a consequence of the uniform Lipschitz continuity and of the uniform boundedness of the $\cos$ and $\sin$ functions. 
\end{pf}

\section{Model derivation} 
\label{sec:model}
We consider the motion of a deformable elastic solid in $\mathcal{D}_a$ which is fixed on its boundary $\Gamma_a:=\partial \mathcal{D}_a$.
In the time interval $(0,T)$, the motion is described by the map
\[
(\mathbf{a},t)\rightarrow \mathbf{\Phi}(\mathbf{a},t)\reviewa{=\mathbf{a}+\mathbf{u}(\mathbf{a},t)}\in \mathbb{R}^3, \quad (\mathbf{a},t)\in \mathcal{D}_{aT}\pier{\, :=  \mathcal{D}_{a}\times (0,T)},
\]
with 
\[
\reviewa{\mathbf{u}(\mathbf{a},0)=\mathbf{0}}\;\; \pier{\text{for}}\;\mathbf{a}\in \mathcal{D}_a \quad \text{and} \quad \reviewa{\mathbf{u}(\mathbf{a},t)=\mathbf{0}} \;\; \text{for} \; \mathbf{a}\in \Gamma_a.
\]
We assume that the motion is not compatible, \pier{i.e.,} there exist a dislocation tensor $\mathbf{Z}\in \mathcal{M}_{\diver}$, with $\mathbf{Z}(\mathbf{a},0)=\mathbf{0}$ \pier{for} $\mathbf{a}\in \mathcal{D}_a$ and $\mathbf{Z}(\mathbf{a},t)=\mathbf{0}$ for $\mathbf{a}\in \Gamma_a$, such that
\begin{equation}
\label{eqn:3}
\reviewa{
\grad{\mathbf{u}}=(\mathbf{R}\mathbf{W}-\mathbf{I})- \curl \mathbf{Z},
}
\end{equation}
where $\mathbf{R}\in \mathcal{SO}$ is the rotation tensor and $\mathbf{W}\in \mathcal{S}$ is the stretch tensor associated to the deformation gradient tensor, \abramo{with $\mathbf{R}\mathbf{W}(\mathbf{a},0)=\mathbf{I}$ for $\mathbf{a}\in \mathcal{D}_a$, $\mathbf{R}(\mathbf{a},t)=\mathbf{W}(\mathbf{a},t)=\mathbf{I}$ 
\pco{{}for $\mathbf{a}\in\Gamma_a\times (0,T)$}}. 
\reviewa{
Since the $\grad, \diver, \curl$ operators are applied to second order tensors row-wise, we observe that the existence of the decomposition \eqref{eqn:3} is a consequence of the application of Theorem \ref{thm:hhd}, in particular of formula \eqref{hhdir}, to the row vectors of the involved tensors. Given $\mathbf{R}$ and $\mathbf{W}$, the components $\boldsymbol{\Phi}=\mathbf{a}+\mathbf{u}$ and $\mathbf{Z}$ in the decomposition \eqref{eqn:3} may be obtained as in \eqref{elldir}, i.e. solving elliptic problems with Dirichlet boundary conditions derived by applying the divergence and the curl operators to \eqref{eqn:3}, ending with the kinematic relations:
\begin{equation}
\label{eqn:4}
\Delta {\mathbf{\Phi}}=\diver\left(\mathbf{R}\mathbf{W}\right),
\end{equation}
endowed with the boundary conditions ${\mathbf{\Phi}}(\mathbf{a},t)=\mathbf{a}$ {for $(\mathbf{a},t)\in\Gamma_a\times (0,T)$}, and
\begin{equation}
\label{eqn:5}
- \Delta {\mathbf{Z}}=\curl \left(\mathbf{R}\mathbf{W}\right), \quad \diver{\mathbf{Z}}=\mathbf{0},
\end{equation}
endowed with the boundary condition ${\mathbf{Z}}(\mathbf{a},t)=\mathbf{0}$ for $(\mathbf{a},t)\in\Gamma_a\times (0,T)$.
\begin{rem}
In our theoretical framework, describing the deformation of an elastic solid which is fixed on its boundary, we have considered homogeneous Dirichlet boundary conditions for the dislocations $\mathbf{Z}$. We may think of dislocations as motions at the microscopic level, for instance motions of atoms in lattices, which result from macroscopic motions. Hence, no macroscopic motion on the boundary of the solid implies that the dislocations do not change on the boundary during the motion. In this situation, the components $\boldsymbol{\Phi}=\mathbf{a}+\mathbf{u}$ and $\mathbf{Z}$ in the decomposition \eqref{eqn:3} are obtained by solving elliptic problems with homogeneous Dirichlet boundary conditions as in \eqref{elldir}.
We observe that dislocations may also result from microscopic actions, for instance radiative actions (example in nuclear plants radiations produce damage), for 
instance thermal actions (example in shape memory alloys, the atoms move due to the thermal evolution). Then an external action may 
produce a flux of dislocations on the fixed boundary of the solid. In the latter situation, the component $\mathbf{Z}$ in the decomposition \eqref{eqn:3} may be obtained by solving an elliptic problem with Neumann boundary conditions as in \eqref{ellneum}, while the component $\boldsymbol{\Phi}=\mathbf{a}+\mathbf{u}$ is still obtained by solving an elliptic problem with Dirichlet boundary conditions as in \eqref{elldir}.
\end{rem}
}
Taking moreover the time derivative of \eqref{eqn:4} and \eqref{eqn:5}, introducing also the velocity vector field $\mathbf{U}:=\overbigdot{\mathbf{\Phi}}\in \mathcal{V}$ and the angular velocity tensor $\boldsymbol{\Omega}:=\overbigdot{\mathbf{R}}\mathbf{R}^T\in \mathcal{A}$, we obtain the \pierhhb{kinematic} relations:
\begin{equation}
\label{eqn:6}
\Delta \mathbf{U}=\diver\left(\mathbf{R}\overbigdot{\mathbf{W}}+\boldsymbol{\Omega}\mathbf{R}\mathbf{W}\right),
\end{equation}
endowed with the boundary condition $\mathbf{U}(\mathbf{a},t)=\mathbf{0}$ on $\Gamma_a$, and
\begin{equation}
\label{eqn:7}
-\Delta \overbigdot{\mathbf{Z}}=\curl \left(\mathbf{R}\overbigdot{\mathbf{W}}+\boldsymbol{\Omega}\mathbf{R}\mathbf{W}\right), \quad \diver\overbigdot{\mathbf{Z}}=\mathbf{0},
\end{equation}
endowed with the boundary condition $\overbigdot{\mathbf{Z}}(\mathbf{a},t)=\mathbf{0}$ on $\Gamma_a$.

We derive the model equations from the principle of virtual power, which gives the \michhhb{equations of motion} for the linear and angular momenta expressed in terms of the \pierhhb{kinematic} variables and  \michhhb{internal force} tensors. We then constitutively assign the form of the \michhhb{internal force} tensors, in terms of the \pierhhb{kinematic} variables, in order for \pier{the system} to satisfy the Clausius--Duhem dissipative equality.
We start by defining the set $\mathfrak{C}$ of virtual velocities. Given $\mathbf{R}\in \mathcal{SO}$, $\mathbf{W}\in \mathcal{S}$, we define, for any $t\in (0,T)$, the set
\begin{align}
\label{eqn:8}
& \displaystyle \mathfrak{C}:=\biggl\{\left(\mathbf{V},\nwhat{\mathbf{W}},\nwhat{\boldsymbol{\Omega}},\nwhat{\mathbf{Z}}\right)\in (\mathcal{V},\mathcal{M},\mathcal{M},\mathcal{M}_{\diver})\ \biggl| \quad  \nwhat{\mathbf{W}}|_{\Gamma_a}=\nwhat{\boldsymbol{\Omega}}|_{\Gamma_a}=\boldsymbol{0},  \;\; \reviewa{\grad\nwhat{\mathbf{W}}|_{\Gamma_a}=\grad\nwhat{\boldsymbol{\Omega}}|_{\Gamma_a}=\mathbf{0}}, \notag\\
& \displaystyle \qquad\quad
\begin{cases}
\Delta \mathbf{V}=\diver\left(\mathbf{R}\nwhat{\mathbf{W}}+\nwhat{\boldsymbol{\Omega}}\mathbf{R}\mathbf{W}\right),\\
\mathbf{V}=\mathbf{0} \quad \text{on} \; \Gamma_a\; ,
\end{cases}
\begin{cases}
-P_L\Delta \nwhat{\mathbf{Z}}=\curl \left(\mathbf{R}\nwhat{\mathbf{W}}+\nwhat{\boldsymbol{\Omega}}\mathbf{R}\mathbf{W}\right),\\
\nwhat{\mathbf{Z}}=\mathbf{0} \quad \text{on} \; \Gamma_a.
\end{cases}
\!\!\!\biggr\}
\end{align}
The virtual velocities then satisfy the following constraint, which, \abramohh{similarly to \eqref{eqn:3}, is a consequence of \eqref{hhdir} applied row-wise:} 
\begin{equation}
\label{eqn:9}
\grad{\mathbf{V}}=\mathbf{R}\nwhat{\mathbf{W}}+\nwhat{\boldsymbol{\Omega}}\mathbf{R}\mathbf{W}-\curl \nwhat{\mathbf{Z}}.
\end{equation}
We observe that the set $\mathfrak{C}$ of virtual velocities is defined in terms of the variables $\mathbf{R}$ and $\mathbf{W}$, and hence depend on the solutions of the \michhhb{equations of motion}. We can formally write
\begin{equation}
\label{eqn:10}
\mathbf{V}=-\mathcal{G}_L\ast\diver\left(\mathbf{R}\nwhat{\mathbf{W}}+\nwhat{\boldsymbol{\Omega}}\mathbf{R}\mathbf{W}\right),
\end{equation}
and
\begin{equation}
\label{eqn:11}
\nwhat{\mathbf{Z}}=\mathcal{G}_{L,\diver}\ast\curl\left(\mathbf{R}\nwhat{\mathbf{W}}+\nwhat{\boldsymbol{\Omega}}\mathbf{R}\mathbf{W}\right).
\end{equation}
Given solutions with regularity, for a.e. $t\in (0,T)$, $\mathbf{R}\in H^2(\mathcal{D}_a;\mathbb{R}^{3\times 3})\cap \mathcal{SO}$, $\mathbf{W}\in H^2(\mathcal{D}_a;\mathbb{R}^{3\times 3})\cap \mathcal{S}$, and choosing $\nwhat{\mathbf{W}},\nwhat{\boldsymbol{\Omega}}\in H^1(\mathcal{D}_a;\mathbb{R}^{3\times 3})$, from elliptic regularity theory and with the assumed regularity of $\Gamma_a$ we get that $\mathbf{V}\in H^2(\mathcal{D}_a;\mathbb{R}^{3})\cap H_{0}^1(\mathcal{D}_a;\mathbb{R}^{3})$ and $\nwhat{\mathbf{Z}}\in H^2(\mathcal{D}_a;\mathbb{R}^{3\times 3})\cap H_{0,\diver}^1(\mathcal{D}_a;\mathbb{R}^{3})$.

We now introduce the virtual power of internal forces $p_{\text{int}}(\mathcal{D}_a,\mathcal{C})$, the virtual power of external forces $p_{\text{ext}}(\mathcal{D}_a,\mathcal{C})$ and the virtual power of acceleration forces $p_{\text{acc}}(\mathcal{D}_a,\mathcal{C})$, defined in terms of $\mathcal{D}_a$ and of an element $\mathcal{C}\in \mathfrak{C}$. The principle of virtual power then states that
\begin{equation}
\label{eqn:12}
p_{\text{acc}}(\mathcal{D}_a,\mathcal{C})=p_{\text{int}}(\mathcal{D}_a,\mathcal{C})+p_{\text{ext}}(\mathcal{D}_a,\mathcal{C}) \quad \forall \mathcal{C}\in \mathfrak{C}.
\end{equation}

\noindent
The virtual power of internal forces is defined as
\begin{align}
\label{eqn:13}
& \displaystyle p_{\text{int}}(\mathcal{D}_a,\mathcal{C}):=-\int_{\mathcal{D}_a}\left(\boldsymbol{\Pi}:\grad \mathbf{V}+\mathbf{X}:: \grad \nwhat{\mathbf{W}}+\mathbf{Y}::: \grad \grad \nwhat{\mathbf{W}}\right) \notag \\
& \displaystyle +\frac{1}{2}\int_{\mathcal{D}_a}\left(\mathbf{M}:\nwhat{\boldsymbol{\Omega}}-\boldsymbol{\Lambda}::\grad \nwhat{\boldsymbol{\Omega}}-\mathbf{C}:::\grad \grad \nwhat{\boldsymbol{\Omega}}\right)+\int_{\mathcal{D}_a} \ap{\boldsymbol{\Gamma}:\curl \nwhat{\mathbf{Z}}},
\end{align}
where $\boldsymbol{\Pi}$ is the Piola--Kirchhoff--Boussinesq stress tensor, $\mathbf{M}$ represents the momentum, $\boldsymbol{\Lambda}$ the momentum flux and $\mathbf{C}$ the flux of the momentum flux. The quantities $\mathbf{X}, \mathbf{Y}, \boldsymbol{\Gamma}$ are new \michhhb{internal force} tensors associated to the kinematic variables $\mathbf{W}$ and $\mathbf{Z}$. \abramohh{In particular, $\boldsymbol{\Gamma}$ is an internal force accounting for the evolution
of the dislocations.} The virtual power of external forces is defined as
\begin{equation}
\label{eqn:14}
p_{\text{ext}}(\mathcal{D}_a,\mathcal{C}):=\int_{\mathcal{D}_a}\boldsymbol{\mathcal{W}}_{\text{ext}}:\nwhat{\mathbf{W}}+ \int_{\mathcal{D}_a}\boldsymbol{\Omega}_{\text{ext}}:\nwhat{\boldsymbol{\Omega}},
\end{equation}
where $\boldsymbol{\mathcal{W}}_{\text{ext}}$ and $\boldsymbol{\Omega}_{\text{ext}}$ are external forces, \abramohhb{possibly depending on $\mathbf{W}$ and $\boldsymbol{R}$}, which perform work by stretching and rotating the system, respectively. \abramohhb{Note that the reader may expect a factor $\frac{1}{2}$ in front of the second integral in \eqref{eqn:14}; in fact, we choose to incorporate this factor in the definition of $\boldsymbol{\Omega}_{\text{ext}}$.}
\abramohhb {\begin{rem}
We note that, since $\mathbf{V}$ can be expressed in terms of $\nwhat{\mathbf{W}}$ and $\nwhat{\boldsymbol{\Omega}}$ through \eqref{eqn:10}, the expression \eqref{eqn:14} may include external powers for classical body forces like gravity, follower forces depending on the solutions \cite{antman}, pressure contributions depending on ${\rm cof} \mathbf{W}$, and so on.
\end{rem}
} 
Finally, the virtual power of acceleration forces is defined as 
\begin{align}
\label{eqn:15}
& \displaystyle \notag  p_{\text{acc}}(\mathcal{D}_a,\mathcal{C}):=\int_{\mathcal{D}_a} 
\biggl( \frac{d\mathbf{U}}{dt} \cdot \mathbf{V} + \twobigdot{\mathbf{W}}:
\nwhat{\mathbf{W}}+\overbigdot{\boldsymbol{\Omega}}:\nwhat{\boldsymbol{\Omega}}+\reviewa{\grad \grad} \twobigdot{\mathbf{W}}: \reviewa{\grad \grad} \nwhat{\mathbf{W}}\\
& \displaystyle+
\reviewa{\grad \grad}  \overbigdot{\boldsymbol{\Omega}}:\reviewa{\grad \grad}  \nwhat{\boldsymbol{\Omega}}\biggr).
\end{align}
\abramohh{As discussed in the Introduction, higher order terms in the virtual power of acceleration forces are introduced to be able to deal with a point with inertia, which requires regularity in space and time of the angular velocity and acceleration variables.}
Using \eqref{eqn:9} in \eqref{eqn:13} we obtain that
\begin{align}
\label{eqn:16}
& \displaystyle p_{\text{int}}(\mathcal{D}_a,\mathcal{C}):=-\int_{\mathcal{D}_a}\left(\mathbf{R}^T\boldsymbol{\Pi}:\nwhat{\mathbf{W}}+\mathbf{X}:: \grad \nwhat{\mathbf{W}}+\mathbf{Y}::: \grad \grad \nwhat{\mathbf{W}}\right) \notag \\
& \displaystyle \notag +\frac{1}{2}\int_{\mathcal{D}_a}\left((\mathbf{M}-2\boldsymbol{\Pi}\mathbf{W}\mathbf{R}^T):\nwhat{\boldsymbol{\Omega}}-\boldsymbol{\Lambda}::\grad \nwhat{\boldsymbol{\Omega}}-\mathbf{C}:::\grad \grad \nwhat{\boldsymbol{\Omega}}\right) \notag
\\ 
&+\int_{\mathcal{D}_a}(\curl \boldsymbol{\Gamma}+\curl \boldsymbol{\Pi}):\nwhat{\mathbf{Z}},
\end{align}
where in the last term we have used integration by parts and the boundary conditions for $\nwhat{\mathbf{Z}}$.

We rewrite the first term in \eqref{eqn:15} employing \eqref{eqn:10} \pier{and then} obtaining
\begin{align}
\label{eqn:17}
& \displaystyle -\int_{\mathcal{D}_a}\frac{d\mathbf{U}}{dt}\cdot \left(\mathcal{G}_L\ast\diver\left(\mathbf{R}\nwhat{\mathbf{W}}+\nwhat{\boldsymbol{\Omega}}\mathbf{R}\mathbf{W}\right)\right)\notag\\
&\displaystyle =-\int_{\mathcal{D}_a}\left(\mathcal{G}_L\ast\frac{d\mathbf{U}}{dt}\right)\cdot \diver\left(\mathbf{R}\nwhat{\mathbf{W}}+\nwhat{\boldsymbol{\Omega}}\mathbf{R}\mathbf{W}\right)=\int_{\mathcal{D}_a}\grad \left(\mathcal{G}_L\ast\frac{d\mathbf{U}}{dt}\right): \left(\mathbf{R}\nwhat{\mathbf{W}}+\nwhat{\boldsymbol{\Omega}}\mathbf{R}\mathbf{W}\right)\pier{,}
\end{align}
where in the last term we have integrated by parts and used the boundary conditions for $\mathcal{G}_L\ast\frac{d\mathbf{U}}{dt}$. We also use \eqref{eqn:11} in the last term of \eqref{eqn:16} \pier{and deduce that}
\begin{align}
\label{eqn:18}
& \displaystyle \int_{\mathcal{D}_a}(\curl \boldsymbol{\Gamma}+\curl \boldsymbol{\Pi}):\nwhat{\mathbf{Z}}=\int_{\mathcal{D}_a}\curl(\boldsymbol{\Gamma}+\boldsymbol{\Pi}):\left(\mathcal{G}_{L,\diver}\ast\curl\left(\mathbf{R}\nwhat{\mathbf{W}}+\nwhat{\boldsymbol{\Omega}}\mathbf{R}\mathbf{W}\right)\right) \notag \\
&\displaystyle \notag \int_{\mathcal{D}_a}\mathcal{G}_{L,\diver}\ast \curl(\boldsymbol{\Gamma}+\boldsymbol{\Pi}):\curl\left(\mathbf{R}\nwhat{\mathbf{W}}+\nwhat{\boldsymbol{\Omega}}\mathbf{R}\mathbf{W}\right)\\
&= \displaystyle\int_{\mathcal{D}_a}\curl\left(\mathcal{G}_{L,\diver}\ast \curl(\boldsymbol{\Gamma}+\boldsymbol{\Pi})\right):\left(\mathbf{R}\nwhat{\mathbf{W}}+\nwhat{\boldsymbol{\Omega}}\mathbf{R}\mathbf{W}\right),
\end{align}
where in the last term we have integrated by parts and used the boundary conditions for $\mathcal{G}_{L,\diver}\ast \curl(\boldsymbol{\Gamma}+\boldsymbol{\Pi})$. Inserting \eqref{eqn:14}-\eqref{eqn:18} in \eqref{eqn:12} and integrating by parts, the principle of virtual power becomes: given $\mathbf{R}\in \mathcal{SO}$, $\mathbf{W}\in \mathcal{S}$,
\begin{align}
\label{eqn:19}
& \displaystyle \notag \int_{\mathcal{D}_a}\biggl(\mathbf{R}^T\grad \left(\mathcal{G}_L\ast\frac{d\mathbf{U}}{dt}\right)+\twobigdot{\mathbf{W}}+\pier{{}\reviewa{\diver \Delta \grad}{}} \twobigdot{\mathbf{W}}-\mathbf{R}^T\curl\left(\mathcal{G}_{L,\diver}\ast \curl(\boldsymbol{\Gamma}+\boldsymbol{\Pi})\right)\pier{\biggr){}\colon \nwhat{\mathbf{W}}}\\
& \displaystyle \notag
\pier{{}+\int_{\mathcal{D}_a}(\mathbf{R}^T\boldsymbol{\Pi}
 -\diver \mathbf{X}+\diver \diver \mathbf{Y})\colon \nwhat{\mathbf{W}}}
+\int_{\mathcal{D}_a}\biggl(\grad \left(\mathcal{G}_L\ast\frac{d\mathbf{U}}{dt}\right)\mathbf{W}\mathbf{R}^T
+\overbigdot{\boldsymbol{\Omega}}+\pier{{}\reviewa{\diver \Delta \grad}{}} \overbigdot{\boldsymbol{\Omega}}\pier{\biggr)\colon \nwhat{\boldsymbol{\Omega}}}
\\
&\displaystyle \notag\pier{{}-\int_{\mathcal{D}_a} \frac12 \Bigl({} 2 \curl\left(\mathcal{G}_{L,\diver}\ast \curl(\boldsymbol{\Gamma}+\boldsymbol{\Pi})\right)\mathbf{W}\mathbf{R}^T +(\mathbf{M}-2\boldsymbol{\Pi}\mathbf{W}\mathbf{R}^T)+\diver \boldsymbol{\Lambda}-\diver \diver \mathbf{C}\Bigr)}\colon \nwhat{\boldsymbol{\Omega}}
\\
& \displaystyle \notag+\int_{\Gamma_a}\left(\mathbf{X}-\diver \mathbf{Y}\reviewa{-\Delta \grad \twobigdot{\mathbf{W}}}\right)\mathbf{N}:\nwhat{\mathbf{W}}+\int_{\Gamma_a}\left(\mathbf{Y}+\reviewa{\grad \grad \twobigdot{\mathbf{W}}}\right)\mathbf{N}::\grad\nwhat{\mathbf{W}}\\
& \displaystyle \notag +\frac{1}{2}\int_{\Gamma_a}\left(\boldsymbol{\Lambda}-\diver \mathbf{C}\reviewa{-\Delta \grad \overbigdot{\boldsymbol{\Omega}}}\right)\mathbf{N}:\nwhat{\boldsymbol{\Omega}}+\frac{1}{2}\int_{\Gamma_a}\left(\mathbf{C}+\reviewa{\grad \grad \overbigdot{\boldsymbol{\Omega}}}\right)\mathbf{N}::\grad\nwhat{\boldsymbol{\Omega}}\\
& \displaystyle =\int_{\mathcal{D}_a}\mathbf{W}_{\text{ext}}\colon \nwhat{\mathbf{W}}+\int_{\mathcal{D}_a}\boldsymbol{\Omega}_{\text{ext}}:\nwhat{\boldsymbol{\Omega}},
\end{align}
for all \abramo{virtual velocities} $\nwhat{\mathbf{W}},\nwhat{\boldsymbol{\Omega}}$, where $\mathbf{N}$ is the outward normal to $\Gamma_a$.
Assuming regularity of the integrands in \eqref{eqn:19}, and \reviewa{considering the boundary conditions assigned to the virtual velocities $\nwhat{\mathbf{W}}$ and $\nwhat{\boldsymbol{\Omega}}$}, the principle of virtual power implies the following  \michhhb{equations}, valid in $\mathcal{D}_{aT}$, which are coupled to the kinematic relations \eqref{eqn:4} and \eqref{eqn:5}:

\begin{equation}
\label{eqn:20}
\begin{cases}
\displaystyle \mathbf{R}^T\grad \left(\mathcal{G}_L\ast\frac{d\mathbf{U}}{dt}\right)+\twobigdot{\mathbf{W}}+\pier{{}\reviewa{\diver \Delta \grad}{}} \twobigdot{\mathbf{W}}-\mathbf{R}^T\curl\left(\mathcal{G}_{L,\diver}\ast \curl(\boldsymbol{\Gamma}+\boldsymbol{\Pi})\right)+\mathbf{R}^T\boldsymbol{\Pi}\\
\displaystyle -\diver \mathbf{X}+\diver \diver \mathbf{Y}=\mathbf{W}_{\text{ext}},\\ \\
\displaystyle \reviewa{\mathbf{W}=\mathbf{I}, \;\; \grad \mathbf{W}=\mathbf{0}} \quad \text{on} \; \Gamma_a\times (0,T),\\ \\
\displaystyle \grad \left(\mathcal{G}_L\ast\frac{d\mathbf{U}}{dt}\right)\mathbf{W}\mathbf{R}^T+\overbigdot{\boldsymbol{\Omega}}+\pier{{}\reviewa{\diver \Delta \grad}{}}\overbigdot{\boldsymbol{\Omega}}-\curl\left(\mathcal{G}_{L,\diver}\ast \curl(\boldsymbol{\Gamma}+\boldsymbol{\Pi})\right)\mathbf{W}\mathbf{R}^T \\
\displaystyle -\frac{1}{2}(\mathbf{M}-2\boldsymbol{\Pi}\mathbf{W}\mathbf{R}^T)-\frac{1}{2}\diver \boldsymbol{\Lambda}+\frac{1}{2}\diver \diver \mathbf{C}=\boldsymbol{\Omega}_{\text{ext}},\\ \\
\displaystyle \reviewa{\boldsymbol{\Omega}=\mathbf{0}, \;\; \grad \boldsymbol{\Omega}=\mathbf{0}} \quad \text{on} \; \Gamma_a\times (0,T),\\ \\
\overbigdot{\mathbf{R}}=\boldsymbol{\Omega}\mathbf{R},\\ \\
\Delta \pier{\mathbf{\Phi}}=\diver\left(\mathbf{R}\mathbf{W}\right), \quad  \pier{\mathbf{\Phi}(\mathbf{a}, t) =\mathbf{a} \quad \text{for} \; (\mathbf{a}, t) \in \Gamma_a\times (0,T)},
\\ \\
- P_L\Delta {\mathbf{Z}}=\curl \left(\mathbf{R}\mathbf{W}\right),\quad \pier{\mathbf{Z}(\mathbf{a}, t)=\mathbf{0} \quad \text{for} \; (\mathbf{a}, t) \in \Gamma_a\times (0,T)}.
\end{cases}
\end{equation}
Concerning boundary conditions, note that if we don't impose homogeneous Dirichlet boundary conditions on $\nwhat{\mathbf{W}}$ and $\nwhat{\boldsymbol{\Omega}}$, then, in view of the boundary terms in \eqref{eqn:19}, we may set homogeneous Neumann boundary conditions of the form \reviewa{\[\mathbf{X}\mathbf{N}=\left(\diver \mathbf{Y}\right)\mathbf{N}=\left(\Delta \grad \twobigdot{\mathbf{W}}\right)\mathbf{N}=\mathbf{Y}\mathbf{N}=\left(\grad \grad \twobigdot{\mathbf{W}}\right)\mathbf{N}=\mathbf{0}\] and \[\boldsymbol{\Lambda}\mathbf{N}=\left(\diver \mathbf{C}\right)\mathbf{N}=\left(\Delta \grad \overbigdot{\boldsymbol{\Omega}}\right)\mathbf{N}=\mathbf{C}\mathbf{N}=\left(\grad \grad \overbigdot{\boldsymbol{\Omega}}\right)\mathbf{N}=\mathbf{0}\] on $\Gamma_a\times (0,T)$.}
\newline
We now assign general constitutive assumptions for $\boldsymbol{\Pi},\mathbf{M},\mathbf{X},\mathbf{Y},\boldsymbol{\Lambda},\mathbf{C},\boldsymbol{\Gamma}$ in order for \eqref{eqn:20} to satisfy the Clausius--Duhem dissipative equality in isothermal situations, which has the form
\abramo{
\begin{equation}
\label{eqn:21}
\frac{d\psi}{dt}+\left(\frac{d D}{d\overbigdot{C}}(\overbigdot{C}),\overbigdot{C}\right)=-p_{\text{int}}(\mathcal{D}_a,\overbigdot{C}),
\end{equation}
}%
where $\overbigdot{C}:=(\overbigdot{\mathbf{W}},\boldsymbol{\Omega},\overbigdot{\mathbf{Z}})$ is the actual velocity, $\psi$ is the free energy of the system and $D$ is the dissipation potential. We assume the following form for the free energy of the system:
\begin{align}
\label{eqn:22}
\displaystyle \notag &\psi(\mathbf{W},\mathbf{R},\mathbf{Z}):=\frac{1}{2}\|\mathbf{W}-\mathbf{I}\|^2+\nwhat{\psi}(\mathbf{W})+\frac{1}{2}\pier{\|}\grad \mathbf{W}\pier{\|}^2+\frac{1}{2}\pier{\|}\grad \mathbf{R}\pier{\|}^2\\
& \displaystyle +\pier{\int_{\mathcal{D}_a}\reviewm{k} |\curl \mathbf{Z}|}+\frac{1}{2}\pier{\|}\curl \mathbf{Z}\pier{\|}^2+\frac{\alpha_1}{2}\pier{\|}\grad \grad \mathbf{W}\pier{\|}^2,
\end{align}  
where $k\geq 0$ is a material parameter, whose meaning will be specified later. \abramohh{We observe that the particular choice for the part of the free energy depending on $\mathbf{Z}$ will induce a constitutive law for $\boldsymbol{\Gamma}+\boldsymbol{\Pi}$ representing conditional compatibility, as discussed in the Introduction and in the Remark \ref{rem:remsigma}} . Moreover, $\alpha_1 \geq 0$ is a \michhhb{physical \ap{coefficient for the second gradient contribution}}, and 
\begin{align}
\label{eqn:23}
\nwhat{\psi}(\mathbf{W}):=\int_{\mathcal{D}_a}I_{\text{SPD}_{\alpha}}(\mathbf{W}),
\end{align}
where $I_{\text{SPD}_{\alpha}}$ is the indicator function of the set
\begin{align}
\label{eqn:24}
\text{SPD}_{\alpha}:=\{\mathbf{W}\in \pier{Sym(\mathbb{R}^{3\times 3})}: \text{det}\mathbf{W}\geq \alpha^3,\; \text{tr}(\text{cof}\mathbf{W})\geq 2\alpha^2,\; \text{tr}\mathbf{W}\geq 3\alpha\}.
\end{align}
\pier{If $\alpha> 0$ the elements of  $\text{SPD}_{\alpha}$ are positive definite matrices due to the constraints on 
the determinant and the other quantities in \eqref{eqn:24}. More precisely, the tensor elements of the set 
\eqref{eqn:24} are characterized by the fact that all their eigenvalues are \abramohhb{not smaller than $\alpha$ at the same time}. 
Also, let us point out that, for
$\mathbf{W} \in L^2(\mathcal{D}_a;\mathbb{R}^{3\times 3})$, 
\begin{align*}
\nwhat{\psi}(\mathbf{W})=\int_{\mathcal{D}_a}I_{\text{SPD}_{\alpha}}(\mathbf{W}) = 
\begin{cases}
0 \quad &\hbox{if } \, \mathbf{W} \in \text{SPD}_{\alpha} \, \hbox{ a.e. in } \,\mathcal{D}_a,\\
+\infty \quad &\hbox{otherwise}. 
\end{cases}
\end{align*}
Then, the} functional \eqref{eqn:23} may be written also as
\begin{align}
\label{eqn:25}
\nwhat{\psi}(\mathbf{W})=\int_{\mathcal{D}_a}I_{S}(\mathbf{W})+\int_{\mathcal{D}_a}I_{C_{\alpha}}(\mathbf{W}),
\end{align}
where \pier{$I_{S}$} is the indicator function of the set of symmetric matrices and \pier{$I_{C_{\alpha}}$} is the indicator function of the set
\begin{align}
\label{eqn:26}
\text{C}_{\alpha}:=\{\mathbf{W}\in M(\mathbb{R}^{3\times 3}): \text{det}\mathbf{W}\geq \alpha^3,\; \text{tr}(\text{cof}\mathbf{W})\geq 2\alpha^2,\; \text{tr}\mathbf{W}\geq 3\alpha\}
\end{align}
for $\alpha> 0$. \pier{Let us introduce for future convenience the notation}
\begin{equation}
\label{eqn:22bis}
 \pier{\psi_D(\mathbf{A}):= \pier{\int_{\mathcal{D}_a} \reviewm{k}|\mathbf{A}}| +\frac{1}{2}\pier{\|}\mathbf{A}\pier{\|}^2, \quad \hbox{for all } \,\mathbf{A}\in L^2(\mathcal{D}_a;\mathbb{R}^{3\times 3}).}
\end{equation}

\begin{rem}
\label{rem:2}
The set \pier{\rm $\text{SPD}_{\alpha}$} defined in \eqref{eqn:24} is closed and convex for all $\alpha \geq 0$, hence the indicator function  \pier{\rm $I_{\text{SPD}_{\alpha}}(\cdot)$} is a convex and l.s.c. function. Moreover, the function $I_{S}(\mathbf{W})+I_{C_{\alpha}}(\mathbf{W})$ is also a convex and l.s.c. function. The proofs of these properties can be found e.g. in \cite{MM}.
\end{rem}
\abramohhb{The fact that the free energy \eqref{eqn:22} is convex implies the derivation of constitutive laws for the material that are monotone. From a mechanical point of view, roughly speaking this means that the more you push the more the material is affected by the deformation.}
\newline
\pier{Moreover, we} assume the following form for the dissipation potential of the system, \abramohhb{containing viscous contributions}:
 \begin{align}
\label{eqn:27}
& \displaystyle \notag D(\overbigdot{\mathbf{W}},\boldsymbol{\Omega}):=\frac{1}{2}\pier{\|}\overbigdot{\mathbf{W}}\pier{\|}^2+\frac{1}{2}\pier{\|}\grad \boldsymbol{\Omega}\pier{\|}^2 + \psi_A (   \boldsymbol{\Omega} ) \\
& \displaystyle  +\frac{\alpha_2}{2}\pier{\|}\grad\overbigdot{\mathbf{W}}\pier{\|}^2 +\frac{\alpha_3}{2}\pier{\|}\grad \grad\overbigdot{\mathbf{W}}\pier{\|}^2+\frac{\alpha_4}{2}\pier{\|}\grad \grad \boldsymbol{\Omega}\pier{\|}^2,
\end{align}
where $\alpha_2,\alpha_3,\alpha_4\geq0$ are \michhhb{physical parameters}\pier{,
$$ \psi_A (   \boldsymbol{\Omega} ):=\int_{\mathcal{D}_a}I_A(\boldsymbol{\Omega})$$
and $I_{A}$} is the indicator function of the set of antisymmetric matrices.
\begin{rem}
\label{rem:1}
The terms proportional to \pier{the non-negative constants $\alpha_1,\ldots,\alpha_4$} in \eqref{eqn:22} and \eqref{eqn:27} introduce higher order gradient and time derivative \michhhb{terms} in \pier{the system} dynamics, and they will be activated (\pier{i.e.,} they will be taken different from zero) only when high regularity in space and time will be required to prove existence of a solution to the system. 
\end{rem}

\noindent
Using \eqref{eqn:16}, \eqref{eqn:22} and \eqref{eqn:27} in \eqref{eqn:21}, we obtain the following constitutive assumptions
\begin{equation}
\label{eqn:29}
\displaystyle \mathbf{R}^T\boldsymbol{\Pi}=\pier{\mathbf{W}-\mathbf{I} 
+\boldsymbol{\chi}_{\alpha}+\overbigdot{\mathbf{W}}},
\end{equation}
where 
\pier{$\boldsymbol{\chi}_{\alpha} \in \partial \pier{\nwhat{\psi}(\mathbf{W})}$;}
\begin{equation}
\label{eqn:30}
\displaystyle \mathbf{M}=2\boldsymbol{\Pi}\mathbf{W}\mathbf{R}^T-2\boldsymbol{\mathit{S}},
\end{equation}
where $\boldsymbol{\mathit{S}}\in \pier{{}\partial \psi_{A}(\boldsymbol{\Omega})}$;
\begin{equation}
\label{eqn:31}
\displaystyle \boldsymbol{\Sigma}:=-(\boldsymbol{\Gamma}+\boldsymbol{\Pi})\in \partial \psi_D(\curl \mathbf{Z})
= \pier{%
\begin{cases}
\displaystyle k\frac{\curl \mathbf{Z}}{|\curl\mathbf{Z}|}+\curl\mathbf{Z}\quad \hbox{if }\, |\curl\mathbf{Z}| \not= 0,\\
\hbox{any } \, \mathbf{M}_D, \, \hbox{ with } |\mathbf{M}_D|\leq k, \, \hbox{ if }\, |\curl\mathbf{Z}| = 0;
\end{cases}
}
\end{equation}
\begin{equation}
\label{eqn:32}
\begin{cases}
 \displaystyle \mathbf{X}=\grad \mathbf{W}+\alpha_2\grad \overbigdot{\mathbf{W}};\\
 \displaystyle  \mathbf{Y}=\alpha_1\grad\grad \mathbf{W}+\alpha_3\grad\grad \overbigdot{\mathbf{W}};\\
 \displaystyle \boldsymbol{\Lambda}=(\grad \mathbf{R})\mathbf{R}^T+\grad \boldsymbol{\Omega};\\
\displaystyle \mathbf{C}=\alpha_4\grad\grad \boldsymbol{\Omega}.
\end{cases}
\end{equation}

\begin{rem}
\label{rem:remsigma}
We observe that $\partial \psi_D$ is a maximal monotone operator \pier{in $L^2(\mathcal{D}_a;\mathbb{R}^{3\times 3})$} and \pier{\eqref{eqn:31} entails that}
\begin{equation}
 \label{eqn:remsigma}
 \curl \mathbf{Z}=\boldsymbol{0} \quad \text{if \pier{and only if}}\;\; |\boldsymbol{\Sigma}|\leq k.
\end{equation}
Hence, if the norm of the reaction term $\boldsymbol{\Gamma+\Pi}$ in \eqref{eqn:20} is lower or equal than \pier{the} threshold $k$, the dislocation tensor $\curl \mathbf{Z}$ in \eqref{eqn:3} is the null tensor, and the motion is compatible.
Moreover, we observe that the constitutive assumptions \eqref{eqn:29}--\eqref{eqn:32} \pier{comply with} the principle of objectivity, \pier{that is, the property
\[
p_{\text{int}}(\mathcal{D}_a,C_{\text{rigid}})=0
\] 
is satisfied for the} rigid virtual velocities, \pier{i.e.,} for $\nwhat{\mathbf{W}}=\boldsymbol{0}$, $\nwhat{\boldsymbol{\Omega}}=\mathbf{A}$, \ap{$\grad \mathbf{V}=\nwhat{\boldsymbol{\Omega}}\grad \boldsymbol{\Phi}$ and $\curl \nwhat{\mathbf{Z}}=\nwhat{\boldsymbol{\Omega}}\curl \mathbf{Z}$}, for any \ap{spatially} constant \pier{tensor} $\mathbf{A}\in \text{Skew}(\mathbb{R}^{3\times 3})$. Indeed, \pier{formally} we have that
\begin{align*}
&p_{\text{int}}(\mathcal{D}_a,C_{\text{rigid}})=-\int_{\mathcal{D}_a}\mathbf{S}:\mathbf{A}+\int_{\mathcal{D}_a} \left(\boldsymbol{\Gamma}+ \boldsymbol{\Pi}\right):  \curl \nwhat{\mathbf{Z}}
=-\int_{\mathcal{D}_a}\pier{\partial \psi_D(\curl \mathbf{Z})}:\mathbf{A} \curl{\mathbf{Z}}\\
&=-\int_{\mathcal{D}_a}\pier{\partial \psi_D(\curl \mathbf{Z})}\left(\curl{\mathbf{Z}}\right)^T:\mathbf{A}=0,
\end{align*}
\reviewa{since the last integrand is the scalar product of a symmetric tensor with a skew-symmetric tensor.}
We remark that \ap{even in presence of incompatibility, the principle of objectivity is satisfied.}
\end{rem}

\pier{\begin{rem}
\label{rem:pier1}
The subdifferential \pier{$ \partial \pier{\nwhat{\psi}}$}  is a maximal monotone operator \pier{in $L^2(\mathcal{D}_a;\mathbb{R}^{3\times 3})$} as well, and the inclusion $\boldsymbol{\chi}_{\alpha} \in \partial \pier{\nwhat{\psi}(\mathbf{W})}$ means that
\begin{align}
\label{pier1}
&\mathbf{W} \, \hbox{ belongs to the domain of } \, \partial \nwhat{\psi} \, \hbox{ and} \notag \\
& \int_{\mathcal{D}_a} \boldsymbol{\chi}_{\alpha}: ( \nwhat{\mathbf{W}} - \mathbf{W}) + \nwhat{\psi} (\mathbf{W}) \leq \nwhat{\psi} (\nwhat{\mathbf{W}} )
\quad \hbox{for all }  \, \nwhat{\mathbf{W}} \in L^2(\mathcal{D}_a;\mathbb{R}^{3\times 3}).
\end{align}
In view of \eqref{eqn:23}--\eqref{eqn:24} and \eqref{eqn:25}--\eqref{eqn:26}, it is not difficult to show that
\eqref{pier1} can be equivalently rewritten as
\begin{align}
\label{pier2}
&\mathbf{W} \in L^2(\mathcal{D}_a;\mathbb{R}^{3\times 3}), \quad \mathbf{W} \in \pier{Sym(\mathbb{R}^{3\times 3})} \, \hbox{ almost everywhere in }\, \mathcal{D}_a,  \hbox{ and} \notag \\
& \int_{\mathcal{D}_a} \boldsymbol{\chi}_{\alpha}: ( \nwhat{\mathbf{W}} - \mathbf{W}) + 
\int_{\mathcal{D}_a} I_{C_\alpha} (\mathbf{W})  \leq \int_{\mathcal{D}_a} I_{C_\alpha} (\nwhat{\mathbf{W}} )
\notag\\
&\quad \hbox{for all symmetric matrices }  \, \nwhat{\mathbf{W}} \in L^2(\mathcal{D}_a;\mathbb{R}^{3\times 3}).
\end{align}
Then, it becomes clear that, setting 
$$ \psi_{C_\alpha} ( \nwhat{\mathbf{W}} ):=\int_{\mathcal{D}_a}I_{C_\alpha}(\nwhat{\mathbf{W}}) , $$
the inclusion $\boldsymbol{\chi}_{\alpha} \in \partial \pier{\nwhat{\psi}(\mathbf{W})}$ can be formulated as  
\begin{align}
\label{pier3}
\mathbf{W} \in L^2(\mathcal{D}_a;\mathbb{R}^{3\times 3}), \quad \mathbf{W} \in \pier{Sym(\mathbb{R}^{3\times 3})} \, \hbox{ a.e. in }\, \mathcal{D}_a,  \hbox{ and } \, \boldsymbol{\chi}_{\alpha}\in \partial \psi_{C_\alpha} (\mathbf{W}).
\end{align}
\end{rem}}%

\noindent
Inserting \eqref{eqn:29}--\eqref{eqn:32} in \eqref{eqn:20} we finally obtain
\begin{equation}
\label{eqn:33}
\begin{cases}
\displaystyle \mathbf{R}^T\grad \left(\mathcal{G}_L\ast\frac{d\mathbf{U}}{dt}\right)+\twobigdot{\mathbf{W}}+\pier{{}\reviewa{\diver \Delta \grad}{}}\twobigdot{\mathbf{W}}+\mathbf{R}^T\curl\left(\mathcal{G}_{L,\diver}\ast\left(\curl 
\boldsymbol{\Sigma}\right)\right)+\mathbf{W}-\mathbf{I}\\
\displaystyle \pier{{}+\boldsymbol{\chi}_{\alpha}+\overbigdot{\mathbf{W}}}-\Delta \mathbf{W}-\alpha_2\Delta \overbigdot{\mathbf{W}}+\alpha_1\pier{{}\reviewa{\diver \Delta \grad}{}} \mathbf{W}+\alpha_3\pier{{}\reviewa{\diver \Delta \grad}{}} \overbigdot{\mathbf{W}}=\mathbf{W}_{\text{ext}},\\ \\
\displaystyle \grad \left(\mathcal{G}_L\ast\frac{d\mathbf{U}}{dt}\right)\mathbf{W}\mathbf{R}^T+\overbigdot{\boldsymbol{\Omega}}+\pier{{}\reviewa{\diver \Delta \grad}{}}\overbigdot{\boldsymbol{\Omega}}+\curl\left(\mathcal{G}_{L,\diver}\ast\left(\curl \boldsymbol{\Sigma}\right)\right)\mathbf{W}\mathbf{R}^T \\
\displaystyle{}+\boldsymbol{\mathit{S}}-\frac{1}{2}\diver \left((\grad \mathbf{R})\mathbf{R}^T\right)-\frac{1}{2}\Delta\boldsymbol{\Omega}+\frac{1}{2}\alpha_4\pier{{}\reviewa{\diver \Delta \grad}{}}\boldsymbol{\Omega}=\boldsymbol{\Omega}_{\text{ext}},
\\ \\
\displaystyle \pier{\boldsymbol{\chi}_{\alpha} \in \partial \pier{\nwhat{\psi}(\mathbf{W})}, \quad  
\boldsymbol{\mathit{S}}\in \partial \psi_{A}(\boldsymbol{\Omega}), \quad  
\boldsymbol{\Sigma}\in \partial \psi_D(\curl \mathbf{Z}),}\\ \\
\overbigdot{\mathbf{R}}=\boldsymbol{\Omega}\mathbf{R},\\ \\
\Delta \pier{\mathbf{\Phi}}=\diver\left(\mathbf{R}\mathbf{W}\right),\\ \\
-P_L\Delta {\mathbf{Z}}=\curl \left(\mathbf{R}\mathbf{W}\right),
\end{cases}
\end{equation}
valid in $\mathcal{D}_{aT}$, with boundary conditions
\begin{align}
\label{eqn:34}
& \displaystyle \notag \mathbf{W}=\mathbf{I},\; \overbigdot{\mathbf{W}}=\boldsymbol{0},\; \reviewa{\grad \mathbf{W}=\mathbf{0}} \quad \text{on} \; \Gamma_a\times (0,T),\\
& \displaystyle \notag  \mathbf{R}=\mathbf{I},\;\boldsymbol{\Omega}=\boldsymbol{0},\; \reviewa{\grad \boldsymbol{\Omega}=\mathbf{0}} \quad \text{on} \; \Gamma_a\times (0,T),\\ 
& \displaystyle \pier{\mathbf{\Phi}(\mathbf{a}, t) =\mathbf{a}, \; \mathbf{Z}(\mathbf{a}, t)=\mathbf{0} \quad \text{for} \; (\mathbf{a}, t) \in \Gamma_a\times (0,T)},
\end{align}
 and initial conditions
\begin{equation}
\label{eqn:35}
\mathbf{W}(\mathbf{a},0)=\mathbf{I}, \; \overbigdot{\mathbf{W}}(\mathbf{a},0)=\mathbf{0}, \; \mathbf{R}(\mathbf{a},0)=\mathbf{I}, \; \boldsymbol{\Omega}(\mathbf{a},0)=\mathbf{0}, \; \mathbf{Z}(\mathbf{a},0)=\mathbf{0}
\quad \pier{\text{for} \; \mathbf{a} \in \mathcal{D}_a}.
\end{equation}

{
\subsection{\reviewa{An Example: the reaction to the compatibility condition}}
We consider the case in which the evolution is simply given by 
\begin{gather*}
\boldsymbol{\Phi} (\mathbf{a},t)=\mathbf{a}, \,\mathbf{U}(\mathbf{a},t)=\mathbf{0}, \, \mathbf{R}(\mathbf{a},t)\mathbf{=I}, \, \mathbf{ W}(\mathbf{a},t)\mathbf{=I,\ }, \, \mathbf{Z}(\mathbf{a},t)=\mathbf{0.}
\end{gather*}%
This yields a solution of equations \eqref{eqn:33}$_1$, \eqref{eqn:33}$_2$, \eqref{eqn:33}$_3$ and \eqref{eqn:33}$_6$ if 
\begin{gather*}
\mathbf{\mathbf{W}}_{\text{ext}}={\rm Sym} (\curl\left\{ \mathcal{G}_{L,%
\diver}\ast \curl \mathbf{\Sigma }\right\} ) , \\
\boldsymbol{\Omega}_{\text{ext}}={\rm Skew} (\curl \left\{ \mathcal{%
G}_{L,\diver}\ast \curl \mathbf{\Sigma }\right\} ),
\end{gather*}%
are given by the internal \ap{force} \reviewa{$\mathbf{\Sigma }$ satisfying the property that}%
\begin{gather*}
\left\vert \mathbf{\Sigma }(\mathbf{a},t)\right\vert \leq k,
\end{gather*}%
which is assumed to be known.
The external actions do not work and do not result in motion. They have no
macroscopic effect. But they produce dislocations which modify the internal
stress state $\mathbf{\Sigma }(\mathbf{a},t)$. To produce a motion the external
actions have to be increased.
{\reviewm{
\newline
 Let matrix $\mathbf{\hat{N}\in M(\mathbb{R}^{3\times 3})}$ be giving independent virtual stretch and angular velocities%
\begin{equation*}
\mathbf{\hat{W}=}{\rm Sym}\left\{ \mathbf{\hat{N}}\right\} ,~\mathbf{\hat{\Omega}}%
 \mathbf{=}{\rm Skew}\left\{ \mathbf{\hat{N}}\right\}
,
\end{equation*}
 satisfying the related boundary conditions of (\eqref{eqn:34})
\begin{equation*}
 \displaystyle \notag  \hat{\mathbf{W}}=\boldsymbol{0},\; \hat{ \boldsymbol{\Omega}}=\boldsymbol{0},\; \reviewa{\grad \hat{\mathbf{W}}=\grad\hat{ \boldsymbol{\Omega}}=\boldsymbol{0}} \quad \text{on} \; \Gamma_a.
\end{equation*}
In this example, the virtual power of the external forces is equal to the opposite of the virtual power of the internal forces
\begin{gather*}
\forall \mathbf{\hat{N},} \\
\int_{\mathcal{D}_{a}}\left(\mathbf{\mathbf{W}}^{ext}:\mathbf{\hat{W}}+\mathbf{%
\mathbf{\Omega }}^{ext}:\mathbf{\hat{\Omega}}\right)=\int_{\mathcal{D}_{a}}(%
\mathbf{\mathbf{W}}^{ext}+\mathbf{\mathbf{\Omega }}^{ext}):\mathbf{\hat{N}}
\\
=\int_{\mathcal{D}_{a}}\curl\left\{ \reviewa{\mathcal{G}_{L,\diver}}\ast \curl\mathbf{%
\Sigma }\right\} :\mathbf{\hat{N}}.
\end{gather*}}
}
\newline
Note that we have%
\begin{equation*}
\mathbf{\mathbf{W}}_{\text{ext}}+\boldsymbol{\Omega}_{\text{ext}}=\curl %
\left\{ \mathcal{G}_{L,\diver}\ast \curl \mathbf{\Sigma }\right\}.
\end{equation*}%
{\reviewm{For}} a virtual velocity $\mathbf{V}$ with stretch \pierhhb{and angular velocities
$${ \rm Sym}\left\{ \grad\mathbf{V}\right\} \hbox{ and }\ 
{\rm Skew}\left\{ \grad\mathbf{V}\right\},$$
respectively, \reviewa{with the proper boundary conditions,} we have that} 
\begin{gather*}
\int_{\mathcal{D}_{a}}\left(\mathbf{\mathbf{W}}_{\text{ext}}:{\rm Sym}(\grad%
\mathbf{V})+\boldsymbol{\Omega}_{\text{ext}}:{\rm Skew}(%
\grad \mathbf{V})\right) \\
=\int_{\mathcal{D}_{a}}\curl \left\{ \mathcal{G}_{L,\diver}\ast 
\curl \mathbf{\Sigma }\right\} :\grad \mathbf{V}=0.
\end{gather*}%
Then $\curl \left\{ \mathcal{G}_{L,\diver}\ast \curl \mathbf{%
\Sigma }\right\} $ is a reaction: a reaction to the compatibility conditions. Its power is null for any virtual velocity  $\mathbf{\hat{N}}$ which is a 
gradient, i.e., for virtual velocities
$\mathbf{\hat{W}}$,
 $\mathbf{\hat{\Omega}}$  
which satisfy the compatibility conditions.
 Its power may be non null for virtual velocities which do not satisfy the compatibility conditions. 
\reviewm{It has the usual property of a reaction to a kinematic constraint:
it is normal to the linear set of the compatible virtual matrices velocities $\mathbf{\hat{N}}$ which are the gradients of virtual velocities. And it does not work in the actual evolution.
} Inside $\mathcal{D}_{a}$ the {\reviewm{actual}} densities of power due to the evolution of the
dislocations are not null. But their total sum is null.
\section{Quasi-stationary case}
\label{sec:analysis}
In this section we study the existence of solutions to \eqref{eqn:33} in the quasi-stationary case, \pier{i.e.,} considering $p_{acc}(\mathcal{D}_a,C)=0$ for all $C\in \mathcal{C}$ and thus neglecting the inertia terms in \eqref{eqn:33}$_1$ and \eqref{eqn:33}$_2$. We \pier{deal with} the case $k>0$ and \pier{let $\alpha_1,\dots,\alpha_4=0$, then we} study the existence and regularity of a global in time weak solution, \abramohhb{which will be proved to be also a strong solution}. 
\pier{Consider} the following reduced version of system~\eqref{eqn:33}:
\begin{equation}
\label{eqn:36}
\begin{cases}
\displaystyle \mathbf{R}^T\curl\left(\mathcal{G}_{L,\diver}\ast\left(\curl \boldsymbol{\Sigma}\right)\right)+ \mathbf{W}-\mathbf{I}\pier{{}+\boldsymbol{\chi}_{\alpha}}+\overbigdot{\mathbf{W}}-\Delta \mathbf{W}
\pco{{}={}}  \mathbf{W}_{\text{ext}}(\mathbf{W},t),\\ \\
\displaystyle \curl\left(\mathcal{G}_{L,\diver}\ast\left(\curl 
\boldsymbol{\Sigma}\right)\right)\mathbf{W}\mathbf{R}^T+\boldsymbol{\mathit{S}}-\frac{1}{2}\diver \left((\grad \mathbf{R})\mathbf{R}^T\right)-\frac{1}{2}\Delta\boldsymbol{\Omega}=\abramo{\boldsymbol{\Omega}_{\text{ext}}(\mathbf{R},t)},
\\ \\
\displaystyle \pier{\boldsymbol{\chi}_{\alpha} \in \partial \pier{\nwhat{\psi}(\mathbf{W})}, \quad  
\boldsymbol{\mathit{S}}\in \partial \psi_{A}(\boldsymbol{\Omega}), \quad  
\boldsymbol{\Sigma}\in \partial \psi_D(\curl \mathbf{Z}),}\\ \\
\Delta \pier{\mathbf{\Phi}}=\diver\left(\mathbf{R}\mathbf{W}\right),\\ \\
- P_L\Delta {\mathbf{Z}}=\curl \left(\mathbf{R}\mathbf{W}\right),
\end{cases}
\end{equation}
valid in $\mathcal{D}_{aT}$, with boundary conditions
\begin{equation}
\label{eqn:37}
\begin{cases}
\displaystyle \mathbf{W}\pier{{}={}} \mathbf{R}=\mathbf{I},\; \boldsymbol{\Omega}=\boldsymbol{0} \quad \text{on} \; \Gamma_a\times (0,T),\\ 
\displaystyle  \pier{\mathbf{\Phi}(\mathbf{a}, t) =\mathbf{a}, \; \mathbf{Z}(\mathbf{a}, t)=\mathbf{0} \quad \text{for} \; (\mathbf{a}, t) \in \Gamma_a\times (0,T)},
\end{cases}
\end{equation}
 and initial conditions
\begin{equation}
\label{eqn:38}
\mathbf{W}(\mathbf{a},0)=\mathbf{W}_0(\mathbf{a}), \; \mathbf{R}(\mathbf{a},0)=\mathbf{R}_0(\mathbf{a})
 \pier{ \quad \text{for} \; \mathbf{a}  \in \mathcal{D}_a},
\end{equation}
\pier{where} $\mathbf{W}_0 \pier{{}={}} \mathbf{R}_0=\mathbf{I}$ on $\Gamma_a$.
\noindent
For simplicity, \pier in the following 
we will take} $\mathbf{R}_0=\mathbf{I}$.
We observe that, given $\boldsymbol{\Omega}\in \mathcal{A}$, the differential equation \eqref{eqn:36}$_4$ and the initial condition in \eqref{eqn:38}, with $\mathbf{R}_0=\mathbf{I}$, uniquely \pier{define} a rotation tensor 
\[
\mathbf{R}(\mathbf{a},t)=\mathrm{e}^{\int_0^t\boldsymbol{\Omega}(\mathbf{a},s)ds}\quad \pier{\text{for} \; (\mathbf{a}, t) \in \mathcal{D}_a\times (0,T)}.
\]
Since $\boldsymbol{\Omega}\in \mathcal{A}$, we have that $\mathbf{R}:\mathbf{R}=\mathrm{e}^{\int_0^t\boldsymbol{\Omega}(\mathbf{a},s)ds}\mathrm{e}^{-\int_0^t\boldsymbol{\Omega}(\mathbf{a},s)ds}:\mathbf{I}=3$, hence
\begin{equation}
\label{eqn:39}
\mathbf{R}
\in L^{\infty}(\mathcal{D}_{aT},\mathbb{R}^{3\times 3}).
\end{equation}
We introduce the variable $\boldsymbol{\Theta}(\mathbf{a},t):=\int_0^t\boldsymbol{\Omega}(\mathbf{a},s)ds$, 
\pier{$(\mathbf{a}, t) \in \mathcal{D}_a\times (0,T),$}
and rewrite \pier{the system}~\eqref{eqn:36} as
\begin{equation}
\label{eqn:48}
\begin{cases}
\displaystyle \mathrm{e}^{-\boldsymbol{\Theta}}\curl\left(\mathcal{G}_{L,\diver}\ast
\left(\curl \boldsymbol{\Sigma}\right)\right)+\mathbf{W}-\mathbf{I}\pier{{}+\boldsymbol{\chi}_{\alpha}}+\overbigdot{\mathbf{W}}-\Delta \mathbf{W}\pco{{}={}} \mathbf{W}_{\text{ext}}(\mathbf{W},t),\\ \\
\displaystyle \curl\left(\mathcal{G}_{L,\diver}\ast
\left(\curl \boldsymbol{\Sigma}\right)\right)\mathbf{W}\mathrm{e}^{-\boldsymbol{\Theta}} +\boldsymbol{\mathit{S}}-\frac{1}{2}\Delta \boldsymbol{\Theta}-\frac{1}{2}\Delta \overbigdot{\boldsymbol{\Theta}}=\abramo{\boldsymbol{\Omega}_{\text{ext}}(\boldsymbol{\Theta},t)},
\\ \\
\displaystyle \pier{\boldsymbol{\chi}_{\alpha} \in \partial \pier{\nwhat{\psi}(\mathbf{W})}, \quad  
\boldsymbol{\mathit{S}}\in \partial \psi_{A}(\abramonew{\overbigdot{\boldsymbol{\Theta}}}), \quad  
\boldsymbol{\Sigma}\in \partial \psi_D(\curl \mathbf{Z}),}\\ \\
\displaystyle \Delta \pier{\mathbf{\Phi}}=\diver\left(\mathrm{e}^{\boldsymbol{\Theta}}\mathbf{W}\right),\\ \\
\displaystyle - P_L\Delta {\mathbf{Z}}=\curl \left(\mathrm{e}^{\boldsymbol{\Theta}}\mathbf{W}\right),
\end{cases}
\end{equation}
with boundary conditions
\begin{equation}
\label{eqn:49}
\begin{cases}
\displaystyle \mathbf{W}=\mathbf{I},\;\boldsymbol{\Theta}=\overbigdot{ \boldsymbol{\Theta}}=\boldsymbol{0} \quad \text{on} \; \Gamma_a\times (0,T),\\ 
\displaystyle  \pier{\mathbf{\Phi}(\mathbf{a}, t) =\mathbf{a}, \; \mathbf{Z}(\mathbf{a}, t)=\mathbf{0} \quad \text{for} \; (\mathbf{a}, t) \in \Gamma_a\times (0,T)},
\end{cases}
\end{equation}
 and initial conditions
\begin{equation}
\label{eqn:50}
\mathbf{W}(\mathbf{a},0)=\mathbf{W}_0\pier{(\mathbf{a})}, \; \boldsymbol{\Theta}(\mathbf{a},0)=\boldsymbol{0} \pier{ \quad \text{for} \; \mathbf{a}  \in \mathcal{D}_a}.
\end{equation}
We observe that an initial condition for $\mathbf{Z}$ can be defined by assuming that \eqref{eqn:48}$_5$ is valid for $t=0$, \pier{i.e.,}
\begin{equation}
\label{eqn:50z}
\mathbf{Z}(\mathbf{a},0)= (\mathcal{G}_{L,\diver}\ast \curl \mathbf{W}_0 ) \pier{(\mathbf{a})}\pier{ \quad \text{for} \; \mathbf{a}  \in \mathcal{D}_a}.
\end{equation}
In the case \pier{$\mathbf{W}_0=\mathbf{I}$, then we have} $\mathbf{Z}(\mathbf{a},0)=\boldsymbol{0}$ $\pier{ \ \text{for} \; \mathbf{a}  \in \mathcal{D}_a}$. 

\pier{Note that the equations \eqref{eqn:48}$_1$ and  \eqref{eqn:48}$_2$ are coupled. Also, since $\psi_A$ is defined as the integral of the indicator function $I_{A}$  of the set of antisymmetric matrices, and 
 $\boldsymbol{\mathit{S}}$ should satisfy
$\boldsymbol{\mathit{S}}\in \partial \psi_{A}(\boldsymbol{\Omega})$, that is, $\boldsymbol{\mathit{S}}\in \partial I_{A}(\boldsymbol{\Omega})$ a.e. in $ \mathcal{D}_a$, then $\boldsymbol{\mathit{S}}$ can be recovered a posteriori in terms of the symmetric part of \eqref{eqn:48}$_2$, that~is, 
\begin{equation}
\label{eqn:52}
\boldsymbol{\mathit{S}}= - \rm{Sym}\left(\curl\left(\mathcal{G}_{L,\diver}\ast
\left(\curl \boldsymbol{\Sigma}\right)\right)\mathbf{W}\mathrm{e}^{-\boldsymbol{\Theta}}\right).
\end{equation}}%
\abramonew{
\begin{rem}
 \label{rem:symskew}
 In the case in which $\mathbf{W} \in \mathcal{S}$ and $\boldsymbol{\Theta}\in \mathcal{A}$, the system \eqref{eqn:48} becomes
\begin{equation}
\label{eqn:48symskew}
\begin{cases}
\displaystyle \text{Sym}\left(\mathrm{e}^{-\boldsymbol{\Theta}}\curl\left(\mathcal{G}_{L,\diver}\ast
\left(\curl \boldsymbol{\Sigma}\right)\right)\right)+\mathbf{W}-\mathbf{I}\pier{{}+\boldsymbol{\chi}_{\alpha}}+\overbigdot{\mathbf{W}}-\Delta \mathbf{W}\pco{{}={}} \mathbf{W}_{\text{ext}}(\mathbf{W},t),\\ \\
\displaystyle \text{Skew}\left(\curl\left(\mathcal{G}_{L,\diver}\ast
\left(\curl \boldsymbol{\Sigma}\right)\right)\mathbf{W}\mathrm{e}^{-\boldsymbol{\Theta}}\right) -\frac{1}{2}\Delta \boldsymbol{\Theta}-\frac{1}{2}\Delta \overbigdot{\boldsymbol{\Theta}}=\abramo{\boldsymbol{\Omega}_{\text{ext}}(\boldsymbol{\Theta},t)},
\\ \\
\displaystyle \pier{\boldsymbol{\chi}_{\alpha} \in \partial \pier{\nwhat{\psi}(\mathbf{W})}, \quad  
\boldsymbol{\Sigma}\in \partial \psi_D(\curl \mathbf{Z}),}\\ \\
\displaystyle \Delta \pier{\mathbf{\Phi}}=\diver\left(\mathrm{e}^{\boldsymbol{\Theta}}\mathbf{W}\right),\\ \\
\displaystyle - P_L\Delta {\mathbf{Z}}=\curl \left(\mathrm{e}^{\boldsymbol{\Theta}}\mathbf{W}\right),
\end{cases}
\end{equation}
where the inclusion $\boldsymbol{\chi}_{\alpha} \in \partial \nwhat{\psi}(\mathbf{W})$ is expressed as in \eqref{pier2}, with boundary conditions
\begin{equation}
\label{eqn:49symskew}
\begin{cases}
\displaystyle \mathbf{W}=\mathbf{I},\;\boldsymbol{\Theta}=\overbigdot{ \boldsymbol{\Theta}}=\boldsymbol{0} \quad \text{on} \; \Gamma_a\times (0,T),\\ 
\displaystyle  \pier{\mathbf{\Phi}(\mathbf{a}, t) =\mathbf{a}, \; \mathbf{Z}(\mathbf{a}, t)=\mathbf{0} \quad \text{for} \; (\mathbf{a}, t) \in \Gamma_a\times (0,T)},
\end{cases}
\end{equation}
 and initial conditions
\begin{equation}
\label{eqn:50symskew}
\mathbf{W}(\mathbf{a},0)=\mathbf{W}_0\pier{(\mathbf{a})}, \; \boldsymbol{\Theta}(\mathbf{a},0)=\boldsymbol{0} \pier{ \quad \text{for} \; \mathbf{a}  \in \mathcal{D}_a}.
\end{equation}
\end{rem}
}
We state now the main theorem of the present paper. We start by introducing the following assumptions on the data:
\begin{itemize}
  \item[\textbf{A1}:] $\mathcal{D}_a \subset \mathbb{R}^3$ is a bounded domain and the boundary $\Gamma_a$ is of class $C^3$;
  \item[\textbf{A2}:] The initial datum \pier{has} the regularity $\mathbf{W}_0\in H^1(\mathcal{D}_a;\mathbb{R}^{3\times 3})\cap \mathcal{S}$, with $\mathbf{W}_0\in \pier{\rm \hbox{SPD}}_{\alpha}$ \pier{almost everywhere in $\mathcal{D}_a$ for a given $\alpha>0$, \pier{and with $\mathbf{W}_0\pier{{}={}}\mathbf{I}$ on $ \Gamma_a\times (0,T)$}};
  \item[\textbf{A3}:] The forcing term \pier{$\pier{\mathbf{W}_{\text{ext}}}: H^1(\mathcal{D}_a;Sym(\mathbb{R}^{3\times 3}))\times (0,T)\rightarrow L^2(\mathcal{D}_a;Sym(\mathbb{R}^{3\times 3}))$ is measurable in $t\in (0,T)$ and Lipschitz continuous in $\mathbf{W} \in H^1(\mathcal{D}_a;Sym(\mathbb{R}^{3\times 3}))$, \pier{and it satisfies} $\pier{\mathbf{W}_{\text{ext}}}(\boldsymbol{0}, t)=\boldsymbol{0}$ for all $t\in [0,T]$ and
  \[
   \pier{\|}\pier{\mathbf{W}_{\text{ext}}}(\mathbf{W}_1,t)-\pier{\mathbf{W}_{\text{ext}}}(\mathbf{W}_2,t)\pier{\|}_{L^2(\mathcal{D}_a;Sym(\mathbb{R}^{3\times 3}))}\leq L\pier{\|}\mathbf{W}_1-\mathbf{W}_2\pier{\|}_{H^1(\mathcal{D}_a;Sym(\mathbb{R}^{3\times 3}))},
  \]
  for a.e. $t\in (0,T)$, for all $\mathbf{W}_1,\mathbf{W}_2\in H^1(\mathcal{D}_a;Sym(\mathbb{R}^{3\times 3}))$ and for some $L\in \mathbb{R}$.} 
  \abramo{Similarly, the forcing term ${\boldsymbol{\Omega}_{\text{ext}}}: H^1(\mathcal{D}_a;Skew(\mathbb{R}^{3\times 3}))\times (0,T)\rightarrow L^2(\mathcal{D}_a;Skew(\mathbb{R}^{3\times 3}))$ is measurable in $t\in (0,T)$ and Lipschitz continuous in $\boldsymbol{\Theta} \in H^1(\mathcal{D}_a;Skew(\mathbb{R}^{3\times 3}))$, {and it satisfies} ${\boldsymbol{\Omega}_{\text{ext}}}(\boldsymbol{0}, t)=\boldsymbol{0}$ for all $t\in [0,T]$ and
  \[
   \pier{\|}\pier{\boldsymbol{\Omega}_{\text{ext}}}(\boldsymbol{\Theta}_1,t)-\pier{\boldsymbol{\Omega}_{\text{ext}}}(\boldsymbol{\Theta}_2,t)\pier{\|}_{L^2(\mathcal{D}_a;Skew(\mathbb{R}^{3\times 3}))}\leq G\pier{\|}\boldsymbol{\Theta}_1-\boldsymbol{\Theta}_2\pier{\|}_{H^1(\mathcal{D}_a;Skew(\mathbb{R}^{3\times 3}))},
  \]
  for a.e. $t\in (0,T)$, for all $\boldsymbol{\Theta}_1,\boldsymbol{\Theta}_2\in H^1(\mathcal{D}_a;Skew(\mathbb{R}^{3\times 3}))$ and for some $G\in \mathbb{R}$.}  
  \end{itemize}
\begin{thm}
 \label{thm:1}
 Let assumptions \textbf{A1}-\textbf{A3} be satisfied.  Then, for any $T>0$ there is a \pco{sextuplet}
 $(\mathbf{W},\boldsymbol{\Theta},\pier{{}\boldsymbol{\chi}_{\alpha}{}}, \boldsymbol{\Sigma},\pier{\mathbf{\Phi}},\mathbf{Z})$, with
\begin{align}
\label{eqn:57thm}
\notag &\mathbf{W}\in L^{\infty}(0,T;H^1(\mathcal{D}_a;Sym(\mathbb{R}^{3\times 3})))\\
&\qquad {}\cap H^1(0,T;L^2(\mathcal{D}_a;Sym(\mathbb{R}^{3\times 3})))\cap L^2(0,T;H^2(\mathcal{D}_a,Sym(\mathbb{R}^{3\times 3}))),
\end{align}
and \pier{\rm $\mathbf{W}(\mathbf{a},t)\in \text{SPD}_{\alpha}$} for a.e. $(\mathbf{a},t)\in \mathcal{D}_{aT}$,  
\begin{equation}
\label{eqn:58thm}
\boldsymbol{\Theta}\in H^1(0,T;H^2(\mathcal{D}_a,Skew(\mathbb{R}^{3\times 3}))),
\end{equation}
\begin{equation}
 \label{eqn:58tristhm}
 \pier{{}\boldsymbol{\chi}_{\alpha}{}}\in L^{2}(0,T;L^2(\mathcal{D}_a;\mathbb{R}^{3\times 3})),
\end{equation}
\begin{equation}
 \label{eqn:58bisthm}
 \boldsymbol{\Sigma}\in L^{\infty}(0,T;L^2(\mathcal{D}_a;\mathbb{R}^{3\times 3})),
\end{equation}
\begin{equation}
\label{eqn:59thm}
\pier{\mathbf{\Phi}\in L^{\infty}(0,T;H^2(\mathcal{D}_a,\mathbb{R}^{3})\cap H^1(\mathcal{D}_a;\mathbb{R}^{3}))\cap L^{2}(0,T;H^3(\mathcal{D}_a;\mathbb{R}^{3})),}
\end{equation}
\begin{equation}
\label{eqn:60thm}
\mathbf{Z}\in L^{\infty}(0,T;H^2(\mathcal{D}_a,\mathbb{R}^{3\times 3})\cap H_{0,\diver}^1(\mathcal{D}_a,\mathbb{R}^{3\times 3}))\cap L^{2}(0,T;H^3(\mathcal{D}_a;\mathbb{R}^{3})),
\end{equation}
which solves \pier{the system \abramonew{\eqref{eqn:48symskew}--\eqref{eqn:50symskew}} with equations and conditions satisfied almost everywhere}. 
\end{thm}

\begin{pf}
Let us introduce the finite dimensional spaces which will be used to formulate the Galerkin ansatz to approximate the solutions of \pier{the system} \eqref{eqn:48}--\eqref{eqn:50}. Let $\{\xi_i\}_{i\in \mathbb{N}}$ be the eigenfunctions of the Laplace operator with homogeneous Dirichlet boundary conditions, \pier{i.e.,}
\[
-\Delta \xi_i=\gamma_i \xi_i \quad \text{in} \; \mathcal{D}_a, \quad \xi_i =0 \quad \text{on} \; \Gamma_a,
\]
with \pier{$0<\gamma_0\leq \gamma_1 \leq \dots \leq \gamma_m\to \infty$}. The sequence $\{\xi_i\}_{i\in \mathbb{N}}$ can be chosen as an orthonormal basis in $L^2(\mathcal{D}_a)$ and an orthogonal basis in $H^1(\mathcal{D}_a)$, and, thanks to Assumption~\pier{\textbf{A1}, $\{\xi_i\}_{i\in \mathbb{N}}\subset \pier{H^2}(\mathcal{D}_a)$.} 

We then introduce the functions $\{\mathbf{S}_{6k+i+j+\pier{n_i}}\}_{k\in \mathbb{N}; i,j=0,\dots,2; j\geq i}$ defined by
\[
 \mathbf{S}_{6k+i+j+(i>0)}:=\xi_k\left(\mathbf{e}_i\otimes \mathbf{e}_j+\mathbf{e}_j\otimes \mathbf{e}_i\right),
\]
where $\mathbf{e}_i, i=0, \dots, 2$ are the elements of the canonical basis of $\mathbb{R}^3$, and $\pier{n_i}$ is $0$ when $i=0$ or $1$ when $i>0$. We observe that, given $k\in \mathbb{N}$, the elements $\mathbf{S}_{6k+i+j+\pier{n_i}}$ span the $6$-th dimensional linear eigenspace of symmetric tensors associated to the eigenvalue $\gamma_k$. 
We also introduce the projection operator
\[
PS_m:\pier{H^1(\mathcal{D}_a;\mathbb{R}^{3\times 3})} \to \text{span}\{\mathbf{S}_0,\mathbf{S}_1,\dots,\mathbf{S}_{6m+5}\}.
\]

We moreover introduce the functions $\{\mathbf{A}_{3k+i+j-1}\}_{k\in \mathbb{N}; i,j=0,\dots,2; j> i}$ defined by
\[
 \mathbf{A}_{3k+i+j-1}:=\xi_k\left(\mathbf{e}_i\otimes \mathbf{e}_j-\mathbf{e}_j\otimes \mathbf{e}_i\right).
\]
We observe that, given $k\in \mathbb{N}$, the elements $\mathbf{A}_{3k+i+j-1}$ span the $3$-th dimensional linear eigenspace of antisymmetric tensors associated to the eigenvalue $\gamma_k$. 
We then introduce the projection operator
\[
PA_m:\pier{H^1(\mathcal{D}_a;\mathbb{R}^{3\times 3})}\to \text{span}\{\mathbf{A}_0,\mathbf{A}_1,\dots,\mathbf{A}_{3m+2}\}.
\]
We make the Galerkin ansatz
\begin{equation}
\label{galerkinansatz}
\pier{\mathbf{W}_m(\mathbf{a}, t)=\mathbf{I}+\sum_{i=0}^{6m+5}x_i^m(t)\mathbf{S}_i(\mathbf{a}), \quad 
\boldsymbol{\Theta}_m(\mathbf{a}, t)=\sum_{i=0}^{3m+2}y_i^m(t)\mathbf{A}_i(\mathbf{a}),
\quad (\mathbf{a}, t) \in \mathcal{D}_a \times (0,T),}
\end{equation}
with 
\begin{align*}
&\mathbf{S}_i\in \pier{H^2}(\mathcal{D}_a;Sym(\mathbb{R}^{3\times 3}))\cap H_0^1(\mathcal{D}_a;Sym(\mathbb{R}^{3\times 3})), \\
&\mathbf{A}_i\in \pier{H^2}(\mathcal{D}_a;Skew(\mathbb{R}^{3\times 3}))\cap H_0^1(\mathcal{D}_a;Skew(\mathbb{R}^{3\times 3})),
\end{align*}
to approximate the solutions $\mathbf{W}$ and $\boldsymbol{\Theta}$ of \pier{the system \eqref{eqn:48}--\eqref{eqn:50}}. \abramo{We observe that through the Galerkin ansatz \eqref{galerkinansatz} we are enforcing by construction that $\mathbf{W}_m\in \mathcal{S}$ and $\boldsymbol{\Theta}_m\in \mathcal{A}$, \abramonew{hence the system \eqref{eqn:48}--\eqref{eqn:50} is equivalent to the system \eqref{eqn:48symskew}--\eqref{eqn:50symskew}}.}
\pier{We consider a Faedo--Galerkin approximation of a regularized version of \eqref{eqn:48symskew}, \abramo{with solutions expressed in the form \eqref{galerkinansatz} and}
where the convex functions (cf. Remark~\ref{rem:pier1}) $\psi_{C_{\alpha}}$ and $\psi_D$ and their subdifferentials $\partial \psi_{C_{\alpha}}$ and $\partial \psi_D$ are replaced by the Moreau--Yosida \pier{approximations} $ \psi_{C_{\alpha}}^{\lambda}$ and $\psi_D^{\lambda}$, $\partial \psi_{C_{\alpha}}^{\lambda}$ and $\partial \psi_D^{\lambda}$, depending on a regularization parameter $\lambda >0$. 
We refer to, e.g., \cite[pp.~28 and~39]{brezismm}) for definitions and properties of these approximations, recalling simply that if  $f:L^2(\mathcal{D}_a;\mathbb{R}^{3\times 3})\to [0,+\infty] $ is a proper convex lower semicontinuous function and $\partial f$ denotes its subdifferential, then  
\[
 \partial f^{\lambda} :=\frac{I-\left(I+\lambda \partial f \right)^{-1}}{\lambda}, \quad \lambda \in (0,1),
\]
where $I$ here denotes the identity operator}. In particular, $\partial f^{\lambda}$ is a monotone and $\frac{1}{\lambda}$-Lipschitz continuous function. Moreover, \pier{due the special form of $\psi_D$ defined in \eqref{eqn:22bis},
we have that the following bounds are valid uniformly in $\lambda$:
\begin{align}
\label{eqn:dpsi1}
& \displaystyle \frac{1}{2}\pier{\|}\mathbf{A}\pier{\|}^2\leq C+\psi_D^{\lambda}(\mathbf{A}), \quad \hbox{for all } \, \mathbf{A}\in L^2(\mathcal{D}_a;\mathbb{R}^{3\times 3}),\\
\label{eqn:dpsi2}
& \displaystyle \pier{\|}\partial \psi_D^{\lambda}(\mathbf{A})\pier{\|}^2\leq C\bigl( \psi_D^{\lambda}(\mathbf{A})+1\bigr), \quad \hbox{for all } \, \mathbf{A}\in L^2(\mathcal{D}_a;\mathbb{R}^{3\times 3})\pco ,
\end{align}
}%
\pco{where the constant $C$ is also independent of $k$ provided that $0< k \leq \overline{k}$, for some  
$\overline{k}>0$.}
Given \eqref{galerkinansatz}, we define the approximations
\begin{align}
 \label{eqn:sigmammapp}
&\pier{{}\boldsymbol{\chi}_{\alpha,m}{} = \partial \psi_{C_\alpha}^{\lambda}\left(\mathbf{W}_m\right),}\quad 
\boldsymbol{\Sigma}_m=\partial \psi_D^{\lambda}\left(\curl \left(\mathbf{Z}_m\right)\right),
\\
 \label{eqn:phimzmapp}
&\Delta \pier{\mathbf{\Phi}}_m=\diver\left(\mathrm{e}^{\boldsymbol{\Theta}_m}\mathbf{W}_m\right), \; - P_L\Delta {\mathbf{Z}_m}=\curl \left(\mathrm{e}^{\boldsymbol{\Theta_m}}\mathbf{W}_m\right),\notag \\
&\quad \hbox{with } \,   \pier{\mathbf{\Phi}_m(\mathbf{a}, t) =\mathbf{a}, \; \mathbf{Z}_m (\mathbf{a}, t)=\mathbf{0} \quad \text{for} \; (\mathbf{a}, t) \in \Gamma_a\times (0,T)}.
\end{align}
Given the elliptic \pier{problems in \eqref{eqn:phimzmapp} with approximated right hand sides, we then have
\begin{align}
 \label{phimzm}
& \pier{\mathbf{\Phi}}_m(\mathbf{a}, t)= \mathbf{a} -\mathcal{G}_L\ast \diver\left(\mathrm{e}^{\boldsymbol{\Theta}_m}\mathbf{W}_m\right)(\mathbf{a}, t), \notag \\ 
&\mathbf{Z}_m(\mathbf{a}, t)=\mathcal{G}_{L,\diver}\ast \curl\left(\mathrm{e}^{\boldsymbol{\Theta}_m}\mathbf{W}_m\right)(\mathbf{a}, t)
\quad \text{for} \;  (\mathbf{a}, t) \in \mathcal{D}_a\times (0,T). 
\end{align}
}
\abramonew{
We project the equation \eqref{eqn:48symskew}$_1$ for $\mathbf{W}_m$ onto $\text{span}\left\{\mathbf{S}_0,\mathbf{S}_1,\dots,\mathbf{S}_{6m+5}\right\}$, the equation \eqref{eqn:48symskew}$_2$ for $\boldsymbol{\Theta}_m$ onto
$\text{span}\left\{\mathbf{A}_0,\mathbf{A}_1,\dots,\mathbf{A}_{3m+2}\right\},$
with $\mathbf{Z}_m$ defined as in \eqref{phimzm} and \pier{$\boldsymbol{\chi}_{\alpha,m},\, \boldsymbol{\Sigma}_m$ defined in \eqref{eqn:sigmammapp}, obtaining the following Galerkin approximation of \eqref{eqn:48}}:
\begin{equation}
\label{eqn:48m}
\begin{cases}
\int_{\mathcal{D}_a}\text{Sym}\left(\mathrm{e}^{-\boldsymbol{\Theta}_m}\curl\left(
\mathcal{G}_{L,\diver}\ast \curl \left[\partial \psi_D^{\lambda}\left(\curl \left(\mathcal{G}_{L,\diver}\ast \curl\left(\mathrm{e}^{\boldsymbol{\Theta}_m}\mathbf{W}_m\right)\right)\right)\right]\right)\right)\colon \mathbf{S}_i\\
+ \int_{\mathcal{D}_a}\left(\mathbf{W}_m-\mathbf{I}+\partial \pier{{}\psi_{C_\alpha}^{\lambda}(\mathbf{W}_m)}+\overbigdot{\mathbf{W}}_m\right)\colon \mathbf{S}_i+\int_{\mathcal{D}_a}\pier{\grad{}} \mathbf{W}_m::\pier{\grad{}} \mathbf{S}_i
\\=\int_{\mathcal{D}_a}\mathbf{W}_{\text{ext}}(\mathbf{W}_m,t)\colon \mathbf{S}_i,\\ \\
\int_{\mathcal{D}_a}\text{Skew}\left(\left(\curl\left(
\mathcal{G}_{L,\diver}\ast \curl \left[\partial \psi_D^{\lambda}\left(\curl \left(\mathcal{G}_{L,\diver}\ast \curl\left(\mathrm{e}^{\boldsymbol{\Theta}_m}\mathbf{W}_m\right)\right)\right)\right]\right)\mathbf{W}_m\mathrm{e}^{-\boldsymbol{\Theta}_m}\right)\right)\colon \mathbf{A}_j \\
+\frac{1}{2}\int_{\mathcal{D}_a}\pier{\grad{}} \boldsymbol{\Theta}_m::\pier{\grad{}} \mathbf{A}_j+\frac{1}{2}\int_{\mathcal{D}_a}\pier{\grad{}} \overbigdot{\boldsymbol{\Theta}}_m::\pier{\grad{}} \mathbf{A}_j=\int_{\mathcal{D}_a}\boldsymbol{\Omega}_{\text{ext}}\abramo{(\boldsymbol{\Theta}_m,t)}\colon \mathbf{A}_j,\\ \\
\Delta \pier{\mathbf{\Phi}}_m=\diver\left(\mathrm{e}^{\boldsymbol{\Theta}_m}\mathbf{W}_m\right),\\ \\
- P_L\Delta {\mathbf{Z}_m}=\curl \left(\mathrm{e}^{\boldsymbol{\Theta_m}}\mathbf{W}_m\right),
\end{cases}
\end{equation}
}
{in $[0,t]$}, with $0<t\leq T$, for $i=0, \dots, 6m+5$, $j=0, \dots, 3m+2$, \pier{with boundary conditions as in \eqref{eqn:49symskew}} and with initial conditions \pier{(cf. the assumption \textbf{A2})}
\begin{equation}
    \label{eqn:48mic}
   \pier{ \mathbf{W}_m(\mathbf{a},0)=\mathbf{I}+PS_m(\mathbf{W}_0-\mathbf{I})(\mathbf{a}), \quad \boldsymbol{\Theta}_m(\mathbf{a},0)=\boldsymbol{0}, \quad \mathbf{a}\in \mathcal{D}_a.}
\end{equation}
The equations \eqref{eqn:48m}$_1$ and  \eqref{eqn:48m}$_2$ are decoupled from the other equations in \pier{the system} \eqref{eqn:48m} and define a collection of initial value problems for a system of coupled ODEs of the form
\begin{equation}
\label{eqn:48ode}
\begin{cases}
\displaystyle \frac{d}{dt}x_i^m=-(1+\gamma_i)x_i^m+\int_{\mathcal{D}_a}\biggl(-\partial \pier{{}\psi_{C_\alpha}^{\lambda}}\biggl(\pier{{}\mathbf{I} +{}}\sum_lx_l^m\mathbf{S}_l\biggr)+\mathbf{W}_{\text{ext}}\biggl(\pier{{}\mathbf{I} +{}} \sum_lx_l^m\mathbf{S}_l,t\biggr)\biggr)\colon \mathbf{S}_i\\
\displaystyle = \int_{\mathcal{D}_a}PS_m\biggl[\mathrm{e}^{-\sum_ly_l^m\mathbf{A}_l}\curl\biggl(\mathcal{G}_{L,\diver}\ast \curl \biggl[\partial \psi_D^{\lambda}\biggl(\curl \biggl(
\mathcal{G}_{L,\diver}\\
\hskip5cm \displaystyle \ast \curl\biggl(\mathrm{e}^{\sum_ry_r^m\mathbf{A}_r}\biggl(\pier{{}\mathbf{I} +{}}\sum_kx_k^m\mathbf{S}_k\biggr)\biggr)\biggr)\biggr)\biggr]\biggr)\biggr]\colon \mathbf{S}_i,\\ \\
\displaystyle \frac{d}{dt}y_j^m=-y_j^m+\frac{2}{\gamma_j}\int_{\mathcal{D}_a}\boldsymbol{\Omega}_{\text{ext}}\abramo{(\boldsymbol{\Theta}_m,t)}\colon \mathbf{A}_j\\
\displaystyle -\frac{2}{\gamma_j}\int_{\mathcal{D}_a}PA_m\biggl[\curl \biggl(\mathcal{G}_{L,\diver}\ast \curl \biggl[\partial \psi_D^{\lambda}\biggl(\curl \biggl(
\mathcal{G}_{L,\diver}\\
\hskip1cm \displaystyle \ast \curl\biggl(\mathrm{e}^{\sum_ry_r^m\mathbf{A}_r}\biggl(\pier{{}\mathbf{I} +{}}\sum_kx_k^m\mathbf{S}_k\biggr)
\biggr)\!\biggr)\!\biggr)\!\biggr]\!\biggr) \biggl(\pier{{}\mathbf{I} +{}}\sum_kx_k^m\mathbf{S}_k\biggr)\mathrm{e}^{-\sum_py_p^m\mathbf{A}_p}\biggr]\colon \mathbf{A}_j, \\ \\
\displaystyle x_i^m(0)=\int_{\mathcal{D}_a}\biggl(\mathbf{W}_0-\mathbf{I}\biggr)\colon \mathbf{S}_i, \, \  y_j^m(0)=0, \quad i=0, \dots, 6m+5, \quad j=0, \dots, 3m+2.
\end{cases}
\end{equation}
Due to Assumptions \textbf{A3}, to the Lipschitz \pier{continuity of $\partial \pier{{}\psi_{C_\alpha}^{\lambda}}$ and $\partial \psi_D^{\lambda}$} and to the regularity in space of the functions $\mathbf{S}_i,  \mathbf{A}_j$, \pier{the system}~\eqref{eqn:48ode} is a coupled system of \pier{first-order} ODEs in the variables $x_i^m, \, y_j^m$, with a right hand side which is measurable in time and continuous in the independent variables. Then,  we can apply the Carath\'eodory's existence theorem to infer that there exist a sufficiently small $t_1$ with $0<t_1\leq T$ and a local solution $(x_i^m,y_j^m)$ of \eqref{eqn:48ode}, for $i=0, \dots, 6m+5$, $j=0, \dots, 3m+2$, which is absolutely continuous. \pier{Once we have} a solution to \eqref{eqn:48ode}, 
dealing with \pier{the} elliptic problems with regular right-hand sides in \pier{\eqref{eqn:48m} leads to the elements $\pier{\mathbf{\Phi}}_m$ and $\mathbf{Z}_m$ solving} \eqref{eqn:48m}$_3$ and \eqref{eqn:48m}$_4$, respectively. 

\pier{Next,
thanks to some uniform} estimates, we will extend these solutions by continuity to the interval $[0,T]$ and we will study the limit as $m\to \infty$ and $\lambda\to 0$. In particolar, we will study \pier{in} a first step the limit as $m\to \infty$, and \pier{then} the limit as $\lambda\to 0$ in the latter limit system.

We now deduce a priori estimates, uniform in the discretization parameter $m$ and in the regularization parameter $\lambda$, for the solutions of system \eqref{eqn:48m}, which can be rewritten, combining the equations over $i=0, \dots, 6m+5$ and $j=0, \dots, 3m+2$, as
\begin{equation}
\label{eqn:48m2}
\begin{cases}
\int_{\mathcal{D}_a}\rm{Sym}\left(\mathrm{e}^{-\boldsymbol{\Theta}_m}\curl\left(\mathcal{G}_{L,\diver}\ast \left(
\curl \boldsymbol{\Sigma}_m\right)\right)\right)\colon \nwhat{\mathbf{W}}_m\\
+\int_{\mathcal{D}_a}\left(\mathbf{W}_m-\mathbf{I}+\pier{{}\boldsymbol{\chi}_{\alpha,m}{}}+\overbigdot{\mathbf{W}}_m\right)\colon \nwhat{\mathbf{W}}_m\\
+\int_{\mathcal{D}_a}\pier{\grad{}} \mathbf{W}_m::\pier{\grad{}}\nwhat{\mathbf{W}}_m=\int_{\mathcal{D}_a}\mathbf{W}_{\text{ext}}(\mathbf{W}_m,t)\colon \nwhat{\mathbf{W}}_m,\\ \\
\int_{\mathcal{D}_a}\rm{Skew}\left(\curl\left(\mathcal{G}_{L,\diver}\ast
\left(\curl \boldsymbol{\Sigma}_m\right)\right)\mathbf{W}_m\mathrm{e}^{-\boldsymbol{\Theta}_m}\right)\colon \nwhat{\boldsymbol{\Omega}}_m +
\frac{1}{2}\int_{\mathcal{D}_a}\pier{\grad{}} \boldsymbol{\Theta}_m::\pier{\grad{}} \nwhat{\boldsymbol{\Omega}}_m\\
+\frac{1}{2}\int_{\mathcal{D}_a}\pier{\grad{}} \overbigdot{\boldsymbol{\Theta}}_m::\pier{\grad{}} \nwhat{\boldsymbol{\Omega}}_m=\int_{\mathcal{D}_a}\boldsymbol{\Omega}_{\text{ext}}\abramo{(\boldsymbol{\Theta}_m,t)}\colon \nwhat{\boldsymbol{\Omega}}_m,\\ \\
\pier{{}\boldsymbol{\chi}_{\alpha,m}{} = \partial \psi_{C_\alpha}^{\lambda}\left(\mathbf{W}_m\right),}\quad
\boldsymbol{\Sigma}_m=\partial \psi_D^{\lambda}(\curl \mathbf{Z}_m),\\ \\
\Delta \pier{\mathbf{\Phi}}_m=\diver\left(\mathrm{e}^{\boldsymbol{\Theta}_m}\mathbf{W}_m\right),\\ \\
-P_L\Delta {\mathbf{Z}_m}=\curl \left(\mathrm{e}^{\boldsymbol{\Theta_m}}\mathbf{W}_m\right),
\end{cases}
\end{equation}
for a.e. $t \in [0,t_1]$ \pier{and all} $\nwhat{\mathbf{W}}_m \in \text{span}\left\{\mathbf{S}_0,\mathbf{S}_1,\dots,\mathbf{S}_{6m+5}\right\}$, $\nwhat{\boldsymbol{\Omega}}_m \in \text{span}\left\{\mathbf{A}_0,\mathbf{A}_1,\dots,\mathbf{A}_{3m+2}\right\}$, and with initial conditions defined in \eqref{eqn:48mic}. 

The first a-priori estimate is obtained by taking $\nwhat{\mathbf{W}}_m=\overbigdot{\mathbf{W}}_m$ in \eqref{eqn:48m2}$_1$ and $\nwhat{\boldsymbol{\Omega}}_m=\overbigdot{\boldsymbol{\Theta}}_m$ in \eqref{eqn:48m2}$_2$. Moreover, we take the time derivative of \eqref{eqn:48m2}$_5$, multiply it by $\mathcal{G}_{L,\diver}\ast\left( \curl \boldsymbol{\Sigma}_m\right)$ and integrate over $\mathcal{D}_a$. 
We observe from \eqref{eqn:48m2}$_5$ and from the regularity in space of the functions $\mathbf{S}_i, \mathbf{A}_j$ that $\mathbf{Z}_m\in H^3(\mathcal{D}_a;\mathbb{R}^{3\times 3})$, for any $t\in [0,t_1]$. Hence, from \eqref{eqn:48m2}$_3$ and the Lipschitz continuity of $\partial \psi_D^{\lambda}$ we obtain that $\boldsymbol{\Sigma}_m\in H^1(\mathcal{D}_a;\mathbb{R}^{3\times 3})$, and as a consequence the $L^2(\mathcal{D}_a;\mathbb{R}^{3\times 3})$ scalar product of equation~\pier{\eqref{eqn:48m2}$_5$} with the element $\mathcal{G}_{L,\diver}\ast\left( \curl \boldsymbol{\Sigma}_m\right)\in H^2(\mathcal{D}_a;\mathbb{R}^{3\times 3})$ is well defined for any $t\in [0,t_1]$.
Finally we sum all the previous contributions and integrate in time between $0$ and $t\in [0,t_1]$. Observing that $\overbigdot{\mathbf{W}}_m=\overbigdot{\mathbf{W}}_m^T$ and $\overbigdot{\boldsymbol{\Theta}}_m=-\overbigdot{\boldsymbol{\Theta}}_m^T$, and since $\mathbf{A}\colon \mathbf{B}=\mathbf{A}^T\colon \mathbf{B}^T$ for any $\mathbf{A},\mathbf{B}\in M(\mathbb{R}^{3\times 3})$, we have that
\begin{align*}
&\frac{1}{2}\int_{\mathcal{D}_a}\left(\mathrm{e}^{-\boldsymbol{\Theta}_m}\curl\left(
\mathcal{G}_{L,\diver}\ast \left(\curl\boldsymbol{\Sigma}_m\right)\right)+\left[\curl\left(
\mathcal{G}_{L,\diver}\ast \left(\curl\boldsymbol{\Sigma}_m\right)\right)\right]^T\mathrm{e}^{\boldsymbol{\Theta}_m}\right)\colon \overbigdot{\mathbf{W}}_m\\
&= \int_{\mathcal{D}_a}\mathrm{e}^{-\boldsymbol{\Theta}_m}\curl\left(
\mathcal{G}_{L,\diver}\ast \left(\curl\boldsymbol{\Sigma}_m\right)\right)\colon \overbigdot{\mathbf{W}}_m=\int_{\mathcal{D}_a}
\curl\left(\mathcal{G}_{L,\diver}\ast \left(\curl\boldsymbol{\Sigma}_m\right)\right)\colon \mathrm{e}^{\boldsymbol{\Theta}_m}\overbigdot{\mathbf{W}}_m,
\end{align*}
and
\begin{align*}
&\frac{1}{2}\int_{\mathcal{D}_a}\left(\curl\left( \mathcal{G}_{L,\diver}\ast \left(\curl\boldsymbol{\Sigma}_m\right)\right)\mathbf{W}_m\mathrm{e}^{-\boldsymbol{\Theta}_m}-\mathrm{e}^{\boldsymbol{\Theta}_m}\mathbf{W}_m\left[\curl\left(\mathcal{G}_{L,\diver}\ast \left(\curl\boldsymbol{\Sigma}_m\right)\right)\right]^T\right)\colon \overbigdot{\boldsymbol{\Theta}}_m\\
& =\int_{\mathcal{D}_a}\curl\left(\mathcal{G}_{L,\diver}\ast \left(\curl\boldsymbol{\Sigma}_m\right)\right)\mathbf{W}_m\mathrm{e}^{-\boldsymbol{\Theta}_m}\colon \overbigdot{\boldsymbol{\Theta}}_m\\
&=\int_{\mathcal{D}_a}\curl\left(\mathcal{G}_{L,\diver}\ast \left(\curl\boldsymbol{\Sigma}_m\right)\right)\colon \overbigdot{\boldsymbol{\Theta}}_m\mathrm{e}^{\boldsymbol{\Theta}_m}\mathbf{W}_m.
\end{align*}
Also, the contribution from \pier{\eqref{eqn:48m2}$_5$}, after integration by parts, gives that
\begin{align*}
&\int_{\mathcal{D}_a}-\Delta \overbigdot{\mathbf{Z}}_m\colon \mathcal{G}_{L,\diver}\ast \left(\curl\boldsymbol{\Sigma}_m\right)=\int_{\mathcal{D}_a} \overbigdot{\mathbf{Z}}_m\colon \curl\boldsymbol{\Sigma}_m=\int_{\mathcal{D}_a} \curl \overbigdot{\mathbf{Z}}_m\colon \partial \psi_D^{\lambda}(\curl \mathbf{Z}_m)\\
& =\int_{\mathcal{D}_a}\mathrm{e}^{\boldsymbol{\Theta}_m}\overbigdot{\mathbf{W}}_m\colon \curl\left(
\mathcal{G}_{L,\diver}\ast \left(\curl\boldsymbol{\Sigma}_m\right)\right)+\int_{\mathcal{D}_a}\overbigdot{\boldsymbol{\Theta}}_m\mathrm{e}^{\boldsymbol{\Theta}_m}\mathbf{W}_m\colon \curl\left( \mathcal{G}_{L,\diver}\ast \left(\curl\boldsymbol{\Sigma}_m\right)\right).
\end{align*}
Hence\pier{, for any $t\in [0,t_1]$, we deduce} that
\begin{align}
\label{eqn:53}
& \displaystyle \notag \frac{1}{2}\pier{\|}\mathbf{W}_m-\mathbf{I}\pier{\|}^2
+\pier{\psi_{C_{\alpha}}^{\lambda}}(\mathbf{W}_m)
+\frac{1}{2}\pier{\|}\grad \mathbf{W}_m\pier{\|}^2+\frac{1}{4}\pier{\|} \boldsymbol{\Theta}_m\pier{\|}^2+\frac{1}{4}\pier{\|}\grad \boldsymbol{\Theta}_m\pier{\|}^2\\
& \displaystyle \notag  +\psi_D^{\lambda}(\curl \mathbf{Z}_m)+\int_0^{t_1}\pier{\|}\overbigdot{\mathbf{W}}_m
\pier{\|}^2+\frac{1}{2}\int_0^{t_1}\pier{\|}\grad \overbigdot{\boldsymbol{\Theta}}_m\pier{\|}^2\\
&= \displaystyle \notag \frac{1}{2}\pier{\|}\mathbf{W}_m(0)-\mathbf{I}\pier{\|}^2
+\pier{\psi_{C_{\alpha}}^{\lambda}}(\mathbf{W}_m(0))
+\frac{1}{2}\pier{\|}\grad \mathbf{W}_m(0)\pier{\|}^2
+\frac{1}{4}\pier{\|}\boldsymbol{\Theta}_m(0)\pier{\|}^2\\
& \displaystyle \notag  +\frac{1}{4}\pier{\|}\grad \boldsymbol{\Theta}_m(0)\pier{\|}^2+\psi_D^{\lambda}(\curl\mathbf{Z}_m(0))+\frac{1}{2}\int_{\mathcal{D}_{at_1}}\overbigdot{\boldsymbol{\Theta}}_m\colon {\boldsymbol{\Theta}}_m\\
& \displaystyle \notag+ \int_{\mathcal{D}_{at_1}}\mathbf{W}_{\text{ext}}(\mathbf{W}_m,t)\colon \overbigdot{\mathbf{W}}_m+ \int_{\mathcal{D}_{at_1}}\boldsymbol{\Omega}_{\text{ext}}\abramo{(\boldsymbol{\Theta}_m,t)}\colon \overbigdot{\boldsymbol{\Theta}}_m\\
& \displaystyle\notag \leq C\Vert \mathbf{W}_m(0)\Vert +C+ \frac{1}{2}\int_0^{t_1}\pier{\|}\overbigdot{\mathbf{W}}_m\pier{\|}^2+\frac{1}{4}\int_0^{t_1}\pier{\|}\grad \overbigdot{\boldsymbol{\Theta}}_m\pier{\|}^2\\
&\displaystyle +\abramo{C\int_0^{t_1}G^2\left(\pier{\|} \boldsymbol{\Theta}_m\pier{\|}^2+\pier{\|}\grad \boldsymbol{\Theta}_m\pier{\|}^2\right)}+C\int_0^{t_1}L^2\left(\pier{\|}\mathbf{W}_m-\mathbf{I}\pier{\|}^2+\pier{\|}\grad \mathbf{W}_m\pier{\|}^2\right),
\end{align}
where we added $\frac{1}{4}\frac{d}{dt}\pier{\|} \boldsymbol{\Theta}_m\pier{\|}^2$ to the left and $\frac{1}{2}\int_{\mathcal{D}_{a}}\overbigdot{\boldsymbol{\Theta}}_m\colon {\boldsymbol{\Theta}}_m$ to the right, and used \reviewa{the Cauchy--Schwarz, the Young and the Poincar\'{e} inequalities}
, Assumptions \textbf{A2} and \textbf{A3}. Thanks to the Gronwall lemma, we thus have that
\begin{align}
\label{eqn:54}
& \displaystyle \notag \frac{1}{2}\pier{\|}\mathbf{W}_m-\mathbf{I}\pier{\|}^2
+\pier{\psi_{C_{\alpha}}^{\lambda}}(\mathbf{W}_m)
+\frac{1}{2}\pier{\|}\grad \mathbf{W}_m\pier{\|}^2+\frac{1}{4}\pier{\|} \boldsymbol{\Theta}_m\pier{\|}^2+\frac{1}{4}\pier{\|}\grad \boldsymbol{\Theta}_m\pier{\|}^2\\
& \displaystyle +\psi_D^{\lambda}(\curl \mathbf{Z}_m)+\frac{1}{2}\int_0^{t_1}\pier{\|}\overbigdot{\mathbf{W}}_m\pier{\|}^2+\frac{1}{2}\int_0^{t_1}\pier{\|}\grad \overbigdot{\boldsymbol{\Theta}}_m\pier{\|}^2\leq C,
\end{align}
where the constant in the right hand side of \eqref{eqn:54} depends only on the initial data, on the domain $\mathcal{D}_a$ and not on the discretization parameter $m$ and on the regularization parameter $\lambda$. Thanks to the a priori estimate \eqref{eqn:54}, we may extend {by continuity} the local solution of system \eqref{eqn:48m2} to the interval $[0,T]$.
Using \eqref{eqn:dpsi1} and \eqref{eqn:54} we have that
\begin{equation}
\label{eqn:542}
\sup_{t\in (0,T)}\pier{\|} (\curl \mathbf{Z}_m ) \pier{(t)} \pier{\|}^2\leq C.
\end{equation}
\pier{Moreover, in view of \eqref{eqn:dpsi2}, from \eqref{eqn:48m2}$_3$ and \eqref{eqn:54} it follows} that
\begin{equation}
 \label{eqn:543}
 \sup_{t\in (0,T)}\pier{\|}\boldsymbol{\Sigma}_m \pier{(t)} \pier{\|}^2\leq C.
\end{equation}
We now multiply \pier{the equality $\boldsymbol{\Sigma}_m=\partial \psi_D^{\lambda}(\curl \mathbf{Z}_m)$ in} \eqref{eqn:48m2}$_3$ by $\curl \left(\mathcal{G}_{L,\diver}\ast \left(\curl \boldsymbol{\Sigma}_m\right)\right)\in H^1(\mathcal{D}_a;\mathbb{R}^{3\times 3})$ and integrate over $\mathcal{D}_a$. Employing multiple integration by parts, the Cauchy--Schwarz and Young inequalities and \eqref{eqn:dpsi2}, we obtain that
\begin{align*}
\displaystyle & \int_{\mathcal{D}_a}\boldsymbol{\Sigma}_m:\curl \left(\mathcal{G}_{L,\diver}\ast \left(\curl \boldsymbol{\Sigma}_m\right)\right)=\int_{\mathcal{D}_a}\curl \boldsymbol{\Sigma}_m: \mathcal{G}_{L,\diver}\ast \left(\curl \boldsymbol{\Sigma}_m\right)\\
\displaystyle &=  \pier{\|}\curl \left(\mathcal{G}_{L,\diver}\ast \left(\curl \boldsymbol{\Sigma}_m\right)\right):\curl \left(\mathcal{G}_{L,\diver}\ast \left(\curl \boldsymbol{\Sigma}_m\right)\right)\pier{\|}^2\\
\displaystyle &= \int_{\mathcal{D}_a}\partial \psi_D^{\lambda}(\curl \mathbf{Z}_m):\curl \left(\mathcal{G}_{L,\diver}\ast \left(\curl \boldsymbol{\Sigma}_m\right)\right)\\
\displaystyle &\leq  C\psi_D^{\lambda}(\curl \mathbf{Z}_m)+C+\frac{1}{2}\pier{\|}\curl \left(\mathcal{G}_{L,\diver}\ast \left(\curl \boldsymbol{\Sigma}_m\right)\right):\curl \left(\mathcal{G}_{L,\diver}\ast \left(\curl \boldsymbol{\Sigma}_m\right)\right)\pier{\|}^2.
\end{align*}
Hence, given the estimate \eqref{eqn:54}, we have that 
\begin{equation}
 \label{eqn:544}
 \sup_{t\in (0,T)}\pier{\|}\curl \boldsymbol{\Sigma}_m\pier{(t)}\pier{\|}_{\left(H_{0,\diver}^1\left(\mathcal{D}_a,\mathbf{R}^{3\times 3}\right)\right)^{\prime}}^2\leq C,
\end{equation}
and, from a Lax--Milgram estimate associated to the operator $-P_L\Delta$, 
\begin{equation}
 \label{eqn:545}
\sup_{t\in (0,T)} \pier{\|}\pier{(\mathcal{G}_{L,\diver}\ast (\curl \boldsymbol{\Sigma}_m))}\pier{(t)}\pier{\|}_{H_{0,\diver}^1\left(\mathcal{D}_a,\mathbf{R}^{3\times 3}\right)}^2\leq C.
\end{equation}
The second \pier{a priori} estimate is obtained by taking $\nwhat{\mathbf{W}}_m=-\Delta \mathbf{W}_m$ in \eqref{eqn:48m2}$_1$ and integrating in time between $0$ and $t\in[0,T]$. \pier{Using Assumptions~\textbf{A2}, \textbf{A3} and estimate~\eqref{eqn:54}, we infer} that
\begin{align}
\label{eqn:55}
& \displaystyle \notag \frac{1}{2}\pier{\|}\grad \mathbf{W}_m\pier{\|}^2+\underbrace{\int_0^t (\pier{\grad (\partial I^\lambda_{C_{\alpha}}(\mathbf{W}_m)) ,\grad \mathbf{W}_m)}}_{\geq 0}+\int_0^t\pier{\|}\Delta{\mathbf{W}_m}\pier{\|}^2\\
& \displaystyle\notag \leq C+\frac{1}{2}\int_0^t\pier{\|}\Delta{\mathbf{W}_m}\pier{\|}^2+ C\int_0^tL^2\left(\pier{\|} \mathbf{W}_m\pier{\|}^2+\pier{\|}\grad \mathbf{W}_m\pier{\|}^2\right)\\
&\displaystyle +C\left\|\mathrm{e}^{-\boldsymbol{\Theta}_m}\right\|_{L^{\infty}(\mathcal{D}_{at},\mathbb{R}^{3\times 3})}^2\int_0^t\pier{\|}\curl\left( \mathcal{G}_{L,\diver}\ast \left(\curl \boldsymbol{\Sigma}_m\right)\right)\pier{\|}^2.
\end{align}
\pier{Hence, by \eqref{eqn:54} and \eqref{eqn:545} the right-hand side is under control and then}
\begin{align}
\label{eqn:55bis}
& \displaystyle \frac{1}{2}\pier{\|}\grad \mathbf{W}_m\pier{\|}^2+\frac{1}{2}\int_0^t\pier{\|}\Delta{\mathbf{W}_m}\pier{\|}^2\leq C.
\end{align}
We derive a further \pier{a priori} estimate by taking $\nwhat{\boldsymbol{\Omega}}_m=-\Delta \boldsymbol{\overbigdot{\Theta}}_m$ in \eqref{eqn:48m2}$_2$ and integrating in time between $0$ and $t\in[0,T]$. Using Assumption \textbf{A3} we obtain that
\begin{align}
\label{eqn:56}
& \notag \displaystyle \frac{1}{4}\pier{\|}\Delta \boldsymbol{\Theta}_m\pier{\|}^2+\frac{1}{2}\int_0^t \pier{\|}\Delta \overbigdot{\boldsymbol{\Theta}}_m\pier{\|}^2\\
& \displaystyle \notag  \leq \abramo{C\int_0^{t_1}G^2\left(\pier{\|} \boldsymbol{\Theta}_m\pier{\|}^2+\pier{\|}\grad \boldsymbol{\Theta}_m\pier{\|}^2\right)}+\frac{1}{4}\int_0^t\pier{\|}\Delta \overbigdot{\boldsymbol{\Theta}}_m\pier{\|}^2+C\left\|\mathrm{e}^{-\boldsymbol{\Theta}_m}\right\|_{L^{\infty}(\mathcal{D}_{at},\mathbb{R}^{3\times 3})}^2 \\
&\times \pier{\|}\curl\left( \mathcal{G}_{L,\diver}\ast \left(\curl \boldsymbol{\Sigma}_m\right)\right)\pier{\|}_{L^{\infty}(0,t;L^2(\mathcal{D}_{a},\mathbb{R}^{3\times 3}))}^2\pier{\|}\mathbf{W}_m\pier{\|}_{L^{2}(0,t;L^{\infty}(\mathcal{D}_{a},Sym(\mathbb{R}^{3\times 3})))}^2,
\end{align}
\pier{whence}, using \eqref{eqn:54}, \eqref{eqn:545}, \eqref{eqn:55bis} and the Sobolev embedding $H^2 \hookrightarrow L^{\infty}$ (obtained from \eqref{eqn:2} with \pier{$j=0$, $p=\infty$, $m=r=q=2$}), we \pier{find out} that
\begin{align}
\label{eqn:56bis}
& \displaystyle \frac{1}{4}\pier{\|}\Delta \boldsymbol{\Theta}_m\pier{\|}^2+\frac{1}{4}\int_0^t \pier{\|}\Delta \overbigdot{\boldsymbol{\Theta}}_m\pier{\|}^2\leq C.
\end{align}
Thanks to \eqref{eqn:54}, \eqref{eqn:55bis}, \eqref{eqn:56bis} and \eqref{eqn:1}, from \eqref{eqn:48m2}$_4$, \eqref{eqn:48m2}$_5$ and \pier{the estimates for the time derivatives of} $\pier{\mathbf{\Phi}}_m$ and $\mathbf{Z}_m$ \pier{we arrive at}
\begin{align}
\label{eqn:57}
& \displaystyle \pier{\|}\pier{\mathbf{\Phi}}_m\pier{\|}_{L^{\infty}\left(0,T;H^2(\mathcal{D}_a;\mathbb{R}^3)\right)\cap L^{2}\left(0,T;H^{3}(\mathcal{D}_a,\mathbb{R}^3)\right)\cap H^{1}\left(0,T;H^{1}(\mathcal{D}_a,\mathbb{R}^3)\right)} \leq C,\\
\label{eqn:58}
& \displaystyle \pier{\|}\mathbf{Z}_m\pier{\|}_{L^{\infty}\left(0,T;H^2(\mathcal{D}_a;\mathbb{R}^{3\times 3})\right)\cap L^{2}\left(0,T;H^{3}(\mathcal{D}_a,\mathbb{R}^{3\times 3})\right)\cap H^{1}\left(0,T;H^{1}(\mathcal{D}_a,\mathbb{R}^{3\times 3})\right)}\leq C.
\end{align}
Collecting the \pier{bounds}~\eqref{eqn:54}, \eqref{eqn:55bis}, \eqref{eqn:56bis}, \eqref{eqn:57} and \eqref{eqn:58}, which are uniform in $m$ and $\lambda$, from the Banach--Alaoglu, the Aubin--Lions and 
the \pier{Arzel{\`a}--Ascoli} lemmas, 
we finally obtain the convergence properties, up to subsequences
, which we still label by the index $m$ (without reporting the index $\lambda$), as follows:
\begin{align}
\label{eqn:convwm1} & \mathbf{W}_m \overset{\ast}{\rightharpoonup} \mathbf{W} \quad \text{in} \quad L^{\infty}(0,T;H^1(\mathcal{D}_a;Sym(\mathbb{R}^{3\times 3}))),\\
\label{eqn:convwm2} &  \mathbf{W}_m \rightharpoonup \mathbf{W} \quad \text{in} \quad L^{2}(0,T;H^{2}(\mathcal{D}_a;Sym(\mathbb{R}^{3\times 3})))\cap H^{1}(0,T;L^2(\mathcal{D}_a;Sym(\mathbb{R}^{3\times 3}))),\\
\label{eqn:convwm3} & \mathbf{W}_m \to \mathbf{W} \quad \text{in} \quad C^{0}(\pier{[0,T]};L^p(\mathcal{D}_a;Sym(\mathbb{R}^{3\times 3}))) \cap L^{2}(0,T;W^{1,p}(\mathcal{D}_a;Sym(\mathbb{R}^{3\times 3}))), \notag\\
 &\quad\quad\quad\quad\quad\quad\quad \pier{\hbox{with } \, }p\in [1,6), \;\; \text{and} \;\; \text{a.e. in} \; \; \mathcal{D}_{aT},\\
\label{eqn:convwm3bis} & \mathbf{W}_m \to \mathbf{W} \quad \text{in} \quad L^{2}(0,T;C^0(\overline{\mathcal{D}}_a;Sym(\mathbb{R}^{3\times 3}))), 
\end{align}
\begin{align}
\label{eqn:convtm2} &  \boldsymbol{\Theta}_m \rightharpoonup \boldsymbol{\Theta} \quad \text{in} \quad H^{1}(0,T;H^{2}(\mathcal{D}_a;Skew(\mathbb{R}^{3\times 3}))),\\
\label{eqn:convtm3} & \boldsymbol{\Theta}_m \to \boldsymbol{\Theta} \quad \text{in} \quad C^{0}(\pier{[0,T]};W^{1,p}(\mathcal{D}_a;Skew(\mathbb{R}^{3\times 3}))), \;\; p\in [1,6), \;\; \text{and} \;\; \text{a.e. in} \; \; \mathcal{D}_{aT},\\
\label{eqn:convtm3bis} & \displaystyle \boldsymbol{\Theta}_m\to \boldsymbol{\Theta},\; \mathrm{e}^{\pm \boldsymbol{\Theta}_m} \to \mathrm{e}^{\pm \boldsymbol{\Theta}} \quad \text{uniformly in} \; \; \overline{\mathcal{D}}_{aT},
\end{align}
\begin{align}
\label{eqn:convphim1} & \pier{\mathbf{\Phi}}_m \overset{\ast}{\rightharpoonup} \pier{\mathbf{\Phi}} \quad \text{in} \quad L^{\infty}(0,T;H^2(\mathcal{D}_a;\mathbb{R}^{3})),\\
\label{eqn:convphim2} & \pier{\mathbf{\Phi}}_m \rightharpoonup \pier{\mathbf{\Phi}} \quad \text{in} \quad L^{2}(0,T;H^{3}(\mathcal{D}_a;\mathbb{R}^{3}))\cap H^{1}(0,T;H^{1}(\mathcal{D}_a;\mathbb{R}^{3})),\\
\label{eqn:convphim3} & \pier{\mathbf{\Phi}}_m \to \pier{\mathbf{\Phi}} \quad \text{in} \quad C^{0}(\pier{[0,T]};W^{1,p}(\mathcal{D}_a;\mathbb{R}^{3}))\cap L^{2}(0,T;W^{2,p}(\mathcal{D}_a;\mathbb{R}^{3})), \notag \\
&\quad\quad\quad\quad\quad\quad \pier{\hbox{with } \, }p\in [1,6), \;\; \text{and} \;\; \text{a.e. in} \; \; \mathcal{D}_{aT},
\end{align}
\begin{align}
\label{eqn:convzm1} & \mathbf{Z}_m \overset{\ast}{\rightharpoonup} \mathbf{Z} \quad \text{in} \quad L^{\infty}(0,T;H^2(\mathcal{D}_a;\mathbb{R}^{3\times 3})),\\
\label{eqn:convzm2} & \mathbf{Z}_m \rightharpoonup \mathbf{Z} \quad \text{in} \quad L^{2}(0,T;H^{3}(\mathcal{D}_a;\mathbb{R}^{3\times 3}))\cap H^{1}(0,T;H^{1}(\mathcal{D}_a;\mathbb{R}^{3\times 3})),\\
\label{eqn:convzm3} & \mathbf{Z}_m \to \mathbf{Z} \quad \text{in} \quad C^{0}(\pier{[0,T]};W^{1,p}(\mathcal{D}_a;\mathbb{R}^{3\times 3}))\cap L^{2}(0,T;W^{2,p}(\mathcal{D}_a;\mathbb{R}^{3\times 3})), \notag \\
&\quad\quad\quad\quad\quad\quad \pier{\hbox{with } \, }p\in [1,6), \;\; \text{and} \;\; \text{a.e. in} \; \; \mathcal{D}_{aT},
\end{align}
\begin{align}
\label{eqn:convsigmam1} & \boldsymbol{\Sigma}_m \overset{\ast}{\rightharpoonup} \boldsymbol{\Sigma} \quad \text{in} \quad L^{\infty}(0,T;L^2(\mathcal{D}_a;\mathbb{R}^{3\times 3})),\\
\label{eqn:convsigmam2} & \mathcal{G}_{L,\diver}(\curl \boldsymbol{\Sigma}_m) \overset{\ast}{\rightharpoonup} \mathcal{G}_{L,\diver}(\curl \boldsymbol{\Sigma}) \quad \text{in} \quad L^{\infty}(0,T;H^1(\mathcal{D}_a;\mathbb{R}^{3\times 3})),
\end{align}
as \pier{$m\to \infty$}. We note that \eqref{eqn:convwm3bis} follows from \pier{\eqref{eqn:convwm3} and the compact embedding
$$ W^{1,p}(\mathcal{D}_a;Sym(\mathbb{R}^{3\times 3}))\subset C^0(\overline{\mathcal{D}}_a;Sym(\mathbb{R}^{3\times 3})),$$
holding for $p>3$.
Moreover, 
as 
$$ H^2(\mathcal{D}_a;Skew(\mathbb{R}^{3\times 3})) \, \hbox{ is compactly embedded into } \,C^0(\overline{\mathcal{D}}_a;Skew (\mathbb{R}^{3\times 3})),$$ 
the convergence \eqref{eqn:convtm2} implies a strong convergence in $C^{0}(\pier{[0,T]};C^0(\overline{\mathcal{D}}_a;Skew (\mathbb{R}^{3\times 3})))$, whence \eqref{eqn:convtm3bis} is easily deduced, thanks the continuity of the exponential operator as well.}

With the convergence results \eqref{eqn:convwm1}--\eqref{eqn:convsigmam2}, we can pass to the limit in \pier{the system} \eqref{eqn:48m2} in a first step as $m\to \infty$. Let's take $\nwhat{\mathbf{W}}_m=PS_m(\nwhat{\mathbf{W}})$ and $\nwhat{\boldsymbol{\Omega}}_m=PA_m(\nwhat{\boldsymbol{\Omega}})$, with arbitrary $\nwhat{\mathbf{W}}\in L^2(\mathcal{D}_a;Sym(\mathbb{R}^{3\times 3}))$, $\nwhat{\boldsymbol{\Omega}}\in L^2(\mathcal{D}_a;Skew(\mathbb{R}^{3\times 3}))$, multiply the first two equations by $\omega \in C_c^{\infty}([0,T])$ and integrate over the time interval $[0,T]$. This gives
\begin{equation}
\label{eqn:48m2cont1}
\begin{cases}
\displaystyle \int_0^T\omega \int_{\mathcal{D}_a}\rm{Sym}\left(\mathrm{e}^{-\boldsymbol{\Theta}_m}\curl\left(\mathcal{G}_{L,\diver}\ast \left(
\curl \boldsymbol{\Sigma}_m\right)\right)\right)\colon \nwhat{\mathbf{W}}_m\\
\displaystyle + \int_0^T\omega \int_{\mathcal{D}_a}\left(\mathbf{W}_m-\mathbf{I}+\partial I_{C_{\alpha}}^{\lambda}(\mathbf{W}_m)+\overbigdot{\mathbf{W}}_m\right)\colon \nwhat{\mathbf{W}}_m\\
\displaystyle + \int_0^T\omega \int_{\mathcal{D}_a}\pier{\grad{}} \mathbf{W}_m::\pier{\grad{}}\nwhat{\mathbf{W}}_m=\int_0^T\omega\int_{\mathcal{D}_a}\mathbf{W}_{\text{ext}}(\mathbf{W}_m,t)\colon \nwhat{\mathbf{W}}_m,\\ \\
\displaystyle \int_0^T\omega\int_{\mathcal{D}_a}\rm{Skew}\left(\curl\left(\mathcal{G}_{L,\diver}\ast
\left(\curl \boldsymbol{\Sigma}_m\right)\right)\mathbf{W}_m\mathrm{e}^{-\boldsymbol{\Theta}_m}\right)\colon 
\nwhat{\boldsymbol{\Omega}}_m 
\\
\displaystyle +  \int_0^T\frac{\omega}{2}\int_{\mathcal{D}_a}
\pier{ \Bigl(\pier{\grad{}} 
\boldsymbol{\Theta}_m + \grad \overbigdot{\boldsymbol{\Theta}}_m\Bigr)} ::
\pier{\grad{}} \nwhat{\boldsymbol{\Omega}}_m
=\int_0^T\omega\int_{\mathcal{D}_a}\boldsymbol{\Omega}_{\text{ext}}\abramo{(\boldsymbol{\Omega}_m,t)}\colon \nwhat{\boldsymbol{\Omega}}_m.
\end{cases}
\end{equation}
We observe that 
\begin{equation}
\label{eqn:48m2cont2}
\begin{cases}
PS_m(\nwhat{\mathbf{W}})\to \nwhat{\mathbf{W}} \quad \text{in} \; \;L^2(\mathcal{D}_a;Sym(\mathbb{R}^{3\times 3})),
\\[2mm]
PA_m(\nwhat{\boldsymbol{\Omega}})\to \nwhat{\boldsymbol{\Omega}}\quad \text{in} \; \;L^2(\mathcal{D}_a;Skew(\mathbb{R}^{3\times 3})),
\end{cases}
\end{equation}
as $m\to \infty$.
Thanks to \eqref{eqn:convtm3bis} and \eqref{eqn:48m2cont2}$_1$, we have that $$\mathrm{e}^{\boldsymbol{\Theta}_m}\nwhat{\mathbf{W}}_m\to \mathrm{e}^{\boldsymbol{\Theta}}\nwhat{\mathbf{W}}\ \hbox{ in } \ L^{\infty}(0,T;L^{2}(\mathcal{D}_a,\mathbb{R}^{3\times 3})).$$ Hence, using \eqref{eqn:convsigmam2}, by the product of weak-strong convergence we have that
\begin{align}
\label{lim1}
 &\notag\int_0^T\omega \int_{\mathcal{D}_a}\rm{Sym}\left(\mathrm{e}^{-\boldsymbol{\Theta}_m}\curl\left(\mathcal{G}_{L,\diver}\ast \left(
\curl \boldsymbol{\Sigma}_m\right)\right)\right)\colon \nwhat{\mathbf{W}}_m\\
&\notag
=\int_0^T\omega \int_{\mathcal{D}_a}\curl\left(\mathcal{G}_{L,\diver}\ast \left(
\curl \boldsymbol{\Sigma}_m\right)\right)\colon \mathrm{e}^{\boldsymbol{\Theta}_m}\nwhat{\mathbf{W}}_m\\
& \to \int_0^T\omega \int_{\mathcal{D}_a}\rm{Sym}\left(\mathrm{e}^{-\boldsymbol{\Theta}}\curl\left(\mathcal{G}_{L,\diver}\ast \left(
\curl \boldsymbol{\Sigma}\right)\right)\right)\colon \nwhat{\mathbf{W}},
\end{align}
as $m\to \infty$. \pier{Owing to \eqref{eqn:convwm3} and the Lipschitz continuity of $\partial \psi_{C_\alpha}^{\lambda}$, it turns out that 
$\boldsymbol{\chi}_{\alpha,m}{} = \partial \psi_{C_\alpha}^{\lambda}\left(\mathbf{W}_m\right)$ strongly converges to 
$\boldsymbol{\chi}_{\alpha} := \partial \psi_{C_\alpha}^{\lambda}\left(\mathbf{W}\right)$ say in $C^{0}(\pier{[0,T]};L^2(\mathcal{D}_a;\mathbb{R}^{3\times 3}))$.
Then, on account of \eqref{eqn:convwm1}--\eqref{eqn:convwm3}, \eqref{eqn:48m2cont2}$_1$ and Assumption {\bf A3},} we readily obtain that
\begin{align*}
 &\int_0^T\omega \int_{\mathcal{D}_a}\left(\mathbf{W}_m-\mathbf{I}+\pier{\boldsymbol{\chi}_{\alpha,m}}+\overbigdot{\mathbf{W}}_m\right)\colon \nwhat{\mathbf{W}}_m+\int_0^T\omega \int_{\mathcal{D}_a}\pier{\grad{}} \mathbf{W}_m::\pier{\grad{}}\nwhat{\mathbf{W}}_m\\
 & \to \int_0^T\omega \int_{\mathcal{D}_a}\left(\mathbf{W}-\mathbf{I}+\pier{\boldsymbol{\chi}_{\alpha}}+\overbigdot{\mathbf{W}}\right)\colon \nwhat{\mathbf{W}}-\int_0^T\omega \int_{\mathcal{D}_a}\Delta \mathbf{W}:\nwhat{\mathbf{W}},
\end{align*}
and 
\begin{equation*}
 \int_0^T\omega\int_{\mathcal{D}_a}\mathbf{W}_{\text{ext}}(\mathbf{W}_m,t)\colon \nwhat{\mathbf{W}}_m\to \int_0^T\omega\int_{\mathcal{D}_a}\mathbf{W}_{\text{ext}}(\mathbf{W},t)\colon \nwhat{\mathbf{W}},
\end{equation*}
as $m\to \infty$. Moroever, thanks to \eqref{eqn:convtm3bis}, \eqref{eqn:convwm3bis} and \eqref{eqn:48m2cont2}$_2$, we have that $\nwhat{\boldsymbol{\Omega}}_m \mathbf{W}_m\mathrm{e}^{\boldsymbol{\Theta}_m}\to \nwhat{\boldsymbol{\Omega}} \mathbf{W}\mathrm{e}^{\boldsymbol{\Theta}}$ in $L^{2}(0,T;L^{2}(\mathcal{D}_a,\mathbb{R}^{3\times 3}))$. Hence, using \eqref{eqn:convsigmam2}, by the product of weak-strong convergence and with similar calculations as in \eqref{lim1} we have that
\begin{align}
\label{lim2}
 & \notag \int_0^T\omega\int_{\mathcal{D}_a}\rm{Skew}\left(\curl\left(\mathcal{G}_{L,\diver}\ast
\left(\curl \boldsymbol{\Sigma}_m\right)\right)\mathbf{W}_m\mathrm{e}^{-\boldsymbol{\Theta}_m}\right)\colon \nwhat{\boldsymbol{\Omega}}_m \\
& \to \int_0^T\omega\int_{\mathcal{D}_a}\rm{Skew}\left(\curl\left(\mathcal{G}_{L,\diver}\ast
\left(\curl \boldsymbol{\Sigma}\right)\right)\mathbf{W}\mathrm{e}^{-\boldsymbol{\Theta}}\right)\colon \nwhat{\boldsymbol{\Omega}},
\end{align}
as $m\to \infty$. Thanks to \pier{\eqref{eqn:convtm2}--\eqref{eqn:convtm3}} it is easy to deduce that
\begin{align*}
& \int_0^T
\frac{\omega}{2}\int_{\mathcal{D}_a}\pier{\grad{}} \boldsymbol{\Theta}_m::\pier{\grad{}} \nwhat{\boldsymbol{\Omega}}_m+\int_0^T\frac{\omega}{2}\int_{\mathcal{D}_a}\pier{\grad{}} \overbigdot{\boldsymbol{\Theta}}_m::\pier{\grad{}} \nwhat{\boldsymbol{\Omega}}_m\\
& \to -\int_0^T
\frac{\omega}{2}\int_{\mathcal{D}_a}\Delta \boldsymbol{\Theta}: \nwhat{\boldsymbol{\Omega}}-\int_0^T\frac{\omega}{2}\int_{\mathcal{D}_a}\Delta \overbigdot{\boldsymbol{\Theta}}: \nwhat{\boldsymbol{\Omega}}
\end{align*}
and
\[
 \int_0^T\omega\int_{\mathcal{D}_a}\boldsymbol{\Omega}_{\text{ext}}\abramo{(\boldsymbol{\Omega}_m,t)}\colon \nwhat{\boldsymbol{\Omega}}_m\to \int_0^T\omega\int_{\mathcal{D}_a}\boldsymbol{\Omega}_{\text{ext}}\abramo{(\boldsymbol{\Omega},t)}\colon \nwhat{\boldsymbol{\Omega}}.
\]
For what concerns \pier{the second equality in} \eqref{eqn:48m2}$_3$, thanks to the convexity of $\psi_D^{\lambda}$ we \pier{can express it as 
\begin{equation}
 \label{lim3}
 \int_0^T\!\!\int_{\mathcal{D}_a}(\curl \mathbf{X}-\curl \mathbf{Z}_m):\boldsymbol{\Sigma}_m+\int_0^T\!\!\int_{\mathcal{D}_a}\psi_D^{\lambda}(\curl \mathbf{Z}_m)\leq \int_0^T\!\!\int_{\mathcal{D}_a}\psi_D^{\lambda}(\curl \mathbf{X}),
\end{equation}
for all $\mathbf{X}\in L^2(0,T;H^1(\mathcal{D}_a;\mathbb{R}^{3\times 3}))$. Given the convergence results \eqref{eqn:convsigmam1} and \eqref{eqn:convzm3}, as $\boldsymbol{\Sigma}_m $ weakly converges to 
$\boldsymbol{\Sigma}$ in $L^2(0,T;L^2(\mathcal{D}_a;\mathbb{R}^{3\times 3}))$ and $\psi_D^{\lambda}$ is Lipschitz continuous, we get in the limit that   
\begin{equation}
 \label{lim3bis}
 \int_0^T\!\!\int_{\mathcal{D}_a}(\curl \mathbf{X}-\curl \mathbf{Z}):\boldsymbol{\Sigma}+\int_0^T\!\!\int_{\mathcal{D}_a}\psi_D^{\lambda}(\curl \mathbf{Z})\leq \int_0^T\!\!\int_{\mathcal{D}_a}\psi_D^{\lambda}(\curl \mathbf{X}),
\end{equation}
as $m\to \infty$, for all $\mathbf{X}\in L^2(0,T;H^1(\mathcal{D}_a;\mathbb{R}^{3\times 3}))$}. Finally, we want to 
pass to the limit in \eqref{eqn:48m2}$_4$ and \eqref{eqn:48m2}$_5$ as $m\to \infty$. In order to do so, we observe 
that, since \abramo{$\mathrm{e}^{\boldsymbol{\Theta}_m}\to \mathrm{e}^{\boldsymbol{\Theta}}$ a.e. in $\mathcal{D}_{aT}$, $\grad \boldsymbol{\Theta}_m\to \grad \boldsymbol{\Theta}$ in $C^{0}\bigl(0,T;L^p\bigl(\mathcal{D}_a;\mathbb{R}^{3\times 3\times 3}\pier{\bigr)}\pier{\bigr)}$ and a.e. in $\mathcal{D}_{aT}$ and since $\mathrm{e}^{\boldsymbol{\Theta}_m}$ is uniformly bounded, a generalized form of the Lebesgue convergence theorem gives that}
\begin{equation}
 \label{lim4}
 \pier{\grad{}} \mathrm{e}^{\boldsymbol{\Theta}_m}\to \pier{\grad{}} \mathrm{e}^{\boldsymbol{\Theta}}\quad \text{in} \quad C^{0}\bigl(0,T;L^p\bigl(\mathcal{D}_a;\mathbb{R}^{3\times 3\times 3}\pier{\bigr)}\pier{\bigr)}, \;\; p\in [1,6).
\end{equation}
Hence, using \eqref{lim4}, \eqref{eqn:convwm3}, \eqref{eqn:convwm3bis} and \eqref{eqn:convtm3bis} we can prove that
\begin{align}
 \label{lim5}
 & \diver\pier{\bigl(}\mathrm{e}^{\boldsymbol{\Theta}_m}\mathbf{W}_m\pier{\bigr)}=\pier{\bigl(}\pier{\grad{}} \mathrm{e}^{\boldsymbol{\Theta}_m}\pier{\bigr)}\mathbf{W}_m+\mathrm{e}^{\boldsymbol{\Theta}_m}\diver \mathbf{W}_m \notag \\
 & \to \diver\pier{\bigl(}\mathrm{e}^{\boldsymbol{\Theta}}\mathbf{W}\pier{\bigr)}\quad \text{in} \quad C^{0}\pier{\bigl(}0,T;L^{\frac{p}{2}}\pier{\bigl(}\mathcal{D}_a;\mathbb{R}^{3}\pier{\bigr)}\pier{\bigr)}\cap L^2\pier{\bigl(}0,T;L^{p}\pier{\bigl(}\mathcal{D}_a;\mathbb{R}^{3}\pier{\bigr)}\pier{\bigr)}, \;\; p\in [1,6),
 \end{align}
 and analogously
 \begin{align}
 \label{lim6}
 & \curl \pier{\bigl(}\mathrm{e}^{\boldsymbol{\Theta}_m}\mathbf{W}_m\pier{\bigr)}=\pier{\bigl(}\curl \mathrm{e}^{\boldsymbol{\Theta}_m}\pier{\bigr)}\mathbf{W}_m+\boldsymbol{\epsilon} \hskip1pt \mathrm{e}^{\boldsymbol{\Theta}_m}\pier{\grad{}} \mathbf{W}_m\notag \\
 & \to \curl\pier{\bigl(}\mathrm{e}^{\boldsymbol{\Theta}}\mathbf{W}\pier{\bigr)}\quad \text{in} \quad C^{0}\pier{\bigl(}0,T;L^{\frac{p}{2}}\pier{\bigl(}\mathcal{D}_a;\mathbb{R}^{3\times 3}\pier{\bigr)}\pier{\bigr)}\cap L^2\pier{\bigl(}0,T;L^{p}\pier{\bigl(}\mathcal{D}_a;\mathbb{R}^{3\times 3}\pier{\bigr)}\pier{\bigr)}, \;\; p\in [1,6),
 \end{align}
where $\boldsymbol{\epsilon}$ is the Ricci tensor. With the strong convergence results \eqref{lim5}, \eqref{lim6}, \eqref{eqn:convphim3} and \eqref{eqn:convzm3} we can straightforwardly pass to the limit in \eqref{eqn:48m2}$_4$ and \eqref{eqn:48m2}$_5$ as $m\to \infty$. Collecting all the previous results, we obtain the following limit system, as $m\to \infty$, \pier{in terms of the limit functions that will be now denoted by $\mathbf{W}^{\lambda}, \, \boldsymbol{\Theta}^{\lambda}, \, \boldsymbol{\chi}_{\alpha}^\lambda, \,\boldsymbol{\Sigma}^{\lambda},\, \pier{\mathbf{\Phi}}^{\lambda},\, \mathbf{Z}^{\lambda} $. Here, it is:}
\begin{equation}
\label{eqn:48lambda}
\begin{cases}
\int_{\mathcal{D}_a}\rm{Sym}\pier{\bigl(}\mathrm{e}^{-\boldsymbol{\Theta}^{\lambda}}\curl\pier{\bigl(}\mathcal{G}_{L,\diver}\ast \pier{\bigl(}
\curl \boldsymbol{\Sigma}^{\lambda}\pier{\bigr)}\pier{\bigr)}\pier{\bigr)}\colon \nwhat{\mathbf{W}}\\
+ \int_{\mathcal{D}_a}\pier{\bigl(}\mathbf{W}^{\lambda}-\mathbf{I}+\pier{{}\boldsymbol{\chi}_{\alpha}^\lambda}+{\overbigdot{\mathbf{W}}}{}^\lambda\pier{\bigr)}\colon \nwhat{\mathbf{W}}\\
-\int_{\mathcal{D}_a}\Delta \mathbf{W}^{\lambda}:\nwhat{\mathbf{W}}=\int_{\mathcal{D}_a}\mathbf{W}_{\text{ext}}(\mathbf{W}^{\lambda},t)\colon \nwhat{\mathbf{W}},\\ \\
\int_{\mathcal{D}_a}\rm{Skew}\pier{\bigl(}\curl\pier{\bigl(}\mathcal{G}_{L,\diver}\ast
\pier{\bigl(}\curl \boldsymbol{\Sigma}^{\lambda}\pier{\bigr)}\pier{\bigr)}\mathbf{W}^{\lambda}\mathrm{e}^{-\boldsymbol{\Theta}^{\lambda}}\pier{\bigr)}\colon \nwhat{\boldsymbol{\Omega}}\\ -
\frac{1}{2}\int_{\mathcal{D}_a}\Delta \boldsymbol{\Theta}^{\lambda}: \nwhat{\boldsymbol{\Omega}}-
\frac{1}{2}\int_{\mathcal{D}_a}\Delta {\overbigdot{\boldsymbol{\Theta}}}{}^{\lambda}: \nwhat{\boldsymbol{\Omega}}=\int_{\mathcal{D}_a}\boldsymbol{\Omega}_{\text{ext}}\abramo{(\boldsymbol{\Theta}^{\lambda},t)}\colon \nwhat{\boldsymbol{\Omega}},\\ \\
\pier{\pier{{}\boldsymbol{\chi}_{\alpha}^\lambda{} = \partial \psi_{C_\alpha}^{\lambda}(\mathbf{W}^\lambda),}\quad
\boldsymbol{\Sigma}^\lambda=\partial \psi_D^{\lambda}(\curl \mathbf{Z}^\lambda),}\\ \\
\Delta \pier{\mathbf{\Phi}}^{\lambda}=\diver\pier{\bigl(}\mathrm{e}^{\boldsymbol{\Theta}^{\lambda}}\mathbf{W}^{\lambda}\pier{\bigr)}, 
\\ \\
-P_L\Delta {\mathbf{Z}^{\lambda}}=\curl \pier{\bigl(}\mathrm{e}^{\boldsymbol{\Theta^{\lambda}}}\mathbf{W}^{\lambda}\pier{\bigr)}, 
\end{cases}
\end{equation}
for a.e. $t \in [0,T]$, for \pier{all choices of} $\nwhat{\mathbf{W}} \in L^2(\mathcal{D}_a;Sym(\mathbb{R}^{3\times 3}))$, $\nwhat{\boldsymbol{\Omega}} \in L^2(\mathcal{D}_a;Skew(\mathbb{R}^{3\times 3}))$, 
and with initial conditions \pier{(cf. the assumption \textbf{A2} and \eqref{eqn:48mic})}
\begin{equation}
    \label{eqn:48lic}
   \pier{ \mathbf{W}^\lambda(\cdot,0)=\mathbf{W}_0, \quad \boldsymbol{\Theta}^\lambda (\cdot,0)=\boldsymbol{0} \quad  \hbox{in }
   \, \mathcal{D}_a.}
\end{equation}
In \pier{the system} \eqref{eqn:48lambda} we have restored the index $\lambda$, to indicate the dependence of the solutions from the regularization parameter $\lambda$. We observe, without reporting all the details, that the estimates \eqref{eqn:54}, \eqref{eqn:542}, \eqref{eqn:543}, \eqref{eqn:545}, \eqref{eqn:55bis}, \eqref{eqn:56bis}, 
\eqref{eqn:57} and \eqref{eqn:58} are preserved in the limit as $m\to \infty$, \pier{i.e.,} they are valid for the solutions of \pier{the system} \eqref{eqn:48lambda}. \pier{This allows us to pass to the limit as $\lambda \to 0$, up to subsequences of $\lambda$, in \pier{the system} \eqref{eqn:48lambda}, with similar calculations as the ones employed for the study of the limit problem as $m\to \infty$. On the other hand, since by comparison in \eqref{eqn:48lambda}$_1$ we obtain that 
\begin{equation}
\label{eqn:comparison}
\pier{\boldsymbol{\chi}_{\alpha}^\lambda} \quad \text{is bounded in} \; \; L^2(0,T;L^2(\mathcal{D}_a;Sym(\mathbb{R}^{3\times 3}))),
\end{equation}
uniformly with respect to $\lambda$, consequently 
$$\hbox{$\pier{\boldsymbol{\chi}_{\alpha}^\lambda}\,$  will converge weakly to some $\,\pier{\boldsymbol{\chi}_{\alpha}}\,$  in $\,L^2(0,T;L^2(\mathcal{D}_a;Sym(\mathbb{R}^{3\times 3})))$}$$ 
as $\lambda \to 0$ along a subsequence. This weak convergence, combined with the strong convergence of $\mathbf{W}^\lambda$ to $\mathbf{W}$ in the same space $L^2(0,T;L^2(\mathcal{D}_a;Sym(\mathbb{R}^{3\times 3})))$, and the maximal monotonicity of the subdifferential operator $\partial \psi_{C_\alpha}$ enable us to prove that ${}\boldsymbol{\chi}_{\alpha}{} \in \partial \psi_{C_\alpha} (\mathbf{W})$ a.e. in $\mathcal{D}_{aT}$. Similar considerations can be done for the proof of the other inclusion
$ \boldsymbol{\Sigma} \in \partial \psi_D (\curl \mathbf{Z} )$; these are usual arguments in the framework of the theory of maximal monotone operators, see e.g. \cite[Lemma~2.3, p.~38]{barbu}.
Therefore, passing to the limit as $\lambda \to 0$ in \pier{the system} \eqref{eqn:48lambda}, we obtain that the limit point is a solution of
\eqref{eqn:48symskew}--\eqref{eqn:49symskew} 
with initial conditions \eqref{eqn:50symskew}, and regularity given by \eqref{eqn:57thm}--\eqref{eqn:60thm}. Moreover, by lower semicontinuity we obtain in the limit 
as $\lambda \to 0$ that
\[
\sup_{t\in (0,T)} \int_{\mathcal{D}_a}I_{C_{\alpha}}\pier{(\mathbf{W} (\cdot, t))}\leq C,
\]
whence, due to \eqref{eqn:57thm} as well, we have that $\mathbf{W}(\mathbf{a},t)\in \hbox{SPD}_{\alpha}$ for all $\mathbf{a}\in \mathcal{D}_a$ and a.a. $t\in (0,T)$.}
\end{pf}
\abramohhb{
\begin{rem}
The property that $\mathbf{W}(\mathbf{a},t)\in \hbox{SPD}_{\alpha}$ for all $\mathbf{a}\in \mathcal{D}_a$ and a.a. $t\in (0,T)$, proved in the previous Theorem, implies that the material experiences neither flattening nor crushing and that (see \cite[Theorem 5.5.1]{ciarlet2}) a point which is in the interior of the domain remains inside the domain during the evolution, for a.a. $t\in (0,T)$.
\end{rem}
}

\section{The limiting case}
\label{sec:limit}
In this section we study the limit system of \abramonew{\eqref{eqn:48symskew}--\eqref{eqn:49symskew}} as $k\to 0$. As we will see, in this case the solution of the limit system is unique and continuously depends on initial data. The drawback is that in this case the incompatibility in the system dynamics is always active, contrarily to what happens in the case with $k>0$, as observed in Remark \ref{rem:remsigma}.

In the case $k=0$, from \abramonew{\eqref{eqn:48symskew}$_3$} we have that $\boldsymbol{\Sigma}=\curl \mathbf{Z}$, hence \pier{the system} \abramonew{\eqref{eqn:48symskew}} becomes
\abramonew{
\begin{equation}
\label{eqn:48k0}
\begin{cases}
\displaystyle \text{Sym}\left(\mathrm{e}^{-\boldsymbol{\Theta}}\curl \left(\mathbf{Z}\right)\right)+\mathbf{W}-\mathbf{I}+\pco{{}
\boldsymbol{\chi}_{\alpha}}+\overbigdot{\mathbf{W}}-\Delta \mathbf{W}\pco{{}={}}\mathbf{W}_{\text{ext}}(\mathbf{W},t),\\ \\
\displaystyle \text{Skew}\left(\curl\left(\mathbf{Z}\right)\mathbf{W}\mathrm{e}^{-\boldsymbol{\Theta}}\right)-\frac{1}{2}\Delta \boldsymbol{\Theta}-\frac{1}{2}\Delta \overbigdot{\boldsymbol{\Theta}}=\boldsymbol{\Omega}_{\text{ext}}\abramo{(\boldsymbol{\Theta},t)},\\ \\
\displaystyle \pco{\boldsymbol{\chi}_{\alpha} \in \partial \pier{\nwhat{\psi}(\mathbf{W})},}\\ \\
\displaystyle \Delta \pier{\mathbf{\Phi}}=\diver\left(\mathrm{e}^{\boldsymbol{\Theta}}\mathbf{W}\right),\\ \\
\displaystyle - P_L\Delta {\mathbf{Z}}=\curl \left(\mathrm{e}^{\boldsymbol{\Theta}}\mathbf{W}\right),
\end{cases}
\end{equation}
}
with boundary conditions
\begin{equation}
\label{eqn:49k0}
\begin{cases}
\displaystyle\mathbf{W}=\mathbf{I},\;\boldsymbol{\Theta}=\overbigdot{ \boldsymbol{\Theta}}=\boldsymbol{0} \quad \text{on} \; \Gamma_a\times (0,T),\\ 
\displaystyle  \pier{\mathbf{\Phi}(\mathbf{a}, t) =\mathbf{a}, \; \mathbf{Z}(\mathbf{a}, t)=\mathbf{0} \quad \text{for} \; (\mathbf{a}, t) \in \Gamma_a\times (0,T)}
\end{cases}
\end{equation}
 and initial conditions
\begin{equation}
\label{eqn:50k0}
\mathbf{W}(\mathbf{a},0)=\mathbf{W}_0\pco{(\mathbf{a})}, \; \boldsymbol{\Theta}(\mathbf{a},0)=\boldsymbol{0} \pier{ \quad \text{for} \; \mathbf{a}  \in \mathcal{D}_a}.
\end{equation}

\begin{rem}
 \label{rem:3}
 The variable $\mathbf{Z}$ in \pier{the system} \eqref{eqn:48k0} may be interpreted as the Lagrange multiplier of the compatibility condition $\curl \left(\mathrm{e}^{\boldsymbol{\Theta}}\mathbf{W}\right)=\boldsymbol{0}$, with the addition of an elliptic regularization of the constraint given by the term $- P_L\Delta {\mathbf{Z}}$ in 
 \pco{\eqref{eqn:48k0}$_5$}.
 Indeed, \pier{the system} \eqref{eqn:48k0} may be obtained from the principle of virtual power \eqref{eqn:12} and the dissipative \pco{equality~\eqref{eqn:21}} by enforcing in the expression of the Free Energy~\eqref{eqn:22} the compatibility constraint through a Lagrange multiplier,~\pier{i.e.,}
 \begin{align}
\label{eqn:22z}
\displaystyle \notag &\psi(\mathbf{W},\mathbf{R},\mathbf{Z}):=\frac{1}{2}\pier{\|}\mathbf{W}-\mathbf{I}\pier{\|}^2+\nwhat{\psi}(\mathbf{W})+\frac{1}{2}\pier{\|}\grad \mathbf{W}\pier{\|}^2+\frac{1}{2}\pier{\|}\grad \mathbf{R}\pier{\|}^2+\int_{\mathcal{D}_a}\mathbf{Z}:\curl\left(\mathrm{e}^{\boldsymbol{\Theta}}\mathbf{W}\right).
\end{align}  
\pco{Setting}
 \[
  \mathcal{F}(\mathbf{Z},\mathbf{W},\boldsymbol{\Theta}):=\int_{\mathcal{D}_a}\mathbf{Z}:\curl\left(\mathrm{e}^{\boldsymbol{\Theta}}\mathbf{W}\right),
 \]
 we observe that
 \begin{align*}
  & \left(\frac{\delta \mathcal{F}}{\delta \mathbf{W}},\delta \mathbf{W}\right)=\left(\mathbf{Z},\curl\left(\mathrm{e}^{\boldsymbol{\Theta}}\delta \mathbf{W}\right)\right)=\abramonew{\left(\text{Sym}\left(\mathrm{e}^{-\boldsymbol{\Theta}}\curl \left(\mathbf{Z}\right)\right),\delta \mathbf{W}\right)},\\
  & \left(\frac{\delta \mathcal{F}}{\delta \boldsymbol{\Theta}},\delta \boldsymbol{\Theta}\right)=\left(\mathbf{Z},\curl\left(\delta \boldsymbol{\Theta} \mathrm{e}^{\boldsymbol{\Theta}}\mathbf{W}\right)\right)=\abramonew{\left(\text{Skew}\left(\curl\left(\mathbf{Z}\right)\mathbf{W}\mathrm{e}^{-\boldsymbol{\Theta}}\right),\delta \boldsymbol{\Theta} \right)},\\
  & \left(\frac{\delta \mathcal{F}}{\delta \mathbf{Z}},\delta \mathbf{Z}\right)=\left(\curl\left(\mathrm{e}^{\boldsymbol{\Theta}}\mathbf{W}\right),\delta \mathbf{Z}\right).
 \end{align*}
Moreover, substituting \pier{\eqref{eqn:48k0}$_5$} with the relation $ -\epsilon P_L\Delta {\mathbf{Z}}=\curl \left(\mathrm{e}^{\boldsymbol{\Theta}}\mathbf{W}\right)$, with $0<\epsilon <<1$, \pier{the system} \eqref{eqn:48k0} may be interpreted as a system with \pco{a} penalization of the compatibility condition. 
\end{rem}

We give for \pier{the system} \eqref{eqn:48k0} the following existence and regularity result.

\begin{thm}
 \label{thm:k0}
 Let assumptions \textbf{A1}-\textbf{A3} be satisfied.  Then, for any $T>0$ there is a \pco{quintuplet
 $(\mathbf{W},\boldsymbol{\Theta},\pier{{}\boldsymbol{\chi}_{\alpha}{}}, \pier{\mathbf{\Phi}},\mathbf{Z})$},
 with
\begin{align}
\label{eqn:57thmk0}
&\mathbf{W}\in L^{\infty}(0,T;H^1(\mathcal{D}_a;Sym(\mathbb{R}^{3\times 3})))  \notag \\
&\qquad{}\cap H^1(0,T;L^2(\mathcal{D}_a;Sym(\mathbb{R}^{3\times 3})))\cap L^2(0,T;H^2(\mathcal{D}_a;Sym(\mathbb{R}^{3\times 3}))),
\end{align}
and $\mathbf{W}(\mathbf{a},t)\in $ \pier{{}\rm\text{SPD}}$_{\alpha}$ for a.e. $(\mathbf{a},t)\in \mathcal{D}_{aT}$,  
\begin{equation}
\label{eqn:58thmk0}
\boldsymbol{\Theta}\in H^1(0,T;H^2(\mathcal{D}_a;Skew(\mathbb{R}^{3\times 3}))),
\end{equation}
\begin{equation}
 \label{eqn:58tristhm0}
 \pier{{}\boldsymbol{\chi}_{\alpha}{}}\in L^{2}(0,T;L^2(\mathcal{D}_a;\mathbb{R}^{3\times 3})),
\end{equation}
\begin{equation}
\label{eqn:59thmk0}
\pier{\mathbf{\Phi}\in L^{\infty}(0,T;H^2(\mathcal{D}_a,\mathbb{R}^{3})\cap H^1(\mathcal{D}_a;\mathbb{R}^{3}))\cap L^{2}(0,T;H^3(\mathcal{D}_a;\mathbb{R}^{3})),}
\end{equation}
\begin{equation}
\label{eqn:60thmk0}
\mathbf{Z}\in L^{\infty}(0,T;H^2(\mathcal{D}_a;\mathbb{R}^{3\times 3})\cap H_{0,\diver}^1(\mathcal{D}_a,\mathbb{R}^{3\times 3}))\cap L^{2}(0,T;H^3(\mathcal{D}_a;\mathbb{R}^{3})),
\end{equation}
which solves \pier{the system} \eqref{eqn:48k0}--\eqref{eqn:49k0} for a.e. $\mathcal{D}_{aT}$ with initial conditions~\pco{\eqref{eqn:50k0}}. Moreover, the solution is unique and the following continuous dependence \pco{result holds: given} two solutions $(\mathbf{W}_1,\boldsymbol{\Theta}_1,\pco{{}\boldsymbol{\chi}_{\alpha,1}{}},\pier{\mathbf{\Phi}}_1,\mathbf{Z}_1)$, corresponding to the initial data $(\mathbf{W}_1^0,\boldsymbol{\Theta}_1^0)$, and $(\mathbf{W}_2,\boldsymbol{\Theta}_2,\pco{{}\boldsymbol{\chi}_{\alpha,2}{}},\pier{\mathbf{\Phi}}_2,\mathbf{Z}_2)$, corresponding to the initial data $(\mathbf{W}_2^0,\boldsymbol{\Theta}_2^0)$, there exists a constant $C$ depending only on $\mathcal{D}_a$ such that
\begin{align}
 \label{eqn:contdep}
 & \notag \frac{1}{2}\pier{\|}(\mathbf{W}_1-\mathbf{W}_2)\pco{(t)}\pier{\|}^2+\frac{1}{4}\pier{\|}\grad(\boldsymbol{\Theta}_1-\boldsymbol{\Theta}_2)\pco{(t)}\pier{\|}^2+\int_0^T\biggl(\pier{\|}\mathbf{W}_1-\mathbf{W}_2\pier{\|}_{H_0^1(\mathcal{D}_a,Sym(\mathbb{R}^{3\times 3}))}^2\\
 & \notag \quad {}+  \frac{1}{2}\pier{\|}\grad \left(\boldsymbol{\Theta}_1-\boldsymbol{\Theta}_2\right)\pier{\|}^2+\pier{\|}\mathbf{Z}_1-\mathbf{Z}_2\pier{\|}_{H_{0,\diver}^1(\mathcal{D}_a,\mathbb{R}^{3\times 3})}^2+\pier{\|}\pier{\mathbf{\Phi}}_1-\pier{\mathbf{\Phi}}_2\pier{\|}_{H_0^1(\mathcal{D}_a,\mathbb{R}^{3})}^2\biggr)\\
 & \leq C\left(\frac{1}{2}\pier{\|}\mathbf{W}_1^0-\mathbf{W}_2^0\pier{\|}^2+\frac{1}{4}\pier{\|}\grad \left(\boldsymbol{\Theta}_1^0-\boldsymbol{\Theta}_2^0\right)\pier{\|}^2\right) \quad \pco{\text{for all } t\in [0,T].} 
\end{align}
\end{thm}
\begin{pf}
\pco{In view of Theorem~\ref{thm:1} and its proof, the existence result in the statement is a consequence of a limit procedure as $k\to 0$. Indeed, letting $\overline{k} $ be some fixed parameter, for $0 < k\leq \overline{k}$ we  consider 
the solution $(\mathbf{W}_k,\boldsymbol{\Theta}_k ,\pier{{}\boldsymbol{\chi}_{\alpha, k}{}}, \boldsymbol{\Sigma}_k,\pier{\mathbf{\Phi}}_k ,\mathbf{Z}_k)$ to the system \abramonew{\eqref{eqn:48symskew}--\eqref{eqn:50symskew}}
given by Theorem~\ref{thm:1}. Recalling the properties~\eqref{eqn:dpsi1}--\eqref{eqn:dpsi2} and observing that they still hold for 
$\psi_D$ and $\partial \psi_D$, it turns out that we can reproduce the estimates \eqref{eqn:54}, \eqref{eqn:542}, \eqref{eqn:55bis}, \eqref{eqn:56bis}, \eqref{eqn:57}, \eqref{eqn:58}, \eqref{eqn:comparison} uniformly with respect to $k$. Hence, we are allowed to pass to the limit in the system~\abramonew{\eqref{eqn:48symskew}--\eqref{eqn:50symskew}}, written for $\mathbf{W}_k,\boldsymbol{\Theta}_k ,\pier{{}\boldsymbol{\chi}_{\alpha, k}{}}, \boldsymbol{\Sigma}_k,\pier{\mathbf{\Phi}}_k ,\mathbf{Z}_k$,
as $k\to 0$. The argument is similar to the one developed in Section~4. Here, we deduce in particular that (see~\eqref{eqn:31}) $\boldsymbol{\Sigma} \in \curl \mathbf{Z}$,
that is $\boldsymbol{\Sigma} = \curl \mathbf{Z}$, almost everywhere, where $\boldsymbol{\Sigma}$
and $\mathbf{Z}$ denote the weak and strong limits (cf.~\eqref{eqn:convzm1}--\eqref{eqn:convsigmam1}) of some subsequence of~$\boldsymbol{\Sigma}_k$ and $\mathbf{Z}_k$, respectively.
By eliminating then the variable $\boldsymbol{\Sigma}$, we obtain the claimed existence result \pco{for} a solution of \pier{the system} \eqref{eqn:48k0}--\eqref{eqn:50k0}.}
 
We are thus left to prove the bound \eqref{eqn:contdep}, which also implies the uniqueness of the solution. 
 Let us rewrite equation \eqref{eqn:48k0}$_1$ as
 \begin{align}
  \label{eqn:cd1}
  & \left(\mathrm{e}^{-\boldsymbol{\Theta}}\curl \left(\mathbf{Z}\right),\pco{\mathbf{W}-\nwhat{\mathbf{W}}}\right)+\left(\mathbf{W}-\mathbf{I}+\overbigdot{\mathbf{W}}-\Delta \mathbf{W},\pco{\mathbf{W}-\nwhat{\mathbf{W}}}\right)+\pco{\psi}_{C_{\alpha}}(\mathbf{W})\notag\\
  & \leq \left(\mathbf{W}_{\text{ext}}(\mathbf{W},t),\pco{\mathbf{W}-\nwhat{\mathbf{W}}}\right) 
  + \pco{\psi}_{C_{\alpha}}(\nwhat{\mathbf{W}}),
 \end{align}
valid for any $\nwhat{\mathbf{W}}\in L^2(\mathcal{D}_a;Sym(\mathbb{R}^{3\times 3}))$. Taking $\nwhat{\mathbf{W}}=\mathbf{W}_2$ in the inequality~\eqref{eqn:cd1} for $\mathbf{W}_1$, $\nwhat{\mathbf{W}}=\mathbf{W}_1$ in the inequality~\eqref{eqn:cd1} for $\mathbf{W}_2$, and summing the two inequalities, we obtain that
\begin{align}
  \label{eqn:cd2}
  & \notag \left(\left(\mathrm{e}^{-\boldsymbol{\Theta}_1}-\mathrm{e}^{-\boldsymbol{\Theta}_2}\right)\curl \left(\mathbf{Z}_1\right),\pco{\mathbf{W}_1-\mathbf{W}_2}\right)+\left(\mathrm{e}^{-\boldsymbol{\Theta}_2}\curl \left(\mathbf{Z}_1-\mathbf{Z}_2\right),\pco{\mathbf{W}_1-\mathbf{W}_2}\right)\\
  & \notag +\pier{\|}\mathbf{W}_1-\mathbf{W}_2\pier{\|}_{H_0^1(\mathcal{D}_a,Sym(\mathbb{R}^{3\times 3}))}^2+\frac{1}{2}\frac{d}{dt}\pier{\|}\mathbf{W}_1-\mathbf{W}_2\pier{\|}^2 \\
  & \notag \leq\left(\mathbf{W}_{\text{ext}}(\mathbf{W}_1,t)-\mathbf{W}_{\text{ext}}(\mathbf{W}_2,t),
  \pco{\mathbf{W}_1-\mathbf{W}_2}\right)\\
&\leq L\pier{\|}\mathbf{W}_1-\mathbf{W}_2\pier{\|}_{H_0^1(\mathcal{D}_a,Sym(\mathbb{R}^{3\times 3}))}\pier{\|}\mathbf{W}_1-\mathbf{W}_2\pier{\|},
 \end{align}
where in the last inequality we have used Assumption $\mathbf{A}_3$. Moreover, taking the $L^2$ scalar product of \eqref{eqn:48k0}$_2$ for $\boldsymbol{\Theta}_1$ with $\boldsymbol{\Theta}_1-\boldsymbol{\Theta}_2$ \pco{and} the $L^2$ scalar product of \eqref{eqn:48k0}$_2$ for $\boldsymbol{\Theta}_2$ with $\boldsymbol{\Theta}_1-\boldsymbol{\Theta}_2$, \pco{then} taking the difference between the two contributions, \pco{with the help of} Assumption~$\mathbf{A}_3$ and 
\reviewa{using the Poincar\'e inequality}
, we obtain that
\abramo{
\begin{align}
 \label{eqn:cd3}
& \left(\curl\left(\mathbf{Z}_1-\mathbf{Z}_2\right)\mathbf{W}_1\mathrm{e}^{-\boldsymbol{\Theta}_1},\boldsymbol{\Theta}_1-\boldsymbol{\Theta}_2\right)+\left(\curl\left(\mathbf{Z}_2\right)(\mathbf{W}_1-\mathbf{W}_2)\mathrm{e}^{-\boldsymbol{\Theta}_1},\boldsymbol{\Theta}_1-\boldsymbol{\Theta}_2\right)\notag\\
  & \notag + \left(\curl\left(\mathbf{Z}_2\right)\mathbf{W}_2\left(\mathrm{e}^{-\boldsymbol{\Theta}_1}-\mathrm{e}^{-\boldsymbol{\Theta}_2}\right),\boldsymbol{\Theta}_1-\boldsymbol{\Theta}_2\right)
  \\
& \notag +\frac{1}{2}\pier{\|}\grad \left(\boldsymbol{\Theta}_1-\boldsymbol{\Theta}_2\right)\pier{\|}^2 + \frac{1}{4}\frac{d}{dt}\pier{\|}\pier{\grad{}}(\boldsymbol{\Theta}_1-\boldsymbol{\Theta}_2)\pier{\|}^2 \\
&\notag \leq 
\left(\boldsymbol{\Omega}_{\text{ext}}{(\boldsymbol{\Theta}_1,t)}-\boldsymbol{\Omega}_{\text{ext}}{(\boldsymbol{\Theta}_2,t)},\boldsymbol{\Theta}_1-\boldsymbol{\Theta}_2\right)\\
& \leq G\pier{\|}\boldsymbol{\Theta}_1-\boldsymbol{\Theta}_2\pier{\|}_{H_0^1(\mathcal{D}_a,Skew(\mathbb{R}^{3\times 3}))}\pier{\|}\boldsymbol{\Theta}_1-\boldsymbol{\Theta}_2\pier{\|}
\leq C\pier{\|}\grad\left(\boldsymbol{\Theta}_1-\boldsymbol{\Theta}_2\right)\pier{\|}^2.
\end{align}
}%
\pco{Next, taking} the $L^2$ scalar product of \pco{\eqref{eqn:48k0}$_4$} for $\pier{\mathbf{\Phi}}_1$ with $(\pier{\mathbf{\Phi}}_1-\pier{\mathbf{\Phi}}_2)$, the $L^2$ scalar product of \pco{\eqref{eqn:48k0}$_4$} for $\pier{\mathbf{\Phi}}_2$ with $(\pier{\mathbf{\Phi}}_1-\pier{\mathbf{\Phi}}_2)$, and \pco{subtracting} the two contributions, we \pco{arrive at}
\begin{align}
 \label{eqn:cd4}
 & \notag \pier{\|}\pier{\mathbf{\Phi}}_1-\pier{\mathbf{\Phi}}_2\pier{\|}_{H_0^1(\mathcal{D}_a,\mathbb{R}^3)}^2 \\
& =-\left(\left(\mathrm{e}^{\boldsymbol{\Theta}_1}-\mathrm{e}^{\boldsymbol{\Theta}_2}\right)\mathbf{W}_1,\pier{\grad{}} (\pier{\mathbf{\Phi}}_1-\pier{\mathbf{\Phi}}_2)\right)-\left(\mathrm{e}^{\boldsymbol{\Theta}_2}\left(\mathbf{W}_1-\mathbf{W}_2\right),\pier{\grad{}} (\pier{\mathbf{\Phi}}_1-\pier{\mathbf{\Phi}}_2)\right).
 \end{align}
 Analogously, taking the $L^2$ scalar product of \pco{\eqref{eqn:48k0}$_5$} for $\boldsymbol{Z}_1$ with $(\boldsymbol{Z}_1-\boldsymbol{Z}_2)$ \pco{and for} $\boldsymbol{Z}_2$ with $(\boldsymbol{Z}_1-\boldsymbol{Z}_2)$, \pco{then} the difference between the two contributions \pco{leads to}
\begin{align}
 \label{eqn:cd5}
 & \notag \pier{\|}\boldsymbol{Z}_1-\boldsymbol{Z}_2\pier{\|}_{H_{0,\diver}(\mathcal{D}_a,\mathbb{R}^{3\times 3})}^2\\
&=\left(\left(\mathrm{e}^{\boldsymbol{\Theta}_1}-\mathrm{e}^{\boldsymbol{\Theta}_2}\right)\mathbf{W}_1,\curl (\boldsymbol{Z}_1-\boldsymbol{Z}_2)\right)+\left(\mathrm{e}^{\boldsymbol{\Theta}_2}\left(\mathbf{W}_1-\mathbf{W}_2\right),\curl (\boldsymbol{Z}_1-\boldsymbol{Z}_2))\right).
 \end{align}
 Finally, summing the \pco{inequalities from \eqref{eqn:cd2} to \eqref{eqn:cd5}}, using the multilinear H{\"o}lder inequality, the Young inequality and the Poincar\'e inequality,
we obtain that
 \begin{align*}
&\frac{1}{2}\frac{d}{dt}\pier{\|}\mathbf{W}_1-\mathbf{W}_2\pier{\|}^2+\frac{1}{4}\frac{d}{dt}\pier{\|}\pier{\grad{}}(\boldsymbol{\Theta}_1-\boldsymbol{\Theta}_2)\pier{\|}^2+\pier{\|}\mathbf{W}_1-\mathbf{W}_2\pier{\|}_{H_0^1(\mathcal{D}_a,Sym(\mathbb{R}^{3\times 3}))}^2
\\
&{}+\frac{1}{2}\pier{\|}\grad \left(\boldsymbol{\Theta}_1-\boldsymbol{\Theta}_2\right)\pier{\|}^2+ \pier{\|}\pier{\mathbf{\Phi}}_1-\pier{\mathbf{\Phi}}_2\pier{\|}_{H_{0,\diver}^1(\mathcal{D}_a,\mathbb{R}^3)}^2 +\pier{\|}\boldsymbol{Z}_1-\boldsymbol{Z}_2\pier{\|}_{H_0^1(\mathcal{D}_a,\mathbb{R}^{3\times 3})}^2 \\
&{}\leq \frac{1}{4}\pier{\|}\mathbf{W}_1-\mathbf{W}_2\pier{\|}_{H_0^1(\mathcal{D}_a,Sym(\mathbb{R}^{3\times 3}))}^2+C\pier{\|}\mathbf{W}_1-\mathbf{W}_2\pier{\|}^2+C\pier{\|}\pier{\grad{}} (\boldsymbol{\Theta}_1-\boldsymbol{\Theta}_2)\pier{\|}^2\\
&{}+\left\|\mathrm{e}^{-\boldsymbol{\Theta}_1}-\mathrm{e}^{-\boldsymbol{\Theta}_2}\right\|_{L^6(\mathcal{D}_a,\mathbb{R}^{3\times 3})}\pier{\|}\curl \left(\mathbf{Z}_1\right)\pier{\|}_{L^3(\mathcal{D}_a,\mathbb{R}^{3\times 3\times 3})}\pier{\|}\mathbf{W}_1-\mathbf{W}_2\pier{\|}\\
&{}+  \left\|\mathrm{e}^{-\boldsymbol{\Theta}_2}\right\|_{L^{\infty}(\mathcal{D}_a,\mathbb{R}^{3\times 3})}\pier{\|}\curl \left(\mathbf{Z}_1-\mathbf{Z}_2\right)\pier{\|}\,\pier{\|}\mathbf{W}_1-\mathbf{W}_2\pier{\|}\\
&{}+ \pier{\|}\curl \left(\mathbf{Z}_1-\mathbf{Z}_2\right)\pier{\|}\,\pier{\|}\mathbf{W}_1\pier{\|}_{L^3(\mathcal{D}_a,Sym(\mathbb{R}^{3\times 3}))}\left\|\mathrm{e}^{-\boldsymbol{\Theta}_1}\right\|_{L^{\infty}(\mathcal{D}_a,\mathbb{R}^{3\times 3})}\pier{\|}\boldsymbol{\Theta}_1-\boldsymbol{\Theta}_2\pier{\|}_{L^6(\mathcal{D}_a,Skew(\mathbb{R}^{3\times 3}))}\\
&{}+ \pier{\|}\curl \left(\mathbf{Z}_2\right)\pier{\|}\,\pier{\|}\mathbf{W}_1-\mathbf{W}_2\pier{\|}_{L^3(\mathcal{D}_a,Sym(\mathbb{R}^{3\times 3}))}\left\|\mathrm{e}^{-\boldsymbol{\Theta}_1}\right\|_{L^{\infty}(\mathcal{D}_a,\mathbb{R}^{3\times 3})}\pier{\|}\boldsymbol{\Theta}_1-\boldsymbol{\Theta}_2\pier{\|}_{L^6(\mathcal{D}_a,Skew(\mathbb{R}^{3\times 3}))}\\
&{}+ \pier{\|}\curl \left(\mathbf{Z}_2\right)\pier{\|}\,\pier{\|}\mathbf{W}_2\pier{\|}_{L^6(\mathcal{D}_a,Sym(\mathbb{R}^{3\times 3}))}\left\|\mathrm{e}^{-\boldsymbol{\Theta}_1}-\mathrm{e}^{-\boldsymbol{\Theta}_2}\right\|_{L^{6}(\mathcal{D}_a,\mathbb{R}^{3\times 3})}\pier{\|}\boldsymbol{\Theta}_1-\boldsymbol{\Theta}_2\pier{\|}_{L^6(\mathcal{D}_a,Skew(\mathbb{R}^{3\times 3}))}\\
&{}+ \left\|\mathrm{e}^{\boldsymbol{\Theta}_1}-\mathrm{e}^{\boldsymbol{\Theta}_2}\right\|_{L^3(\mathcal{D}_a,\mathbb{R}^{3\times 3})}\pier{\|}\mathbf{W}_1\pier{\|}_{L^6(\mathcal{D}_a,Sym(\mathbb{R}^{3\times 3}))}
\\
&{}\times\left(\pier{\|}\pier{\mathbf{\Phi}}_1-\pier{\mathbf{\Phi}}_2\pier{\|}_{H_0^1(\mathcal{D}_a,\mathbb{R}^3)}+\pier{\|}\boldsymbol{Z}_1-\boldsymbol{Z}_2\pier{\|}_{H_{0,\diver}^1(\mathcal{D}_a,\mathbb{R}^{3\times 3})}\right)\\
&{}+ \left\|\mathrm{e}^{\boldsymbol{\Theta}_2}\right\|_{L^{\infty}(\mathcal{D}_a,\mathbb{R}^{3\times 3})}\pier{\|}\mathbf{W}_1-\mathbf{W}_2\pier{\|}
\\
&{}\times
\left(\pier{\|}\pier{\mathbf{\Phi}}_1-\pier{\mathbf{\Phi}}_2\pier{\|}_{H_0^1(\mathcal{D}_a,\mathbb{R}^3)}+\pier{\|}\boldsymbol{Z}_1-\boldsymbol{Z}_2\pier{\|}_{H_{0,\diver}^1(\mathcal{D}_a,\mathbb{R}^{3\times 3})}\right).
 \end{align*}
\pco{Hence, by employing} the Sobolev embeddings~\eqref{eqn:2}, the bound~\eqref{eqn:omegalp}, the Young inequality and the regularity results \eqref{eqn:57thmk0}--\eqref{eqn:60thmk0}, integrating moreover in time in the interval $(0,T)$, we obtain that\pco{%
\begin{align*}
&\frac{1}{2}\pier{\|}(\mathbf{W}_1-\mathbf{W}_2)\pco{(t)}\pier{\|}^2+\frac{1}{4}\pier{\|}\grad(\boldsymbol{\Theta}_1-\boldsymbol{\Theta}_2)\pco{(t)}\pier{\|}^2 +\int_0^{\pco t}\biggl(\pier{\|}\mathbf{W}_1-\mathbf{W}_2\pier{\|}_{H_0^1(\mathcal{D}_a,Sym(\mathbb{R}^{3\times 3}))}^2\\
&\quad{}+\frac{1}{2}\pier{\|}\grad \left(\boldsymbol{\Theta}_1-\boldsymbol{\Theta}_2\right)\pier{\|}^2+ \pier{\|}\pier{\mathbf{\Phi}}_1-\pier{\mathbf{\Phi}}_2\pier{\|}_{H_0^1(\mathcal{D}_a,\mathbb{R}^3)}^2 +\pier{\|}\boldsymbol{Z}_1-\boldsymbol{Z}_2\pier{\|}_{H_{0,\diver}^1(\mathcal{D}_a,\mathbb{R}^{3\times 3})}^2\biggr) \\
&{}\leq \frac{1}{2}\pier{\|}\mathbf{W}_1^0-\mathbf{W}_2^0\pier{\|}^2+\frac{1}{4}\pier{\|}\grad \left(\boldsymbol{\Theta}_1^0-\boldsymbol{\Theta}_2^0\right)\pier{\|}^2\\
&\quad{}+\int_0^{\pco{t}}\biggl(\frac{3}{4}\pier{\|}\mathbf{W}_1-\mathbf{W}_2\pier{\|}_{H_0^1(\mathcal{D}_a,Sym(\mathbb{R}^{3\times 3}))}^2+\frac{1}{2}\pier{\|}\pier{\mathbf{\Phi}}_1-\pier{\mathbf{\Phi}}_2\pier{\|}_{H_0^1(\mathcal{D}_a,\mathbb{R}^3)}^2
+\frac{1}{2}\pier{\|}\boldsymbol{Z}_1-\boldsymbol{Z}_2\pier{\|}_{H_{0,\diver}^1(\mathcal{D}_a,\mathbb{R}^{3\times 3})}^2\biggr)
\\
& \quad {} +C\int_0^{\pco{t}} \left(\pier{\|}\mathbf{W}_1-\mathbf{W}_2\pier{\|}^2+\pier{\|}\pier{\grad{}} (\boldsymbol{\Theta}_1-\boldsymbol{\Theta}_2)\pier{\|}^2\right)  \quad \pco{\text{for all } t\in [0,T]}  ,
 \end{align*}
 from which, using a Gronwall argument, we finally show} \eqref{eqn:contdep}.
\end{pf}

\newpage
\section{Conclusion}
\label{sec:conclusions}
\abramohhb{
In this work we introduced a novel model for 
 large deformations, described in terms of the stretch and the rotation tensors as independent variables. This description has a direct geometrical interpretation and the predicted quantities may be experimented. We derived the model from a generalized form of the principle of virtual power, where the virtual velocities depend on the state variables as a consequence of internal \pierhhb{kinematic} constraints associated to the compatibility condition. In our system, the compatibility of the deformation is conditionally valid depending on the magnitude of an internal \ap{force} associated to dislocations, which enters the system as a new independent variable. We assumed a quadratic expression for the free energy density of the system, depending on the stretch, the rotation and the dislocation tensors, containing first and second gradient terms. In order to enforce the positive definiteness of the stretch matrix, we also added to the free energy the indicator function of a closed and convex set whose elements are positive definite symmetric matrices with eigenvalues which are not smaller than a given positive constant at the same time. We then assumed a quadratic form also for the dissipation potential of the system, containing viscous contributions in terms of the time derivative of the stretch tensor and on the angular velocity tensor. The \michhhb{internal forces} in the system,  which are thermodynamically coupled with the virtual velocities, were then chosen in compliance with the Clausius--Duhem dissipative equality. We adopted the strategy to invert the \pierhhb{kinematic} constraints associated to the compatibility condition through Green propagators, expressing the virtual velocities associated to the deformation map and the dislocations in terms of the virtual velocities associated to the stretch matrix and to the rotation, thus reducing the set of independent virtual velocities and eliminating their internal constraints, obtaining a system of integro-differential coupled equations with inclusions. 
\newline
We then developed the analysis of the model in a simplified setting, i.e., considering the quasi-stationary version of the full system where we neglect inertia. Through a Faedo--Galerkin approximation strategy and employing the \michhhb{Moreau--Yosida} regularization of the subdifferential of multivalued functions in the free energy, we proved the existence of a global in time weak solution in three space dimensions for the system, which is \pierhhb{actually} a strong solution, \pierhhb{by} studying the limit problem as the discretization parameter and the \michhhb{Moreau--Yosida} regularization parameter tend to zero. We also proved that everywhere in space and almost everywhere in time the material is neither flattening nor crushing and that a point which is inside its domain at a certain time remains in the interior of the domain at later times. 
\newline
Finally, we considered a limit problem, letting the magnitude of the internal \ap{force} associated to dislocations tend to zero, in which case the deformation becomes incompatible and the equations take the form of a coupled system of PDEs. In the latter situations we obtained stronger analytical results, i.e., we obtained global existence, uniqueness and continuous dependence from data of the strong solution in three space dimensions. 
\newline
In a second contribution we \pierhhb{intend to} study the full model with inertia.
}
\newpage
\section*{Acknowledgments}
\pier{This research activity has been supported by the MIUR-PRIN Grant 2020F3NCPX 
``Mathematics for industry 4.0 (Math4I4)''. AA and PC also acknowledge some support from the GNAMPA (Gruppo Nazionale per l'Analisi Matematica, la Probabilit\`a e le loro Applicazioni) of INdAM (Istituto Nazionale di Alta Matematica)
through the GNAMPA project CUP E53C23001670001, and their affiliation to the GNAMPA. Moreover, PC aims to point out his collaboration,
as Research Associate, to the IMATI -- C.N.R. Pavia, Italy.}
%
\bibliographystyle{plain}
\bibliography{biblio_CHVE} 
\end{document}